\documentclass[11pt]{amsart}
\usepackage{amsmath,amsthm,upref}
\usepackage[dvips]{graphicx}
\usepackage{psfrag}
\usepackage[all]{xy}
\usepackage{a4wide}
\usepackage{epsfig,subfigure}

\numberwithin{equation}{section}

\newcommand\Su{\mathcal S}
\newcommand\Sl{\textrm{SL}_2(\mathbb R)}
\newcommand\Dq{q}
\newcommand\CYL{Cyl(\Su)}
\newcommand\C{\mathbb C}
\newcommand\R{\mathbb R}
\newcommand\Z{\mathbb Z}
\newcommand\Q{\mathcal Q}
\newcommand\HA{\mathcal H}
\newcommand\CP{{\mathbb CP}(1)}
\newcommand\Red{\textrm{Irred}_2}
\newcommand\D{d\hspace{-0.5pt}}

\DeclareMathOperator{\RE}{Re}
\DeclareMathOperator{\IM}{Im}

\newtheorem{Theorem}{Theorem}[section]
\newtheorem{Theorem2}{Theorem}[section]
\newtheorem{Corollary}[Theorem]{Corollary}
\newtheorem{Proposition}[Theorem]{Proposition}
\newtheorem{Lemma}[Theorem]{Lemma}

\newtheorem{TheoClaim}{Claim}
\newtheorem*{Conjecture}{Conjecture}
\newtheorem*{NoNumberTheorem}{Theorem}

\newtheorem*{NoNumberLemma}{Lemma}

\newtheorem*{Claim}{Claim}
\theoremstyle{remark}
\newtheorem{Remark}[Theorem]{Remark}
\theoremstyle{definition}
\newtheorem{Definition}[Theorem]{Definition}
\newtheorem{Paragraphe}{}
\newtheorem{Example}[Theorem]{Example}
\newtheorem*{Notation}{Notation}
\newtheorem*{NoNumberDefinition}{Definition}

\newtheorem*{Convention}{Convention}

\begin{document}
\title[Components of the strata of the moduli spaces]
{Connected components of the strata of the moduli spaces
of quadratic differentials}

\author{Erwan Lanneau}

\address{
Centre de Physique Th\'eorique de Marseille (CPT), UMR CNRS 6207 \newline
Universit\'e du Sud Toulon-Var and \newline
F\'ed\'eration de Recherches des Unit\'es de
Math\'ematiques de Marseille \newline
Luminy, Case 907, F-13288 Marseille Cedex 9, France
}

\email{lanneau@cpt.univ-mrs.fr}

\subjclass[2000]{Primary: 32G15. Secondary: 30F30, 57R30, 37D40}
\keywords{Quadratic differentials, moduli space, measured foliations,
Teich\"umuller geodesic flow}

\date{October 15, 2006}

\renewcommand{\abstractname}{Abstract}

\begin{abstract}
In two fundamental classical papers, Masur~\cite{Masur} and Veech~\cite{Veech:flow} have
independently proved that the Teichm\"uller geodesic flow acts
ergodically on each connected component of each stratum of the
moduli space of quadratic differentials. It is therefore
interesting to have a classification of the ergodic components.
Veech has proved that these strata are not necessarily connected.
In a recent work~\cite{Kontsevich:Zorich}, Kontsevich and Zorich
have completely classified the components in the particular case
where the quadratic differentials are given by the global square of
Abelian differentials. 

Here we are interested in the complementary case. In a previous
paper~\cite{Lanneau}, we have described some particular component,
namely the {\it hyperelliptic} connected components, and showed
that some strata are non-connected. In this paper, we give the
general classification theorem: up to four exceptional cases in
low genera, the strata of meromorphic quadratic differentials are
either connected, or have exactly two connected components. In
this last case, one component is hyperelliptic, the other not.
This result was announced in the paper~\cite{Lanneau}.

Our proof is based on a new approach of the so-called
Jenkins-Strebel differential. We will present and use the notion of
{\it generalized permutations}.
\end{abstract}

\maketitle

\renewcommand{\abstractname}{R\'esum\'e}

\begin{abstract}
Dans des travaux maintenant classiques, Masur~\cite{Masur} et
Veech~\cite{Veech:flow} ont d\'emontr\'e ind\'ependamment que le flot
g\'eod\'esique de Teich\"uller est ergodique  
sur chaque composante connexe de chaque strate de l'espace des modules des 
diff\'erentielles quadratiques. Il devient d\`es lors int\'eressant d'avoir une 
description de ces composantes ergodiques. Veech a montr\'e que ces 
strates ne sont pas n\'ecessairement connexes. Dans un article r\'ecent 
Kontsevich et Zorich~\cite{Kontsevich:Zorich} donnent une description 
compl\`ete des composantes dans le cas particulier o\`u les diff\'erentielles 
quadratiques sont donn\'ees par le carr\'e de diff\'erentielles ab\'eliennes.

Dans cet article nous consid\'erons le cas compl\'ementaire. 
Dans un  pr\'ec\'edent article~\cite{Lanneau} nous montrions que les strates ne sont 
pas forc\'ement connexes. Nous donnions une s\'erie de strates non-connexes poss\'edant 
des composantes connexes {\it hyperelliptiques}. Dans cet article, nous d\'emontrons 
le th\'eor\`eme g\'en\'eral annonc\'e dans~\cite{Lanneau} : except\'e quatre 
cas particuliers en petits genres, les strates de l'espace des modules des 
diff\'erentielles quadratiques ont au plus deux composantes connexes. Les cas 
de non-connexit\'e \'etant d\'ecrits exactement par~\cite{Lanneau} : 
une composante est hyperelliptique, l'autre non.

Notre preuve repose principalement sur une nouvelle approche des 
diff\'erentielles quadratiques de type Jenkins-Strebel \`a savoir la
notion de {\it permutations g\'en\'eralis\'ees}.
\end{abstract}

\setcounter{tocdepth}{1}

\section{Introduction}

The moduli space of compact connected Riemann surfaces $\Su$ of
genus $g$ endowed with an integrable meromorphic quadratic
differential $\Dq$ is a disjoint union $\HA_g \sqcup \Q_g$,
where the isomorphism class of $(\Su,q)$ belongs to $\HA_g$ if
and only if $\Dq$ is the (global) square of a holomorphic Abelian
differential. It can be identified with the cotangent bundle of the
moduli space $\mathcal M_g$ of compact connected smooth complex curves
(see for instance~\cite{Hubbard:Masur}). It carries a natural flow, called the
Teichm\"uller geodesic flow, see for instance~\cite{Masur,Veech:flow}. It has a
natural stratification, whose strata are denoted by $\HA(k_1,\cdots,k_n)=
\HA_g(k_1,\cdots,k_n)$ contained in $\HA_g$ and $\Q(k_1,\cdots,k_n)=\Q_g(k_1,\cdots,k_n)$
contained in $\Q_g$, where $k_1,\dots, k-n$ is the (unordered)
list of multiplicities of the zeroes and poles of the quadratic
differentials. It is well known that the flow preserves
this stratification and that each stratum carries a complex
algebraic orbifold structure of complex dimension
$2g+n-\varepsilon$ (here $\varepsilon=1$ or $-1$ depending
respectively of the strata of Abelian differential or quadratic
differentials). Masur and Smillie~\cite{Masur:Smillie:2} proved that
all of these strata (corresponding to the multiplicities satisfying the
Gauss-Bonnet condition), except four particular cases in low genera,
are non-empty.

The aim of this paper is motivated by a fundamental theorem,
independently proved by Masur and Veech~\cite{Masur},~\cite{Veech:flow},
which asserts that the Teichm\"uller geodesic flow acts ergodically on each
connected component of each stratum (with respect to a finite measure
equivalent to the Lebesgue measure).

Kontsevich and Zorich~\cite{Kontsevich:Zorich} have recently
described the set of connected components for the strata in 
$\HA_g$. In~\cite{Lanneau}, using a construction developed
in~\cite{Kontsevich:Zorich} (hyperelliptic components), we showed
that some strata in $\Q_g$ are non-connected. More precisely, we
presented three series of one discrete-parameter strata which are
non connected; those strata have a connected component consisting of
hyperelliptic curves equipped with an ``hyperelliptic
differential''. This component is called an {\it hyperelliptic
  component}. 

In this paper, we describe the set of connected components of any
stratum of $\Q_g$. The general case stabilizes at genus $5$ and
corresponds to Theorem~\ref{theo:main:1}. Theorem~\ref{theo:main:2}
gives the remaining cases.

\begin{Theorem}
\label{theo:main:1}
Let us fix $g \geq 5$. Each stratum of the moduli space $\Q_g$
having an hyperelliptic connected component has exactly two
connected components: one is hyperelliptic --- the other not; the
detailed list is given in~\cite{Lanneau}.

Any other stratum of the moduli space $\Q_g$ of quadratic
differentials is connected.
\end{Theorem}

\noindent In small genera, there are some exceptional cases
coming from the geometry of genus one and genus two surfaces
(respectively elliptic and hyperelliptic curves). There are also $4$
mysterious cases which appear.

\begin{Theorem}
\label{theo:main:2}
Let us fix $g \leq 4$. The components of the strata of the moduli
space $\Q_g$ fall in the following description:
\begin{itemize}

\item In genera $0$ and $1$, any stratum is connected.

\item In genus $2$, there are two non-connected strata. For these two,
one component is hyperelliptic, the other not. Any other stratum of
$\Q_2$ is connected.

\item In genera $3$ and $4$, each stratum with an hyperelliptic
  connected component has exactly two connected components: one is
  hyperelliptic, the other not. 

\item There are $4$ sporadic strata in genus $3$ and $4$ which are
  non-connected and which do not possess an hyperelliptic component.

\item Any other stratum of $\Q_3$ and $\Q_4$ is connected.

\end{itemize}
\end{Theorem}

\subsection{Precise formulation of the statements}

In order to establish notations and to give a precise statement, we
review basic notions concerning moduli spaces, Abelian differentials
and quadratic differentials. There is an abundant literature on this subject;
for more details and proofs see for instance~\cite{Douady:Hubbard},
\cite{Eskin:Masur:Zorich}, \cite{FaLaPo},
\cite{Hubbard:Masur} \cite{Kontsevich}, \cite{Kontsevich:Zorich},
\cite{Masur},
\cite{Strebel},\cite{Thurston},\cite{Veech:flow},\cite{Veech:moduli},
\cite{Veech:geodesic:flow},$\dots$.
For a nice survey see~\cite{Masur:Tabachnikov} or~\cite{Zorich:Survey}.

\subsubsection{Background}

For $g \geq 1$, we define the moduli space of Abelian differentials
$\HA_g$ as the moduli space of pairs $(\Su,\omega)$ where $\Su$ is a
genus $g$ (closed connected) Riemann surface and $\omega\in
\Omega(\Su)$ a non-zero holomorphic 
$1-$form defined on $\Su$. The term moduli spaces means that we
identify the points $(\Su,\omega)$ and $(\Su',\omega')$ if there
exists an analytic isomorphism $f:\Su \rightarrow \Su'$ such that 
$f^\ast \omega'=\omega$.

For $g \geq 0$, we also define the moduli space of quadratic
differentials $\Q_g$ {\it which are not the global square of
Abelian differentials} as the moduli space of pairs $(\Su,\Dq)$
where $\Su$ is a genus $g$ Riemann surface and $\Dq$ a non-zero
meromorphic quadratic differential defined on $\Su$ such that
$\Dq$ is not the global square of any Abelian differential. In
addition, we assume that $\Dq$ has at most simple poles, if any.
This last condition guaranties that the area of $\Su$ in terms of
the metric determined by $\Dq$ is finite:
$$
\int_\Su |q| < \infty.
$$
We will denote by $\HA(k_1,\dots,k_n)$ the subset of $\HA_g$
consisting of (classes of) pairs $(\Su,\omega)$ such that $\omega$ possesses
exactly $n$ zeroes on $\Su$ with multiplicities $(k_1,\dots,k_n)$. We
will 
also denote by $\Q(k_1,\dots,k_n)$ the subset of $\Q_g$ consisting
of pairs $(\Su,\Dq)$ such that $\Dq$ possesses exactly $n$
singularities on $\Su$ with multiplicities $(k_1,\dots,k_n)$, $k_i \geq -1$.

Note that the Gauss-Bonnet formula implies that the sum of the
multiplicities $\sum k_i$ equals
$2g-2$ in the case of $\HA(k_1,\dots,k_n)$ and $4g-4$ in the case of 
$\Q(k_1,\dots,k_n)$. In Section~\ref{sec:sub:tools}, we will present Thurston's
approach to these surfaces via the theory of measured foliations.

From these definitions, it is a well known part of the
Teichm\"uller theory that these spaces are (Hausdorff) complex
analytic, and in fact algebraic, spaces (see~\cite{Douady:Hubbard}
for a nice description of the stratum $\Q(1,\dots,1)$; see also
\cite{Hubbard:Masur}, \cite{Kontsevich}, \cite{Veech:moduli}).
Basically, one can see that as follows. We first concentrate on the
strata of the moduli spaces $\HA_g$. 

Let $(\Su,\omega^2)$ be a representative of an element in $\HA(k_1,\dots,k_n)$, 
$S$ its underlying topological surface, and 
$P_1,\dots, P_n$ its singular points. Let us denote by ${\rm hol}=
{\rm hol}_{(\Su,\omega)}$ the group morphism $H_1(S,\{P_1\dots,P_n\},{\mathbb
Z})\rightarrow{\mathbb C}$  defined by ${\rm
hol}([\gamma])=\int_\gamma \omega$ for every $1$-cycle $\gamma$ in $S$
relative to $\{P_1\dots,P_n\}$. Fix a basis $(\gamma_1,\dots,
\gamma_{2g+n-1})$ of the free abelian group $H_1(S,\{P_1\dots,P_n\},
\mathbb Z)$. Any other element of $\HA(k_1,\dots,k_n)$ will be
represented by an element having the same underlying surface and the
same singular points. With these notations, the map
$$
\Phi = \left( \begin{array}{ccc} \HA(k_1, k_2,\dots, k_n) &
  \longrightarrow & H^1(S,\{P_1,\dots,P_n\},\mathbb C)  \\  
\Su' &  \longmapsto & (\gamma_1,\dots,\gamma_{2g+n-1})\mapsto
(hol_{\Su'}(\gamma_1) , \dots ,
hol_{\Su'}(\gamma_{2g+n-1})) \end{array} \right) 
$$
is named the {\em period map} and is a local homeomorphism in a
neighbourhood of $(\Su,\omega^2)$. 
Therefore we get a locally one-to-one correspondence
between the corresponding stratum of $\HA_g$ and an open domain
in the vector space $H^1(S,\{P_1,\dots,P_n\};\C) \simeq \mathbb
C^{2g+n-1}$. The change of coordinates are affine maps outside the
singularities of $\HA(k_1, k_2,\dots, k_n)$ and produces after 
a study of the singularities a differentiable orbifold structure on
the strata of $\HA_g$.

Let us now consider the case of a stratum of the moduli space
$\Q_g$. For every $(\Su,\Dq) \in \Q_g$, 
consider the canonical double cover $\pi : \hat \Su
\rightarrow \Su$ such that $\pi^\ast\Dq = \omega^2$ for some holomorphic 
Abelian differential $\omega$ on $\hat \Su$ (see for
instance~\cite{Lanneau:spin}). As above, we consider the period
map between a neighborhood of the point $(\hat \Su,\omega)$
and an open subset of $H^1(\hat S,\{\hat P_1,\dots,\hat P_n\};\C)$. The
covering involution $\tau : \hat \Su \rightarrow \hat \Su$ 
induces an involutive linear map on this cohomology vector space. 
Therefore, this vector space decomposes into two eigenspaces for 
$\tau^\ast$, say $E_{-1}$ and
$E_{+1}$, with eigenvectors $-1$ and $+1$. Abelian differentials
in $E_{-1}$ are precisely those which arise from quadratic
differentials on $\Su$ by pull-back by $\pi$. 
Hence we obtain a one-to-one correspondence between a neighborhood of any point in 
the corresponding strata of $\Q_g$ and an open domain of $E_{-1}
\simeq \mathbb C^{2g+n-2}$.

Next we recall the construction of a measure $\mu$ on each
stratum. For that, the tangent space to $\HA_g$ (respectively $\Q_g$)
at each point contains a lattice:
$$
H^1(\Su,\{P_1,\dots,P_n\};\Z)\oplus i \cdot H^1(\Su,\{P_1,\dots,P_n\};\Z) 
\subset H^1(\Su,\{P_1,\dots,P_n\};\C).
$$
We define the function $A:\HA_g \rightarrow \mathbb R_{+}$ by the
formula $A(\Su,\omega)= \cfrac{i}{2} \int_\Su
\omega\wedge\overline{\omega}$. This is the area
of $\Su$ in terms of the flat metric associated to $\omega$.

The group $\Sl$ acts by linear transformations with
constant coefficients on the pair of real-valued $1$-forms
$\left(\RE(\omega),\IM(\omega)\right)$. In the local affine
coordinates, this action  is the action of $\Sl$ on the coefficient 
field of the cohomology vector space $H^1(S,\{P_1,\dots,P_n\};\C)$. 
From this description, it is clear that the subgroup
$\Sl$ preserves the measure $\mu$ and the function $A$.

On the hypersurface $\HA_g^{(1)}=A^{-1}(1)$, we define the induced
measure by the formula
$$
\mu^{(1)}=\cfrac{\mu}{dA}.
$$
Recall that each stratum carries a complex algebraic orbifold structure
modeled on the first relative cohomology group (see for 
instance~\cite{Kontsevich}). The dimensions are respectively given by
$$
\begin{array}{l}
\dim_\C \HA(k_1,\dots,k_n) = 2g+n-1 \qquad \textrm{where} \qquad
k_1+\dots+k_n=2g-2 \\
\dim_\C \Q(k_1,\dots,k_n) = 2g+n-2 \qquad \textrm{where} \qquad
k_1+\dots+k_n=4g-4.
\end{array}
$$
The action of the 1-parameter subgroup of diagonal matrices
$g_t:=\textrm{diag}(e^{t/2},e^{-t/2})$ presents a particular
interest for our purpose. It gives a measure-preserving flow for
$\mu^{(1)}$, preserving each stratum. This flow is known as the
{\it Teichm\"uller geodesic flow}. Note that orbits under $g_t$
project to Teichm\"uller geodesics on the moduli space of Riemann
surfaces $\mathcal M_g$. The next fundamental result motives our study
(see~\cite{Masur}, \cite{Veech:flow}).

\begin{NoNumberTheorem}[Masur,~Veech]
The Teichm\"uller geodesic flow acts ergodically on each connected
component of each stratum of $\HA^{(1)}_g$ with respect to the 
measure $\mu^{(1)}$, which is finite and in the Lebesgue class.
\end{NoNumberTheorem}

A direct corollary of the finiteness of the
measure $\mu^{(1)}$ on any stratum of $\HA^{(1)}_g$ is a proof
of the conjecture of Keane~\cite{Keane}: almost all
intervals exchange transformations are uniquely ergodic.

\subsubsection{Topology of the Moduli Space}

Following the theorem of Masur and Veech, we are interested in
the classification of the connected components of the strata of
$\HA_g\sqcup \Q_g$. Veech and Arnoux discovered, by direct
calculations in terms of {\it Rauzy classes}, that some strata are
non-connected. They have proved that $\HA(2)$ is connected,
$\HA(4)$ have $2$ connected components and $\HA(6)$ have $3$ connected
components.

Recently, in the context of moduli space of Abelian differential
$\HA_g$, Kontsevich and Zorich (see~\cite{Kontsevich:Zorich}) obtained
the following complete description.

\begin{NoNumberTheorem}[Kontsevich,~Zorich]
Let $g \geq 4$ be any integer. The topology of any stratum of $\HA_g$
is given by the following list:
\begin{itemize}

\item The strata $\HA(2g-2)$ and $\HA(2k,2k)$, for any $k \geq 2$,
have three connected components.

\item Any other stratum $\HA(2k_1,\dots,2k_n)$, for any
$k_i \geq 1$, has two connected components.

\item The stratum $\HA(2k-1,2k-1)$, for any $k \geq 2$, has two
connected components.

\end{itemize}
Any other stratum of Abelian differentials on a surface of genus $g
\geq 4$ is non-empty and connected.
\end{NoNumberTheorem}

The description of connected components for strata of genera $1 \leq g
\leq 3$ is similar to the previous one with some exceptions; 
we do not present the result here 
(see~\cite{Kontsevich:Zorich}). Roughly speaking,
Kontsevich and Zorich use two invariants to obtain
this classification: the parity of the spin structure and the
hyperellipticity.

In~\cite{Lanneau:spin}, we prove that the first invariant extends
trivially to the moduli space $\Q_g$. However, the second produces
non-trivial values. In~\cite{Lanneau}, we classify all strata for
which this second one invariant produces non-trivial values. In order
to present our statement, we will recall briefly this construction in
the coming Section.

In this paper, we will show that this (hyperelliptic)
invariant is complete in genera $g\geq 5$: it classifies precisely
the components of the strata of $\Q_g$. For small genera, we
obtain a similar result with $4$ additional mysterious components.

\begin{Remark}
As a direct corollary of our result and Kontsevich-Zorich's
theorem, we draw that the ergodic components of the Teichm\"uller 
geodesic flow are given by an explicit list. In particular, if 
$n\geq 5$, for any $g\geq 0$, 
any stratum of $\HA \sqcup \Q_g$ with $n$ singularities is connected.
\end{Remark}

\begin{Remark}
This paper achieves the classification of connected components of
the strata of the moduli spaces $\HA \sqcup \Q_g$ announced in the vast program
in~\cite{Kontsevich}. Nevertheless a more precise description of the
strata is not completely understood, even for the simplest
``non-trivial'' case: $\HA(2)$.
\end{Remark}

\begin{Conjecture}[Kontsevich]
Each connected component of the strata of the moduli space $\HA \sqcup \Q_g$ is a 
$K(\pi,1)$, where $\pi$ is a group commensurable with some mapping
class group.
\end{Conjecture}

\subsubsection{Hyperelliptic components}

We will need the following statement. A proof can be found
in~\cite{Kontsevich:Zorich}; see also~\cite{Lanneau}.

\begin{Proposition}
\label{prop:particular:cases}
Any stratum of the moduli space $\Q(k_1,\dots,k_n)$ in genus $0$ is
connected.
\end{Proposition}

Let $\Su_g$ be a (compact, connected) Riemann surface endowed with a 
(integrable, meromorphic) quadratic differential $\Dq_0$ 
which is not the global square of any Abelian
differential. Let $(k_1, \dots, k_n)$ be its singularity pattern.
We do not exclude the case when some of $k_i$ are equal to zero:
by convention this means that we have some marked points. Sometimes 
We shall use the ``exponential'' notation to denote multiple singularities 
(simple poles), for example $\Q(-1^5,1):=\Q(-1,-1,-1,-1,-1,1)$.

Let $\pi : \tilde{\Su}_{\tilde{g}} \rightarrow \Su_g$ be a
(ramified) covering such that the image of any ramification point
of $\pi$ is a marked point, or a zero, or a pole of the quadratic
differential $\Dq_0$. Fix the {\it combinatorial type} of the
covering $\pi$: the degree of the covering, the number of critical
fibers and the ramification index of the points in every critical
fiber. Consider the induced quadratic differential
$\pi^{\ast}\Dq_0$ on $\tilde{\Su}_{\tilde{g}}$; let
$(\tilde{k}_1,\dots,\tilde{k}_m)$ be its singularity pattern.

Deforming slightly the initial point
$(\Su_g,\Dq_0)\in\Q(k_1,\dots,k_n)$, we get a ramified covering
over the deformed Riemann surface with the same combinatorial type
as the covering $\pi$. 
The induced quadratic differential $\pi^{\ast}\Dq$ has the same
singularity pattern $(\tilde{k}_1,\dots,\tilde{k}_m)$ as
$\pi^{\ast}\Dq_0$. Thus we obtain a map from the stratum 
$\Q(k_1, \dots, k_n)$ to the stratum $\Q(\tilde{k}_1,\dots,\tilde{k}_m)$. 
We will denote denote this map by $\Pi$.
In~\cite{Lanneau}, we prove that $\Pi$ is an immersion.

Now we recall an example of this 
construction for strata having four singularities.
Consider a meromorphic quadratic differential $\Dq$ on $\CP$
having the singularity pattern $(2(g-k)-3,2k+1,-1^{2g+2})$, where
$k \geq -1$, $g \geq 1$ and $g-k \geq 2$. Consider a ramified
double covering $\pi$ over $\CP$ having ramification points over
the $2g+2$ poles of $\Dq$, and no other ramification points. We
obtain a hyperelliptic Riemann surface $\tilde{\Su}$ of genus $g$
with a quadratic differential $\pi^{*}\Dq$ on it. A straightforward
computation shows that the induced quadratic differential $\pi^{*}\Dq$ has the
singularity pattern $(2(g-k)-3,2(g-k)-3,2k+1,2k+1)$. Thus we get a map
$$
\Q(2(g-k)-3,2k+1,-1^{2g+2}) \to \Q(2(g-k)-3,2(g-k)-3,2k+1,2k+1),
$$
where $k \geq -1, g \geq 1$ and $g-k \geq 2$. Computing the
dimensions of the strata, we get
\begin{align*}
&\dim_{\C}\Q(2(g-k)-3,2k+1,-1^{2g+2})=2\cdot 0 + (2g + 4) - 2=2g+2\\
&\dim_{\C}\Q(2(g-k)-3,2(g-k)-3,2k+1,2k+1)=2 g + 4 -2 =2g+2
\end{align*}
Since the dimension of the strata coincide, and the mapping is an
immersion, we obtain an open set inside the stratum
$\Q(2(g-k)-3,2(g-k)-3,2k+1,2k+1)$. This remark, together with 
Proposition~\ref{prop:particular:cases} motivate the following definition.
We call the {\it hyperelliptic} connected component 
$\Q^{hyp}(2(g-k)-3,2(g-k)-3,2k+1,2k+1)$ of the stratum 
$\Q^{hyp}(2(g-k)-3,2(g-k)-3,2k+1,2k+1)$, 
consisting of quadratic differentials 
on hyperelliptic curves of genus $g$ such that the quotient produces 
a quadratic differential with singularities data $\Q(2(g-k)-3,2k+1,-1^{2g+2})$.

In~\cite{Lanneau}, we classify all components of this type. This 
motivates the following definition.

\begin{NoNumberDefinition}
By {\it hyperelliptic components} we call the following connected 
components of the following strata of quadratic differentials.
\begin{enumerate}
\item 
\label{hyp} 
The connected component $\Q^{hyp}(2(g-k)-3,2(g-k)-3,2k+1,2k+1)$ of the 
stratum $\Q(2(g-k)-3,2(g-k)-3,2k+1,2k+1)$ consisting of quadratic differentials 
of the stratum $\Q(2(g-k)-3,2k+1,-1^{2g+2})$ where $k \geq -1,
\ g \geq 1, \ g-k \geq 2$. \newline 
The corresponding ramified double covering has ramification points over
the $2g+2$ poles.
\item The connected component $\Q^{hyp}(2(g-k)-3,2(g-k)-3,4k+2)$
of the stratum $\Q(2(g-k)-3,2(g-k)-3,4k+2)$ consisting of quadratic differentials 
of the stratum $\Q(2(g-k)-3,2k,-1^{2g+1})$ where $k \geq 0$, $g
\geq 1$ and $g-k \geq 1$.\newline
The corresponding ramified double covering has ramification points over
the $2g+1$ poles and over the zero of degree $2k$.
\item The connected component $\Q^{hyp}(4(g-k)-6,4k+2)$ of the stratum 
$\Q(4(g-k)-6,4k+2)$ consisting of quadratic differentials 
of the stratum $\Q(2(g-k)-4,2k,-1^{2g})$ where $k \geq 0$, $g \geq 2$ and
$g-k \geq 2$. \newline 
The corresponding double ramified covering has ramification points over all
the singularities.
\end{enumerate}
\end{NoNumberDefinition}

\begin{Remark}
Hyperelliptic connected components of type~(\ref{hyp}) were first
discovered by Kontsevich.
\end{Remark}

\subsubsection{Main result}

We are finally in a position to give a precise statement of our
result.

\begin{Theorem2}
Let $g \geq 5$ be any integer. Let us consider the following families $\mathcal
F_i$, $i=2,3,4$, of strata inside the moduli space of quadratic differentials $\Q_g$.
$$
\begin{array}{llll}
\mathcal F_2 = \{ & \Q(4(g-k)-6\ ,\ 4k+2)  & | & \ \ 0 \leq k \leq g-2 \}  \\
\mathcal F_3 = \{ & \Q(2(g-k)-3\ ,\ 2(g-k)-3\ ,\ 4k+2) & | & \ \ 0 \leq k \leq g-1 \}  \\
\mathcal F_4 = \{ & \Q(2(g-k)-3\ ,\ 2(g-k)-3\ ,\ 2k+1\ ,\ 2k+1) &
|  & -1 \leq k  \leq g-2 \}
\end{array}
$$
Then any stratum listed above has exactly two connected
components: one is hyperelliptic --- the other not.

Any other stratum of the moduli space $\Q_g$ is non-empty and
connected.
\end{Theorem2}

\noindent In small genera, some components are missing compared
to the general case. The complete description is given by the 
following theorem.

\begin{Theorem2}
Let $g \leq 4$ be any non-negative integer. The components of the strata of the moduli
space $\Q_g$ fall into the following description.
\begin{itemize}

\item In genera $0$ and $1$, any stratum is non-empty and
connected, except $\Q(0)$ and $\Q(-1,1)$ which are empty.

\item In genus $2$, there are two non-connected strata:
$\Q(-1,-1,6)$ and $\Q(-1,-1,3,3)$; they have two
components, one component is hyperelliptic the other not. Any other
stratum of $\Q_2$ is connected.

\item In genera $3$ and $4$, each stratum possessing an hyperelliptic
connected component has exactly two connected components: one is
hyperelliptic, the other not. 

\item There are $4$ sporadic strata:
$$
\Q_3(-1,9), \ \Q_3(-1,3,6), \ \Q_3(-1,3,3,3), \ \Q_4(12)
$$
which have exactly two connected components.

\item Any other stratum of $\Q_3$ and $\Q_4$ is connected.

\end{itemize}
\end{Theorem2}

\begin{Remark}
\label{rk:zorich:special}
The fact that the four sporadic components are not connected is due to 
Zorich (see~\cite{Zorich:computation}). Zorich proved this result by 
a direct calculation in terms of {\it extended Rauzy classes}.
It would be interesting to have an
algebraic-geometric proof of the non-connectedness of these strata.
\end{Remark}

\begin{Remark}
Using our approach of quadratic differentials, we get a new proof of a
result of Masur and Smillie~\cite{Masur:Smillie:2} concerning the fact
that the strata $\Q(0)$, $\Q(-1,1)$, $\Q(1,3)$ and $\Q(4)$ are empty.
\end{Remark}

\subsection{Outline of the proof}

The proof involves the dynamic and topology of
measured foliations. We will use the well known fact that
quadratic differentials $\Dq$ on $\Su$ and pair of transverse
measured foliations on $\Su$ define the same objects
(see~\cite{Hubbard:Masur}). 
We will say that a component $\mathcal C_1$ of a strata of $\Q_g$ 
is {\it adjacent}
to a component $\mathcal C_2$ if $\mathcal C_2 \subset
\overline{\mathcal C_1}$ where the closure is taken inside the
whole space $\Q_g$.

\begin{Paragraphe}[Claim A]
Let $\mathcal C$ be any component of any stratum $\Q_g(k_1,\dots,k_n)$ 
with $n\geq 2$ and $g\geq 3$. We assume that $\mathcal C$ is neither a
hyperelliptic component nor a sporadic component
$\Q^{irr}(-1,9)$, $\Q^{irr}(-1,3,6)$ or $\Q^{irr}(-1,3,3,3)$ (to 
be defined later on).
Then $\mathcal C$ is adjacent to a component of the stratum
$\Q(4g-4)$.
\end{Paragraphe}

\begin{Paragraphe}
We will use the following corollary of a theorem of
Konsevich~(\cite{Kontsevich:Zorich}): The number of connected
components of a stratum of $\Q_g$, which are adjacent to a component
of the stratum $\Q(4g-4)$, is bounded by the number of components of
$\Q(4g-4)$. In addition, for a given stratum, there is at most one hyperelliptic or one 
sporadic component. \\
These two remarks, together with Claim A, provide 
us with a bound on the number of 
components of any stratum in terms of the number $r_g$ of components
of $\Q(4g-4)$. More precisely, 
the number of components of a statum containing a 
hyperelliptic component (or a sporadic component) 
is bounded upper by $r_g + 1$. The number of components of a statum 
containing nor a hyperelliptic component neither a sporadic component 
is bounded upper by~$r_g$. 
\end{Paragraphe}

\begin{Paragraphe}
The list of stratum which possess a hyperelliptic component
has been given in~\cite{Lanneau}. In particular, these strata are 
non connected for $g\geq 3$. 
Hence the proof of Theorem~\ref{theo:main:1} is reduced to the proof 
that the equality $r_g=1$ holds for any $g \geq 5$. \\ 
The genera $3$ and $4$ are similar with some additional cases
(sporadic strata). Theorem~\ref{theo:main:2} follows from the equality 
$r_3=1$ and the inequality $r_4 \leq 2$ that we prove in Section~\ref{sec:minimal}.
The small genera cases $g=1,2$ are considered separately; this
is done in Section~\ref{sec:proof}.
\end{Paragraphe}

Therefore the proof of our result is reduced to the one of Claim~A and to the
computation of $r_g$. \medskip

\noindent 
We will say that a saddle connection $\gamma$ on $(\Su,\Dq)$ has {\it multiplicity one} 
(see~\cite{Masur:Zorich}) if one can collapse $\gamma$ to a point to
get a new {\it closed} half-translation surface $\Su'$. 
In particular, this condition is satisfied if $\gamma$ is ``small''
with respect to the other saddle connections in the direction of
$\gamma$. More precisely if
\begin{equation}
\label{eq:reader}
|\gamma| < \cfrac{1}{2}\ |\eta|, \textrm{ for any saddle connection $\eta$ on
} \Su \backslash \gamma \textrm{ in the direction } \overrightarrow{\gamma}
\end{equation}
then $\gamma$ has multiplicity one. We will give a useful criterion on
half-translation surfaces to have such an inequality.

\begin{Paragraphe}
\label{para:adja}
Let $(\Su,\Dq) \in \mathcal C$ and $\gamma$ be a 
multiplicity one saddle connection of $(\Su,\Dq)$. 
Then $\gamma$ can be collapsed 
to a point to get a new half-translation surface $(\Su',\Dq')$. 
The component $\mathcal C'$ containing $(\Su',\Dq')$ is contained 
in $\overline{\mathcal C}$.
\end{Paragraphe}

\begin{Paragraphe}
Let us assume that $\Su$ is decomposed into a single cylinder for the
horizontal flow. It means that any two regular horizontal
geodesics are closed and homologous: they form a family which fill the
surface into a metric cylinder. The boundaries of this cylinder
consist of a set of saddle connections and separatrix loops. The
arrangement of these separatrices is described by a ``generalized
permutation'' $\pi$ (see
Section~\ref{sec:generalized:permutations}). We will note
$(\Su,\Dq)=\Su(\pi)$ for a one cylinder surface.

In Section~\ref{sec:properties:permutations} we give a {\it
    combinatorial criterion} (namely irreducibility) on the
    combinatorics of $\pi$ so that there exists a 
saddle connection satisfying Equation~(\ref{eq:reader}); In particular
this produces a multiplicity one saddle connection.
\end{Paragraphe}

\begin{Paragraphe}
\label{para:excep}
Let $\Su(\pi)\in \mathcal C$ be a point in an arbitrary component
(such points are dense in each stratum of the moduli space
$\Q_g$). We consider the set of surfaces $\Su(\pi_s):=h_s \cdot
\Su(\pi)$ where $h_s = \left( \begin{smallmatrix} 1 & s \\ 0 & 1 
\end{smallmatrix} \right)$.

If $\Su(\pi_s)$ has a multiplicity one saddle connection for any $s$,
we are done thanks to {\bf \ref{para:adja}}. Otherwise the generalized 
permutation $\pi_s$ must obey to some combinatorial conditions. Repeating
that for any $s$ we get that either $\Su(\pi)$ has a multiplicity one 
saddle connection or $\pi$ is an exceptional (or hyperelliptic)
permutation. The last case implies that $\mathcal
C$ is an exceptional (or hyperelliptic) component 
contradicting the assumptions of Claim~A.
\end{Paragraphe}

\begin{Paragraphe}
\label{para:minimal}
The proof of the connectedness of $\Q(4g-4)$ (i.e. of the equality $r_g=1$) for $g\geq 
5$ is done by induction on the genus $g$ of the surfaces. We first show
by a direct computation that the minimal 
stratum in genus $5$ is connected. We then use the surgery
``Bubbling a handle'' for the step of the induction. Precisely, we 
find in each component a surface with a cylinder filled by closed
geodesic such that the boundary is formed by two single
multiplicity one saddle connections. Therefore one can ``erase''
this cylinder to obtain a closed $(g-1)$ half-translation surface.
\end{Paragraphe}

The paper is organized in the following manner. 
In Section~\ref{sec:tools} we remind key results concerning the
geometry of quadratic differentials. Then
Sections~\ref{sec:generalized:permutations}--\ref{sec:properties:permutations}
are devoted to the notion of {\it generalized permutations}. We
develop this useful notion and give some relations between the
combinatorics of $\pi$ and the dynamic of the measured foliation on
$\Su(\pi)$. Finally in Section~\ref{sec:minimal} and
Section~\ref{sec:adjacency} we prove respectively the two
points {\bf \ref{para:minimal}} and {\bf \ref{para:excep}}.

\subsection{Remarks and Acknowledgments}

This paper is the last part of my PhD, the first two parts being
the papers~\cite{Lanneau} and~\cite{Lanneau:spin}. I would like to thank 
Pascal~Hubert and Anton~Zorich for encouraging me to write 
this paper.

I also thank Pascal~Hubert, Maxim~Kontsevich, Howard~Masur, 
Anton~Zorich and the anonymous referee for helpful comments on
preliminary versions of this manuscript.

I thank the Institut de Math\'ematiques in Luminy, the
Max-Planck-Institute for Mathematics in Bonn and the Centre de
Physique Th\'eorique in Marseille for excellent welcome and working
conditions.

\section{Preliminaries and preparatory material}
\label{sec:tools}

In order to establish notations and preparatory material, we
review basic notions concerning quadratic differentials versus
half-translation surfaces. These surfaces have been considered and
studied by numerous authors in various guises, see
say~\cite{Hubbard:Masur}, \cite{Kerckhoff:Masur:Smillie},
\cite{Strebel} for more details; See 
also~\cite{Eskin:Masur:Zorich} and~\cite{Masur:Zorich} for recent
related developments on surgeries about half-translation surfaces.

\subsection{Flat metrics}

\subsubsection{Flat surfaces and geodesics}
\label{sec:sub:tools}

A {\it half-translation surface} is a (compact, connexe, real) genus
$g$ surface equipped with a flat metric (with isolated conical
singularities) such that the holonomy group belongs to $\{\pm
\textrm{Id}\}$. Here holonomy means that the parallel transport of a
vector along a small loop going around a conical point brings the
vector back to itself or to its negative. This implies that all cone
angles are integer multiples of $\pi$. An equivalent definition is the
following. A half-translation surface is a triple $(\Su,\mathcal
U,\Sigma)$ such that $\Su$ is a topological compact connected surface,
$\Sigma$ is a finite subset of $\Su$ (whose elements are called
{\em singularities}) and $\mathcal U = (U_i,z_i)$ is an
atlas of $\Su \setminus \Sigma$ such that the transition maps
$z_j \circ z_i^{-1} : z_i(U_i\cap U_j) \rightarrow
z_j(U_i\cap U_j)$ are translations or half-translation: $z_i = \pm
z_j + const$. This implies that the holonomy belongs to $\{\pm
\textrm{Id}\}$. \medskip

Therefore, we get on $\Su$ a flat metric with conical singularities
located in $\Sigma$ (possibly not all $\Sigma$). We also get a {\it
quadratic differential} defined locally in the coordinates $z_i$ by
the formula 
$\Dq=\D  z_i^2$. This form extends to the points of $\Sigma$ to zeroes, 
poles or marked points (see~\cite{Masur:Tabachnikov}). We will
sometimes use the notation $(\Su, \Dq)$ or simply $\Su$. \medskip

\begin{Remark}
The holonomy is trivial if and only if all transition functions are
translations or equivalently the quadratic differentials $\Dq$ is the
square of an Abelian differential. We will then say that $\Su$ is a
translation surface.
\end{Remark}

Hence the half-translation structure defines on $\Su \setminus \Sigma$
a Riemannian structure; we therefore have notions of directional
foliation, geodesic, length, angle, measure, etc.

\begin{Convention} \textrm{ }
\begin{itemize}

\item the {\it singularities} are the zeroes and poles of $\Dq$.

\item leaves of the directional foliation are called geodesics.

\item geodesics meeting singularities are called {\it separatrices}. A
  geodesic emanating from a singularity 
and going back to the same is called a {\it separatrix loop}. A
geodesic connecting two different singularities is called a {\it
saddle connection}.

\item a geodesic not passing through a singularity is called {\it
     regular}.  
\end{itemize}
\end{Convention}

One can see the equivalence between half-translation surface and
quadratic differential as follows. We start with the first
definition. If we cut a flat surface $\Su$ successively along an
appropriate collection of saddle connections, we can decompose it
into polygons contained in $\mathbb C$. We may then view $\Su$ as
a union of polygons with sides ordered by pairs consisting of 
parallels sides of the same length. The surface $\Su$ is then isometric
to the polygons where we identify these pair by translations or
half-translations (translation post-composed with a central symmetry). Note that we
have endowed each polygon with a complex coordinate. By
construction, the transition functions in these complex
coordinates $z$ have the form
$$
z \mapsto \pm z + const.
$$
Thus any flat surface with the conical singularities removed is
endowed  with a natural complex structure. Moreover, consider a
holomorphic quadratic differential $\Dq=dz^2$ on every polygon.
Since $dz^2=d(\pm z+const)^2$ we  obtain a globally well defined
holomorphic quadratic differential on $\Su'$. It is a direct
calculation to check that the complex structure and the quadratic
differential can be extended to the singularities; the quadratic
differential $\Dq$ extends to a (possibly meromorphic) form on
$\Su$ with zeroes or simple poles at every conical point. Note
that when all transitions function are only translations, the
quadratic differential $\Dq$ can be globally written as
$\Dq=\omega^2$, where $\omega$ is an Abelian differential. In this
case the corresponding foliation is oriented.

Conversely, given a pair $(\Su,\Dq)$, and a point $P \in \Su$ such
that $\Dq(P) \not = 0$, the integral $\int_{z_0}^z \sqrt{\Dq}$
produces a local coordinate $z$ near $P$ such that $\Dq = dz^2$.
Thus $|dz|^2$ defines a flat metric on $\Su$.  At a
singularity of multiplicity $k\geq -1$ the total angle we get is
$(k+2)\pi$. Remark that for regular point of $\Dq$ ($k=0$), one
get regular point of the metric.

\subsubsection{$\Sl$-action}

Given any matrix $A\in \Sl$, we can post-compose the local
coordinate of the charts of our translation atlas on $(\Su,\Dq)$
by $A$. One easily checks that this gives a new half-translation
surface, denoted by $A\cdot (\Su,\Dq)$. In local coordinates, this
gives
$$
A\cdot (\Su,\D z^2) = (\Su,(\D Az)^2).
$$
Hence this produces an $\Sl$-action on these half-translation
surfaces.

\subsubsection{Cylinders}

Closed regular geodesics appear in families of parallel
geodesics of the same length. Such parallel closed geodesic,
typically, do not filled the surface, but only a cylindrical
subset. Each boundary component of such a cylinder is comprised by 
saddle connections.

Generically, each boundary component of a cylinder filled with
closed regular geodesics is a single closed saddle connection. The
converse, however, is false.  A closed saddle connection does not
necessarily bound a cylinder of regular closed geodesics. In fact,
it bounds such a cylinder if and only if the angle at the
singularities between the outgoing and incoming segments is
exactly $\pi$. One calls such a cylinder a {\it simple} cylinder (see
Section~\ref{sec:minimal}). See also Figure~\ref{fig:simple:cyl}. 

\subsubsection{Adjacency}

Let $\mathcal C,\ \mathcal C'$ be two connected components of the
moduli space $\Q_g$. We will say that $\mathcal C'$ is adjacent to
$\mathcal C$ if $\mathcal C \subset \overline{\mathcal C'}$ (the
closure being taken inside the whole space $\Q_g$).

\subsubsection{Example}

We end this section with an example of half-translation surface
with Figure~\ref{fig:suspension:gp}. Identifying pairs of sides of
the polygon by isometries, we get a half-translation surface of
genus $g=1$. Note the form $\D z$ is not globally defined but $\D
z^2$ is; therefore the holonomy is exactly $\{\pm \textrm{Id}\}$.

\subsection{\^Homologous saddle connections}
\label{sec:homologous}

Let $\Su$ be a half-translation surface. We denote by $\pi : \hat
\Su \rightarrow \Su$ the standard orientating double covering so
that $\pi^\ast \Dq = \hat \omega^2$ (see~\cite{Lanneau:spin}). Let
$\tau$ be the induced involution of the covering. Let $\gamma$ be
a {\it compact separatrix} on $\Su$. We consider $\gamma^{+}$ and
$\gamma^{-}$ the two lifts of $\gamma$ by $\pi$. We choose an
orientation of $\gamma$. According to this choice, we define
$$
\hat \gamma=[\gamma^{+}]-[\gamma^{-}]
$$
($\hat \gamma$ is well defined up to a sign). If $P_i$ denote the
singularities of $\Dq$ on $\Su$ and $\hat P_i$ the singularities of
$\omega$ on $\hat \Su$, let $H^{+}_1(\hat \Su,\{\hat P_i\},\mathbb C)$ be
the first homological group invariant with respect to the involution $\tau$ and $H^{-}_1(\hat
\Su,\{\hat P_i\},\mathbb C)$ the first homological group
anti-invariant with respect to the involution $\tau$. We therefore get
$$
\hat \gamma \in H^{-}_1 / \pm
$$

\subsubsection{\^Homologous saddle connections}

Following Masur and Zorich~\cite{Masur:Zorich}, we say that two
compact separatrices $\gamma$ and $\eta$ are {\it \^homologous} if
their corresponding loops $\hat \gamma$ and $\hat \eta$ on $\hat
\Su$ are proportional inside $H^{-}_1$.

\begin{figure}[htbp]
\begin{center}
\psfrag{pi}{$\pi=\left(\begin{array}{ccc} 1 & 1 & 2 \\ 3 & 2 & 3
\end{array} \right)$} \psfrag{1}{$1$}  \psfrag{2}{$2$}
\psfrag{3}{$3$} \psfrag{S}{$\Su(\pi,\lambda)$}
\psfrag{gp}{$\gamma(\pi)$} \psfrag{e}{$\eta$}

\includegraphics[width=14cm]{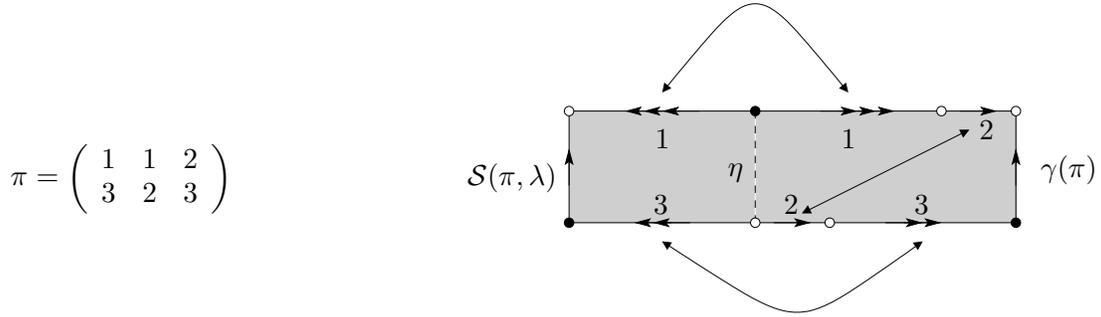}
\end{center}
\caption{ \label{fig:suspension:gp} A suspension over a
permutation $\pi$ with an appropriate admissible vector $\lambda$.
The resulting point $(\Su(\pi,\lambda),\D z^2)$ belongs to the
stratum $\Q(-1,-1,2)$. The black bullets correspond to the poles
and the white bullets to the unique zero of the differential $\Dq$
on $\Su(\pi,\lambda)$.
}
\end{figure}

Note that this definition does not depend of the choice of the
orientation of the geodesics neither the choice of a direction on
$\Su$. Moreover, $\gamma$ and $\eta$ are not supposed to be
homeomorphic to a circle. For instance, one can have a saddle
connection (homeomorphic to a segment) \^homologous to a
separatrix loop (homeomorphic to a circle); see
Figure~\ref{fig:suspension:gp}. The following proposition gives a
necessary condition for two compact separatrices to be \^homologous.

\begin{Proposition}[Masur,~Zorich]
\label{prop:saddle:multiplicity}
Let us assume that $\gamma$ and $\eta$ are two compact separatrices. If
$\gamma$ and $\eta$ are \^homologous then they are parallel and their
lengths are equal or differ by a factor two:
$$
\cfrac{|\gamma|}{|\eta|} \in \left\{ 1,2,\cfrac{1}{2} \right\}
$$
\end{Proposition}

\begin{Example}
In Figure~\ref{fig:suspension:gp}, one can check that the vertical
saddle connection $\gamma(\pi)$ and the vertical separatrix loop
$\eta$ are \^homologous. More precisely, we have 
$\widehat{\gamma(\pi)} = \hat \eta$.
\end{Example}

\subsubsection{Multiplicity}
\label{subsec:mult}

Let $\gamma$ be a compact separatrix on $\Su$. We will say that $\gamma$
has {\it multiplicity $n$}, and we will note $\textrm{mult}(\gamma)=n$,
if there exists exactly $n$ different classes $[\eta]\in H_1$ where
$\eta$ is a compact separatrices \^homologous to $\gamma$.

We will say that a simple cylinder has {\it multiplicity $n$} if the
multiplicity of its boundary (which can be represented by a compact
separatrix) is $n$.

\begin{Lemma}
\label{lm:characteriz:hyper}
Let $(\Su,\Dq)$ be any half-translation surface in any
hyperelliptic connected component of the following type:
$$
\begin{array}{l}
\Q^{hyp}(4(g-k)-6,4k+2) \\
\Q^{hyp}(2(g-k)-3,2(g-k)-3,4k+2) \\
\Q^{hyp}(2(g-k)-3,2(g-k)-3,2k+1,2k+1)
\end{array}
$$
For each of these strata, let $\gamma$ be a saddle connection on $(\Su,\Dq)$
between (respectively)
\begin{itemize}
\item the two zeroes of degree $4(g-k)-6$ and $4k+2$.

\item one of the two zeroes of degree $2(g-k)-3$ and the other
zero of degree $4k+2$.

\item one of the two zeroes of degree $2(g-k)-3$ and one of the
two zeroes of degree $2k+1,2k+1$.
\end{itemize}

Then $\gamma$ has multiplicity at least two.
\end{Lemma}

\begin{proof}[Proof of the lemma]
The surface $\Su$ is hyperelliptic so it is equipped with an
hyperelliptic involution, say $\tau$. Take a saddle connection
$\gamma$ as indicated in the assumptions. By
construction~\cite{Lanneau}, $\gamma_2:=\tau(\gamma) \not =
\gamma$ thus we obtain an other saddle connection in the same
direction of $\gamma$ (and of the same length). Always by
construction we get that $\hat \gamma = \hat \gamma_2$. Therefore
$\gamma_2$ is \^homologous to $\gamma$ and $[gamma]\not = \gamma_2$
hence $\textrm{mult}(\gamma) \leq 2$. The lemma is proven.
\end{proof}

\subsection{Surgeries}
\label{sec:sub:surgeries}

The principal ingredient of our proof is to decrease the complexity of
a stratum. By complexity, we mean the genus $g$ of the surfaces and
the number $n$ of the singularities. The next two subsections describe
how one can increase the complexity. In
Section~\ref{sec:properties:permutations} we will present some results
to get the converse. For details and proofs of the two next sections we refer
to~\cite{Eskin:Masur:Zorich}, \cite{Kontsevich:Zorich} and~\cite{Masur:Zorich}.

\subsubsection{Breaking up a singularity}

Let $(\Su,\Dq)$ be a half-translation surface and $P \in \Su$ a
singularity of $\Dq$. Let $k\pi$ be the conical angle around $P$
with $k\leq -1$. Choose any partition of $k$ into two non-zero
integers $k_1,k_2$ with $k_i \geq -1$. We recall the well known
construction to obtain a new half-translation structure $\Dq'$ on
$\Su$ with the same singularities pattern as $\Dq$ except at the
point $P$; The new half-translation surface will possess two
singularities $P_{k_1}$ and $P_{k_2}$ of multiplicities $k_1$ and
$k_2$. Moreover $\Su'$ will possess a multiplicity one saddle
connection between $P_{k_1}$ et $P_{k_2}$. 
Here we detail the case $k$ odd and the case $k_1,k_2$ even.

\begin{figure}[htbp]
\begin{center}
\psfrag{6p}{\scriptsize $5\pi$}
\psfrag{e}{\scriptsize $\varepsilon$}
\psfrag{d1}{\scriptsize $\delta$}
\psfrag{e+}{\scriptsize $\varepsilon+\delta$}
\psfrag{e-}{\scriptsize $\varepsilon-\delta$}

\includegraphics[width=13cm]{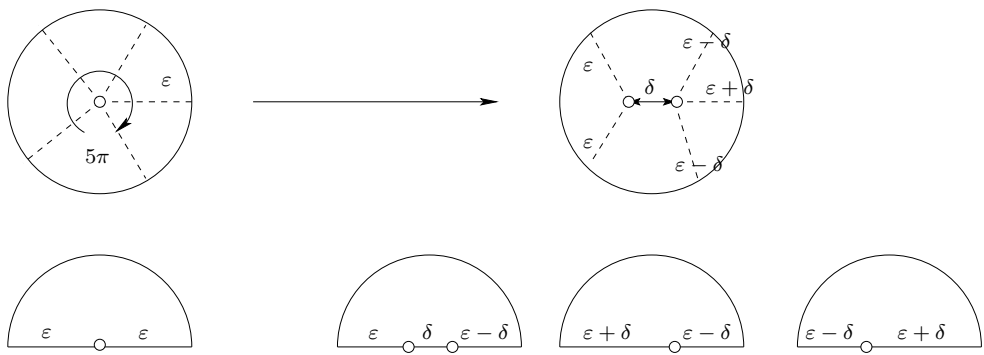}
\end{center}
\caption{ \label{fig:break:3:12} Breaking up a zero of order $3$
into two zeroes of orders $1$ and $2$ correspondingly. Note that
the surgery is local: we do not change the flat metric outside of
the neighborhood of the zero.}
\end{figure}

Consider a small geodesic neighborhood of $P$, that is an
$\varepsilon-$``polydisc'' constructed from $k+2$ half Euclidean
discs of radii $\varepsilon$ glued in their centers;
Figure~\ref{fig:break:3:12} (see also also Figure~$4$ and
Figure~$5$ in~\cite{Lanneau}).

Now, for $\varepsilon$ small enough, there is no critical geodesic
passing through this polydisc, other than the $k+2$ emanating rays
from $P$. Let us remove this polydisc and change continuously
parameters in the following manner. We break up the singularity
$P$ into two singularities of conical angle $(k_1+2)\pi$ and
$(k_2+2)\pi$. This is possible by the choice of $k_i$ (see
Figure~\ref{fig:break:3:12} for details). Now one can re-glue this
polydisc on the surface $\Su$ to obtain the desired new
half-translation structure on $\Su$.
\begin{Remark}
The new geodesic for $\Dq'$ on $\Su$ which connects $P_{k_1}$ and
$P_{k_2}$ has multiplicity one. Indeed one can choose its length
arbitrary small without changing the others lengths thus 
Proposition~\ref{prop:saddle:multiplicity} applies. 

\noindent If $\mathcal C$ denotes the component which contains $(\Su,\Dq)$
and $\mathcal C'$ denotes the component which contains $(\Su,\Dq')$
then $\mathcal C'$ is adjacent to $\mathcal C$.

\noindent The construction we have presented is local: we do not change the metric
outside the $\varepsilon-$polydisc. In the case where the parameters
$k_1$ {\it and} $k_2$ are odd, the construction is global. This is 
the {\it parallelogram construction} (see~\cite{Masur:Zorich}).
\end{Remark}

\begin{Theorem}[Masur,~Zorich]
\label{theo:collapse:singu}
Let $(\Su,\Dq)\in \Q(k_1,\dots,k_n)$ be a point. Let us also assume that
there exists a multiplicity one saddle connection on $\Su$ between
the singularities $P_1\in \Su$ of order $k_1$ and $P_2\in \Su$ of
order $k_2$. We will make that additional assumption 
that $\{k_1,\ k_2\} \not = \{ -1,\ -1 \}$. \\
Then there exists a point $(\Su',\Dq') \in Q(k_1+k_2,k_3,\dots,k_n)$ 
such that one can break up the singularity $P'_{k_1+k_2}\in \Su'$ 
(for $\Dq'$) into two singularities to obtain the initial 
half-translation surface $(\Su,\Dq)$.
\end{Theorem}

\begin{proof}[Proof of Theorem~\ref{theo:collapse:singu}]
Here we address the proof in cases of $k_1,\ k_2$ even or $k_1+k_2$
odd. Assume that $\gamma$ is a multiplicity one saddle connection
between $P_1$ and $P_2$. As usual we will assume that $\gamma$ is
vertical. Using the geodesic flow, we contract it to a short
segment of length $\delta$. Choose any $\varepsilon$ with
$\frac{\delta}{2} < \varepsilon$. Now consider an
$\varepsilon-$polydisc $D(\varepsilon)$ of these two points as
indicated in Figure~\ref{fig:break:3:12}. The assumption on the
multiplicity of $\gamma$ implies that one can and do choose
$\varepsilon$ and $\delta$ small enough so that there are no
critical (vertical) geodesics inside $D(\varepsilon)$ other than
$\gamma$ and the $k_1+2$ and $k_1+2$ verticals emanating from
$P_1$ and $P_2$. We can apply above construction to replace
this polydisc by a new one, where one has glued the two
singularities together. The new surface $(\Su,\Dq')$ we have
constructed satisfies the conclusions of the theorem: 
on can break up the singularity $P=P_1=P_2$ into two to obtain our
initial surface $(\Su,\Dq)$. The theorem is proven.
\end{proof}

\noindent From the previous construction we easily gets the
following useful result.

\begin{Corollary}
\label{cor:continuous}
Let $\Su$ be a flat surface in the stratum $Q(k_1+k_2,k_3,\dots,k_n)$. Let
us assume that $k_1+k_2 \not = k_i$ for all $i=3,\dots,n$. 
Then all surfaces obtained from $\Su$ by breaking up the
singularity of multiplicity $k_1+k_2$ with discrete parameters
$k_1,k_2$ (and arbitrary continuous parameters) belong to the same
connected component of the stratum $Q(k_1,k_2,k_3,\dots,k_n)$. \medskip

\noindent In terms of adjacency: If $\mathcal C_0 \subseteq
\Q(k_1+k_2,\dots,k_n)$ is a component and $\mathcal C_1,\mathcal C_2
\subset \Q(k_1,k_2,\dots,k_n)$ are two components such that $\mathcal C_0
\subset \overline{\mathcal C_1} \cap \overline{\mathcal C_2}$, then 
$\mathcal C_1=\mathcal C_2$.
\end{Corollary}

\subsubsection{``bubbling a handle''}

Let $(\Su,\Dq)$ be a half-translation surface and $P \in \Su$ a
singularity of $\Dq$. 
Let us break up the singularity $P$ into two singularities $P_1,\
P_2 \in \Su$ (see previous section). One gets a new half-translation
structure, say $\Dq_1$ on $\Su$, and a closed saddle connection
$\gamma$ (of length $\delta$) between $P_1$ and $P_2$. Let us cut
this surface along $\gamma$. We obtain a surface with some
boundaries components. 
We identify the two points $P_1,\ P_2$ on this surface to obtain a
surface $\Su_1$ with a boundary component isometric to the union of two
circles, each of length $\delta$. Then let us glue a straight metric
cylinder with the following parameters. The height and the twist are
chosen arbitrarily, and of weight (circumference) is $\delta$. 
The new surface $(\Su',\Dq')$ we get is
a ($\textrm{genus}(\Su)+1$) half-translation surface. The angle
between the new handle is $k_1+2$ (or $k_2+2$ if we consider
the complementary angle). See Figure~\ref{fig:simple:cyl} for an
example.

Let $(\Su,\Dq)$ be any half-translation surface. Assume that the
singularities data of $\Su$ are either $(4g-3,-1)$ or $(4g-4)$. Let
$s$ be any non-negative integer.
Thanks to Corollary~\ref{cor:continuous}, the surfaces
$(\Su',\Dq')$ constructed by the surgery ``bubbling a handle'' at the
unique zero of $\Dq$ in $\Su$ with arbitrary continuous parameters
(height, twist and width of the cylinder) and fix discrete parameter $s$
represent quadratic differentials belonging in the same connected
component of the stratum. This motivates the following definition.

\begin{Definition}
Let $(\Su,\Dq)$ be any half-translation surface, with singularities pattern
$(4g-3,-1)$ or $(4g-4)$ and let $s$ be any non-negative integer. We will
denote by $\Su \oplus s$ a  surface constructed by the surgery
``bubbling a handle'' at the unique zero of $\Dq$ in $\Su$ with
arbitrary continuous parameters (height, twist and width of the
cylinder) and discrete parameter $s$. 
If $\mathcal C$ denote the component containing $(\Su,\Dq)$ then we denote 
the component containing $(\Su',\Dq')$ by $\mathcal C \oplus s$. In
other terms we obtain the two well defined mappings.
$$
\begin{array}{ll}
\textrm{for any } g\geq 2, \qquad \oplus = \left(
\begin{array}{lllll}
\pi_0(\Q(4g-3,-1)) & \times & \mathbb N & \longrightarrow & \pi_0(\Q(4g+1,-1)) \\
& (\mathcal C, s ) & & \longmapsto & \mathcal C \oplus s
\end{array} \right) \\ \\
\textrm{for any } g\geq 3, \qquad \oplus = \left(
\begin{array}{lllll}
\pi_0(\Q(4g-4)) & \times & \mathbb N & \longrightarrow & \pi_0(\Q(4g)) \\
& (\mathcal C, s ) & & \longmapsto & \mathcal C \oplus s
\end{array} \right)
\end{array}
$$
Here $\pi_0(E)$ denotes the set of connected components of the
topological space $E$. Note also that the discrete parameter $s$
corresponding to the angle between the two new sectors can be chosen
modulo $2g$ ($s \in \{1,\dots,2g\}$).
\end{Definition}

\begin{Proposition}
\label{prop:properties:oplus}
Let $\mathcal C$ be any connected component of any above stratum, namely 
$\Q(-1,4g-3)$ or $\Q(4g-4)$. Then the following statements hold. \\ \\
$
\begin{array}{lll}
(1) & \mathcal C \oplus s_1 \oplus s_2 = \mathcal C \oplus s_2 \oplus
s_1 & \textrm{ for any } s_1,\ s_2. \\
(2) & \mathcal C \oplus s_1 \oplus s_2 = \mathcal C \oplus
(s_2-2) \oplus (s_1+2) & \textrm{ for any } s_2 \geq 3. \\
(3) &\mathcal C \oplus s_1 \oplus s_2 = \mathcal C \oplus
(s_2-4) \oplus s_1 & \textrm{ for any }s_2-s_1 \geq 4.
\end{array}
$
\end{Proposition}

\begin{proof}[Proof of Proposition~\ref{prop:properties:oplus}]
The proof uses the description of quadratic differentials in terms 
of separatrices diagrams (see~\cite{Kontsevich:Zorich}) and 
ribbons graphs. For such diagrams, an element of $\Q(4g-4)$ 
(respectively of $\Q(-1,4g-3)$) is presented by a measured foliation 
with $4g-4+2$ (respectively $4g-3+2$) leaves emanating from the 
singularity. Gluing an handle of angle $s\pi$ consists topologically 
to glue two news sectors with an angle of $s\pi$. We thus get 
measured foliation with $4g-4+2+4$ (respectively $4g-3+2+4$) leaves
around the new singularity. Geometrically it is easy to see that 
the two surgeries consisting to first glue the handle of angle $s_1\pi$ 
or to first glue the handle of angle $s_2\pi$ produces surfaces 
belonging to the same component. This is the first point of the 
proposition. \\
The second statement is also clear in terms of diagrams. 
This is illustrated in Figure~\ref{fig:diag}.
\begin{figure}[htbp]
\begin{center}
\psfrag{s1}{\scriptsize $s_1\pi$}
\psfrag{s2}{\scriptsize $s_2\pi$}
\includegraphics[width=8cm]{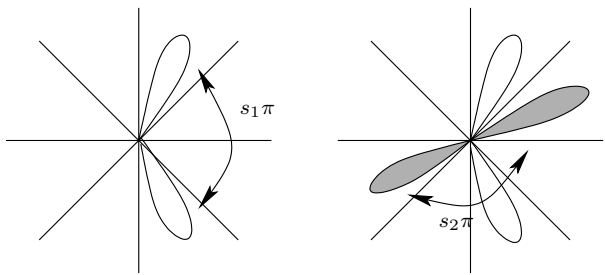}
\end{center}
\caption{ 
\label{fig:diag}
Gluing two handles attached on a singularity. Gluing first 
the ``white'' handle of angle $s_1\pi$ and then the ``grey'' handle 
of angle $s_2\pi$ is equivalent to glue first the ``grey'' 
handle of angle $(s_2-2)\pi$ and then the ``white'' handle of 
angle $(s_1+2)\pi$.
}
\end{figure}

\noindent 
The last statements follows from the same considerations. The proposition 
is proven.
\end{proof}

\subsection{Jenkins-Strebel surfaces}

An important class of flat metrics is given by the so-called
Jenkins-Strebel differentials. We first explain these particular
metrics and then we give a density result.

Let us denote by $\Gamma(\Dq)$ the critical graph of $(\Su,\Dq)$
induced by the horizontal foliation, 
that is the union of all separatrices in the horizontal direction.
It is easy to see that $\Gamma$ is compact if and only if
the horizontal measured foliation of $\Dq$ is completely periodic;
this means that any horizontal geodesic is closed. We then call such
a form a {\it Jenkins-Strebel differential}.

Recall that locally, a stratum of the moduli space $\Q_g$ is
modeled by the first cohomology group with coefficients in
$\mathbb C$:
\begin{multline*}
H^1(\Su,\{P_1,\dots,P_n\};\mathbb Q)\oplus i \cdot H^1(\Su,\{P_1,\dots,P_n\};\mathbb Q) \subset \\
\subset
H^1(\Su,\{P_1,\dots,P_n\};\C)=H^1(\Su,\{P_1,\dots,P_n\};\R\oplus
i\R)\,.
\end{multline*}
Taking forms such that the absolute and relative periods are in
$\mathbb Q \oplus i \mathbb Q$, we get arithmetic surfaces: the
orientating covering is itself a ramified covering over the standard
two torus. Therefore, up to normalize the direction to the horizontal direction,
we get that Jenkins-Strebel differentials are dense in each stratum 
(see~\cite{Douady:Hubbard}, \cite{Hubbard:Masur},
\cite{Kontsevich:Zorich} and~\cite{Strebel}).

\begin{Remark}
In~\cite{Douady:Hubbard}, Douady and Hubbard proved a
stronger result: the Jenkins-Strebel differentials are dense on each
Riemann surface, not just in a stratum.
\end{Remark}

The complement of $\Gamma(\Dq)$ in $\Su$ is a disjoint union of
{\it maximal} periodic components for the horizontal foliation.
These components are isometric to metric straight cylinders,
foliated by regular horizontal leaves. A simple computation with the
Euler characteristic, using the cylinders decomposition, shows that the
number of the cylinders is bounded by $3g-3$.

In~\cite{Masur:79} Masur proved that the set of Jenkins-Strebel
differentials with exactly $r$ cylinders (for any $1 \leq r \leq
3g-3$) is dense inside the principal stratum $\Q(1,\dots,1)$ of
genus $g$. Recently, Kontsevich and
Zorich~\cite{Kontsevich:Zorich} have obtained a similar result. They
have proved that the set of Jenkins-Strebel differentials with exactly one
cylinder is dense inside any stratum of $\HA_g$. Here 
we extend their proof to the case of strata of $\Q_g$.

\begin{Theorem}
\label{theo:density}
The set of quadratic differentials, such that the horizontal
foliation is completely periodic and decomposes the surface into a
unique straight metric cylinder, forms a dense subset of any connected
component of any stratum of $\Q_g$.
\end{Theorem}

\begin{proof}[Proof of Theorem~\ref{theo:density}]
We will only prove the existence of such surfaces on each
connected component of the moduli space $\Q_g$. For the density
result, we refer to the Kontsevich-Zorich's proof. Note that in
this paper, we will only need the existence result.

Let $(\Su,\Dq)$ be a point. Thanks to previous discussion we may
assume that the surface $(\Su,\Dq)$ is an arithmetic surface:
its orientating double covering $\hat \Su$ covers the standard torus 
$\mathbb T^2=\mathbb C/\mathbb Z^2$. 
The vertical foliation on $\Su$ is completely periodic and
decomposes the surface into many vertical cylinders $\mathcal C_i$.

Let us construct a closed regular curve $\gamma$ transverse to
this foliation. The surface $\hat \Su$ is a ramified covering 
$\pi : \hat \Su \rightarrow \mathbb T^2$. Obviously the directional
flow on $\hat \Su$ projects to the directional flow on the two-torus
$\mathbb T^2$. Let us consider a foliation on this torus in the
direction $\theta = 1/b$ with $b$ arbitrary large. The lift of
this foliation allows us to obtain a closed regular geodesic
$\gamma$ on $\hat \Su$, and thus a closed geodesic on $\Su$. This
leaf is transverse to the vertical foliation determined by $\Dq$ and 
$\gamma$ does not contain any singularity of $\Dq$. In addition 
one can choose $\gamma$ in such a way that its length with respect to
the metric defined by $\Dq$ is arbitrary large; Indeed the length of
$\gamma$ is greater than $\sqrt{1+b^2}$.

The closed loop $\gamma$ cuts the boundaries of cylinders $\mathcal
C_i$ many times (this means that $\gamma$ cuts the set of vertical
saddle connections and separatrix loops of $\Dq$). By construction
$\partial \mathcal C_i \backslash \gamma$ is a disjoint union of
vertical intervals. One can and do choose $\gamma$ long enough so
that each components of $\partial \mathcal C_i \backslash
\gamma$ contains at most one singularity of $\Dq$. Now we will
change slightly the transverse structure (in the vertical
direction) to obtain a periodic horizontal foliation with only one
cylinder. We will do that without changing the structure in the
direction of $\gamma$.

We cut the surface along the vertical critical graph $\Gamma(\Dq)$
of $\Dq$ and also along $\gamma$. We obtain a finite union of
parallelograms $R_i$. Up to the $\Sl-$action one may assume that
$\gamma$ is vertical. The set of horizontal sides of $\R_i$ is a
part of $\gamma$ and the set of vertical sides of $\R_i$ is a part
$\Gamma(\Dq)$. By construction in each vertical side of $R_i$
there is at most one singularity of $\Dq$.

Let us construct a new foliation as follows. We conserve all
horizontal parameters and we change vertical parameters in the
following way: we declare that the length of any vertical side of
$R_i$ is $1$ for all $i$. In addition, if there is a singularity
located on a vertical side, we declare that it is located at the
middle of this side. With our above considerations, there is no
contradiction. Finally we obtain a new set of parallelograms $R_i'$
endowed with the natural metric $\D z^2$.

Let $(\Su',\Dq')$ be the flat surface constructed from the new
rectangles $R_i'$ with the corresponding identifications of
vertical and horizontal sides given by gluing described above. We
obtain a new half-translation structure $\Dq'$ on our surface.

The surfaces $(\Su,\Dq)$ and $(\Su',\Dq')$ carry the same
topology. By construction the vertical critical graphs
$\Gamma(\Dq)$ and $\Gamma(\Dq')$ coincide. We just have change
absolute and relative periods of the form $\Dq$. The subvariety of
quadratic differentials sharing the same vertical foliation is
connected and depends continuously on the suitable of deformations
of the vertical foliation (see~\cite{Hubbard:Masur}
and~\cite{Veech:flow}). Thus it implies that the two points
$(\Su,\Dq)$ and $(\Su',\Dq')$ belong in the same connected
component. It is easy to check that in the 
horizontal direction on $\Su'$ for $\Dq'$ the foliation is
completely periodic and decomposes the surface into a single
cylinder.
\end{proof}

\section{Generalized permutations}
\label{sec:generalized:permutations}

In this section, we propose a natural way to encode
Jenkins-Strebel differentials; namely we will introduce the
notion of {\it generalized permutations}.

In this section $(\Su,\Dq)$ will denote a surface which is decomposed
into a unique metric cylinder into the horizontal direction. Note that
the holonomy of $(\Su,\Dq)$ is not assumed to be non-trivial. 

\subsection{Combinatorics of $(\Su,\Dq)$}

Let $\Gamma(\Dq)$ be the horizontal critical graph of $(\Su,\Dq)$. The complement
$\Su \backslash \Gamma(\Dq)$ is a metric cylinder $\CYL$. 
The horizontal saddle connections are labelled by 
$\gamma_1,\dots,\gamma_k$.

One can encode the sequence of saddle connections contained in the
bottom and in the top of $\CYL$, and ordered following
the cyclic ordering of the boundary of the cylinder, by a sequence of
labels into the following manner. 
Each saddle connection $\gamma_i$ is presented
twice on the boundaries of $\CYL$. Let us denote by 
$\gamma_i^1$ and $\gamma_i^2$ these two copies. Hence the top
(respectively the bottom) of $\CYL$ is a sequence of
$\gamma_i^\epsilon$ where $i\in \{1,\dots,k\}$ and $\epsilon=1,2$.
Roughly speaking a generalized permutation $\pi$ is a table with two
lines encoding the sequence of labels of the saddle connections.
The precise definition is the following.

\begin{Definition}
Let $r,l$ be any non-negative integers with the same parity. A
{\it generalized permutation} $\pi$ (of type $(r,l)$) is an
involution, without fixed points, of the set
$\{1,\dots,r,r+1,\dots,r+l\}$.
\end{Definition}
Through this paper, we will represent a generalized permutation by a
table.

\begin{Example}
\label{example:perm}
The generalized permutation (of type $(7,5)$) given by
$\pi(1)=10,\ \pi(2)=12,\ \pi(3)=5,\ \pi(4)=7,\ \pi(5)=3,\ \pi(6)=11,
\ \pi(7)=4,\ \pi(8)=9,\ \pi(9)=8,\ \pi(10)=1,\ \pi(11)=6,\ \pi(12)=2$
is represented by the table 
$$
\pi= \left(
\begin{array}{ccccccc}
1 & 2 & 3 & 4 & 3 & 5 & 4 \\
6 & 6 & 1 & 5 & 2
\end{array}
\right)
$$
in a natural way.
\end{Example}
The term ``generalized'' is justified by the fact that a classical permutation
$\pi_1$ of the symmetric group $S_k$ is a generalized permutation with
$l=r=k$ and
$$
\pi(i) = \left\{
\begin{array}{lr}
k + \pi_1^{-1}(i), & \textrm{for } i \leq k \\
\pi_1(i-k), & \textrm{for } i > k
\end{array}
\right.
$$
In the present paper, we are interested by the quadratic
differentials which are not the global square of any Abelian
differential. Thus, in order to avoid ``true'' permutations, we
require the following technical condition
\begin{equation}
\label{eq:tech:cond}
\textrm{there exist } i_0 \leq r \textrm{ and } j_0 \geq r+1
\textrm{ such that } \pi(i_0) \leq r \textrm{ and }\pi(j_0) \geq
r+1
\end{equation}

\subsection{Admissible vectors}
\label{sec:admissible}

Let $\pi$ denote a generalized permutation of type $(r,l)$.

\begin{Definition}
We say that $\lambda \in \mathbb R_{+}^{r+l}$ is an admissible
vector (for $\pi$) if
\begin{displaymath}
\label{eq:rationnal}
\begin{cases}
\lambda_i = \lambda_{\pi(i)} \textrm{ for all } i=1,\dots,r+l \\
\\
\sum_{i=1}^{r} \lambda_i = \sum_{j=1}^l \lambda_{r+j}
\end{cases}
\end{displaymath}
\end{Definition}
Note that for the ``true'' permutations $\pi$ the vector
$$
(\lambda_1,\dots,\lambda_r,\lambda_{\pi(1)},\dots,\lambda_{\pi(r)})
$$
is admissible for any $\lambda_i > 0$.

\subsection{Suspension over a generalized permutation}
\label{sec:suspension}

Let $\pi$ be a generalized permutation and let $\lambda$ be any
admissible vector for $\pi$. We denote by $w$ (width of the
cylinder) the quantity
\begin{equation*}
w:=\sum_{i=1}^{r} \lambda_i = \sum_{j=1}^l \lambda_{r+j}
\end{equation*}
Let $\mathcal C=[0,w] \times [0,1]$ be an Euclidean cylinder 
endowed with the form $\D z^2$. Let us consider the partition of the
top (respectively bottom) of $\mathcal C$ into horizontal intervals of
length $\lambda_i$ for $i=1,\dots,r$ (respectively
$i=r+1,dots,r+l$). Now we identify the horizontal interval labeled 
$i$ with the horizontal interval $\pi(i)$ for all $i$ into the
following manner. If the two intervals are presented twice on a side,
we identify them by a centrally symmetry and otherwise we identify
them by a translation. 

The resulting space is a Riemann surface, denoted by
$\Su(\pi,\lambda)$, endowed with a natural quadratic differential
$\Dq = \D z^2$. We call this half-translation surface $(\Su,\Dq)$
the {\it suspension} over the element $(\pi,\lambda)$.
\begin{Notation}
The surface $\Su=\Su(\pi,\lambda)$ decomposes into a single
cylinder in the horizontal direction. By construction there exists 
in the vertical direction a compact separatrix on this
surface. We will denote it by $\gamma(\pi) \subset \Su$.
\end{Notation}

\begin{Remark}
If $\pi$ a true permutation then the surface $\Su=\Su(\pi,\lambda)$
is a {\it translation surface}. The converse is true.
\end{Remark}

\begin{Lemma}
The surfaces $\Su(\pi,\lambda)$ with fixed parameter $\pi$ and
arbitrary parameters $\lambda$ belong to the 
same connected component.
\end{Lemma}

\begin{proof}[Proof of the lemma]
The lengths $\lambda_i$ of the horizontal intervals correspond to
the absolute and to the relative periods of the corresponding form
$\Dq$ on $\Su$. Thus the lemma is a direct consequence of the
local description of the orbifoldic structure of the strata in terms
of the cohomological coordinates.
\end{proof}
This construction implies a simple but important fact: we can
encode the set of connected components using generalized
permutations. Given a permutation, it determines completely the
type of the singularities and hence a stratum. In addition, above
lemma shows that it also determines the connected component of the
stratum as well. The set of generalized permutations, for a fix
stratum is obviously finite. Thus it gives an independent proof of
a theorem of Veech~\cite{Veech:flow}.

\begin{NoNumberTheorem}[Veech]
The set of connected components of a stratum $\Q(k_1,\dots,k_n)$
of the moduli space $\Q_g$ of meromorphic quadratic differentials
is finite.
\end{NoNumberTheorem}

\subsection{Cyclic order and horocyclic flow}
\label{sec:cyclic}

The class of Jenkins-Strebel quadratic differentials is stable under
the horocyclic flow $h_s=\left( \begin{smallmatrix} 1 & s \\ 0 & 1
\end{smallmatrix} \right)$. More precisely the corresponding suspended
surfaces
$$
h_s \cdot \Su(\pi,\lambda) = \Su(\pi',\lambda')
$$
are related into the following way. Elements of $\pi'$ corresponds to
the elements of $\pi$ with a ``rotating'' of the first line and the
second line. We will say that $\pi$ is equivalent to $\pi'$ and we
will note $\pi \sim \pi'$. 
For example the permutation of Example~\ref{example:perm}
is equivalent to the following one
$$
\pi = \left(
\begin{array}{ccccccc}
1 & 2 & 3 & 4 & 3 & 5 & 4 \\
6 & 6 & 1 & 5 & 2
\end{array}
\right) \sim
\left(
\begin{array}{ccccccc}
1 & 2 & 3 & 4 & 3 & 5 & 4 \\
1 & 5 & 2 & 6 & 6
\end{array}
\right) = \pi'
$$
Note that this relation preserves the stratum and also the
connected component.
\begin{figure}[htbp]
\begin{center}
\psfrag{1}{$1$}  \psfrag{h}{$h_s$}
\psfrag{e}{$\scriptstyle \textrm{action of the mapping class group}$}

\includegraphics[width=12cm]{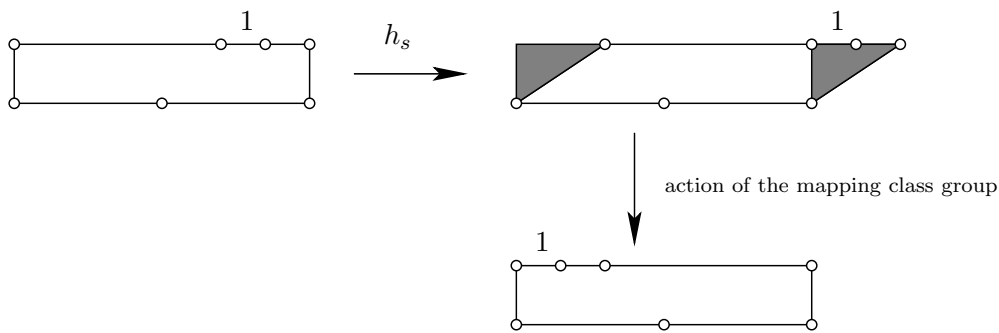}
\end{center}
\caption{ 
\label{fig:horocycle} 
Action of the horocyclic flow $h_s$ on a half-translation surface.}
\end{figure}

\section{Representative elements}
\label{sec:representative: elements}

In this section, we give a ``bestiary'' of half-translations
surfaces in terms of generalized permutations (see previous 
Section~\ref{sec:generalized:permutations}).

\subsection{Hyperelliptic connected components}
\label{sec:permu:pii}

Let $r,l$ be two non-negative integers. Let $k=r+l$. We construct the
generalized ``hyperelliptic'' permutation $\Pi_1(r,l)$ of type $(k+2,k+2)$
into the following way.
$$
\Pi_1(r,l) = \left( \begin{array}{cccccccc}
0_1 & 1 & \dots & r & 0_1 & r+1 & \dots & r+l \\
r+l & \dots & r+1 & 0_2 & r & \dots & 1 & 0_2
\end{array} \right)
$$
A direct computation of the type of conical angles identifies
the stratum which contain surfaces $\Su(\Pi_1(r,l),\lambda)$.
It depends of the parity of the two integers $r$ and $l$. 
We easily establish the next lemma.

\begin{Lemma}
\label{lm:strate:perm:n-m}
Let $(\Su,\Dq)$ be the half-translation surface given by
$\Su(\Pi_1(r,l),\lambda)$ for any admissible vector $\lambda$. If
$r$ and $l$ are odd then $\Dq$ has two singularities. If $r$ and
$l$ have different parities then $\Dq$ has three singularities. If
$r$ and $l$ are even then $\Dq$ has four singularities. The
following table gives the type of the singularities in terms of
$r$ and $l$:
$$
\begin{array}{|c|c|c|}
\hline
r & l & \textrm{stratum which contains $(\Su,\Dq)$} \\
\hline
2k+1 & 2(g-k)-3 & \Q( \ 4k+2 \ , \ 4(g-k)-6 \ ) \\
2k+2 & 2(g-k)-3 & \Q( \ 2k+1 \ , \ 2k+1 \ , \ 4(g-k)-6 \ ) \\
2k+2 & 2(g-k)-2 & \Q( \ 2k+1 \ , \ 2k+1 \ ,\  2(g-k)-3 \ , \ 2(g-k)-3 \ ) \\
\hline
\end{array}
$$
\end{Lemma}
According to~\cite{Lanneau}, each above stratum contains an hyperelliptic connected component.
We have the following lemma.

\begin{Lemma}
\label{lm:hyper:pi:n-m}
For any $\lambda$ and any integers $r,\ l$, the surface
$\Su(\Pi_1(r,l),\lambda)$ belongs to the hyperelliptic connected
component (of the corresponding stratum).
\end{Lemma}

\begin{proof}[Proof of Lemma~\ref{lm:hyper:pi:n-m}]
Here we present the proof of the first case; that is $r$ and $l$
are odd. The other cases are similar and left to the reader.
Take $r=2k+1$ and $l=2(g-k)-3$. We consider the rectangle
$$
\mathcal R \qquad = \qquad \left] -\cfrac{r+l}{2}-1\ ,\ \cfrac{r+l}{2}+1
\right[ \qquad \times \qquad ]-1,1[
$$
Let $\tau : \mathcal R \rightarrow \mathcal R$ be the involution of
$\mathcal R$ given by
$\tau(x,y)=(-x,-y)$. The combinatorics of $\Pi_1$ implies that
$\tau$ induces a {\it global} involution of the surface
$\Su=\Su(\Pi_1(r,l),\lambda)$; we still denote it by $\tau$. By
Lemma~\ref{lm:strate:perm:n-m} the surface $(\Su,\Dq)$ belongs to
the stratum $\Q(4k+2,4(g-k)-6)$ therefore the Gauss-Bonnet
formula implies that $\Su$ has genus $g$.

Recall that the hyperelliptic component of this stratum is, by
definition, the image of the map
\begin{gather*}
\Q(2k,2(g-k)-4,{-1}^{2g}) \to \Q(4k+2,4(g-k)-6) \\
(\CP,\Dq_0) \mapsto (\Su_0,\pi^{\ast}\Dq_0)
\end{gather*}
where $\pi : \Su_0 \rightarrow \CP$ is a double ramified covering.
The locus of ramification being the zeroes and poles of $\Dq_0$.

In order to prove that $(\Su,\Dq)$ belongs to the hyperelliptic
component, we have to construct a double ramified  covering $\pi :
\Su \to \mathbb P^1$ and a quadratic differential $\Dq_0$ on
the sphere such that $\pi^{\ast}\Dq_0=\Dq$.
Let us count the number of fixed points of the map $\tau$:

\begin{itemize}

\item there are $r+l$ fixed points of $\tau$ on the horizontal
sides of $R$ located at the middle of the intervals (precisely at
the middle of separatrix loops).

\item there is a fixed point located at the middle of the vertical
side.

\item there is a fixed point at $(0,0)$.

\item there are $2$ fixed points which corresponds to the two
zeroes of $\Dq$.

\end{itemize}

Therefore the total number of fixed points of $\tau$ on $\Su_0$ are
$$
r+l+1+1+2 = 2k+1 + 2(g-k)-3 + 4 = 2g+2
$$
The Riemann-Hurwitz formula implies that the genus of
$\Su \left/ (x \sim \tau(x)) \right.$ is zero. Let us consider the
projection map
$$
\pi :  \Su \to \Su \left/ (x \sim \tau(x)) \right.\simeq \CP.
$$
It is easy to check that above covering gives the desired
map. Lemma~\ref{lm:hyper:pi:n-m} is proven.
\end{proof}

\begin{Remark}
The permutation $\Pi_2(r,l)$ of type $(2r,2l)$ given by:
$$
\Pi_2(r,l) = \left( \begin{array}{cccccc}
1 & \dots & r & 1 & \dots & r \\
r+1 & \dots & r+l & r+1 & \dots & r+l
\end{array} \right)
$$
furnishes also a representative element for hyperelliptic connected
components.
\end{Remark}

Through the proof of Theorem~\ref{theo:anexe:1} on adjacency of
strata, we will get a characterization of hyperelliptic components. We
will obtain the nice description.
\begin{Theorem}
Let $(\Su,\Dq)$ be any point into any hyperelliptic component. Let us assume
that $\Su$ is decomposed into a unique metric cylinder for the
horizontal flow. Then there exist $i\in \{1,2\}$, non-negative
integers $r,l$, an admissible vector $\lambda$ and $s \in \mathbb R$ such that
$$
(\Su,\Dq) = h_s \cdot \Su \left( \Pi_i(r,l),\lambda \right)
$$
\end{Theorem}

\subsection{Irreducible connected components}

Here we give representative elements of the sporadic components
discussed in Theorem~\ref{theo:main:2}: 
we will denote them by the {\it irreducible} connected components.

\begin{Definition}
\label{def:irr}
The irreducible connected components are defined to be the components
containing the elements $\Su(\pi,\lambda)$ given by the next table.
$$
\begin{tabular}{|c|l|}
\hline
irreducible components & Representatives elements \\
\hline

$\Q^{irr}(-1,9)$ &
$ \left( \begin{array}{cccccc}
0& 1 & 2 & 3 & 4 & 0 \\
4 & 3 & 2 & 5 & 1 & 5
\end{array} \right)
$  \\ \hline

$\Q^{irr}(-1,3,6)$ & $ \left( \begin{array}{ccccccc}
0 & 1 & 2 & 3 & 4 & 5 & 0 \\
5 & 4 & 3 & 2 & 6 & 1 & 6
\end{array} \right)
$ \\ \hline

$\Q^{irr}(-1,3,3,3)$ & $ \left( \begin{array}{cccccccc}
0 & 1 & 2 & 3 & 4 & 5 & 6 & 0 \\
6 & 5 & 3 & 2 & 7 & 4 & 1 & 7
\end{array} \right)
$  \\ \hline

$\Q^{irr,I}(12)$ & $ \left(
\begin{array}{ccccccc}
1 & 2 & 3 & 4 & 2 & 5 & 6 \\
1 & 4 & 5 & 7 & 6 & 7 & 3
\end{array} \right)
$  \\ \hline

$\Q^{irr,II}(12)$ & $ \left(
\begin{array}{ccccccc}
1 & 2 & 3 & 4 & 3 & 5 & 6 \\
1 & 5 & 7 & 4 & 2 & 6 & 7
\end{array} \right)
$ \\ \hline
\end{tabular}
$$
\end{Definition}

\section{Dynamical properties of $\Su(\pi,\lambda)$ versus
  combinatorics of $\pi$} 
\label{sec:properties:permutations}

An important part of our proof is to find
surfaces with multiplicity one saddle connection. For that we will
give a combinatorial condition such that 
Proposition~\ref{prop:saddle:multiplicity} applies.

\subsection{Irreducibility}

There exists ``bad'' permutations $\pi$ such that, for any $\lambda$,
the saddle connection $\gamma(\pi)$ has multiplicity at least two on
$\Su(\pi,\lambda)$. More precisely the two next cases can occur.
$$
\forall \lambda , \ \exists \textrm{ saddle connection }\eta \textrm{
  such that } |\eta| = |\gamma(\pi_1)| \textrm{ on } \Su(\pi_1,\lambda)
$$
and
$$
\forall \lambda , \ \exists\textrm{ saddle connection } \eta \textrm{
  such that } |\eta| = 2|\gamma(\pi_2)| \textrm{ on } \Su(\pi_2,\lambda)
$$
The first class of ``bad'' permutations will lead to the notion of
weak reducibility. The second class will lead to the $\Red$ notion.

\subsubsection{Weak irreducibility}

We say that $\pi$ (of type $(r,l)$) is {\it weakly reducible} if
there exist $1 \leq i_0 < r$ and $r+1 \leq j_0 < p=r+l$ such that
at least one of the following two conditions holds.
\begin{itemize}

\item $\pi(\{1,\dots,i_0 \}) = \{r+1,\dots,j_0\}$ or
$\pi(\{i_0+1,\dots,r \}) = \{r+j_0+1,\dots,p\}$

\item each $1 \leq i \leq r$ with $1 \leq \pi(i) \leq r$
satisfies $i \leq i_0$ and $\pi(i) > i_0$. All other $i \leq i_0$
with $\pi(i) \geq r+1$ satisfy $\pi(i) \leq j_0$. \\
each $r+1 \leq j$ with $r+1 \leq \pi(j)$ satisfies $j \leq j_0$
and $\pi(j) > j_0$. All other $r+1 \leq j \leq j_0$ with $\pi(j)
\leq r$ satisfy $\pi(j) \leq i_0$.

\end{itemize}
We will say that $\pi$ is {\it weakly irreducible} if $\pi$ is not
weakly reducible. 

\begin{Lemma}
\label{rk:length:1}
Let $\Su=\Su(\pi,\lambda)$ be a half-translation surface. Then 
$\pi$ is weakly irreducible if and only if there exists a 
full Lebesgue measure set of admissible vectors 
$\lambda$ such that any vertical separatrix $\eta$ 
(different from $\gamma(\pi))$ on $\Su$ 
satisfies $|\eta| \geq 2\cdot |\gamma(\pi)|$.
\end{Lemma}

The proof is obvious.

\begin{Notation}
In order to clarify the situation, for a weakly reducible permutation
$\pi$, we will denote by a vertical segment the position of the
corresponding elements $i_0,\ j_0$. For instance the permutation
$$
\left(
\begin{array}{ccc|ccc}
1 & 2 & 3 & 4 & 3 & 5  \\
6 & 1 & 2 & 6 & 5 & 4
\end{array}
\right)
$$
is weakly reducible with corresponding $i_0=3$ and $j_0=9$.
\end{Notation}

\subsubsection{The condition $\Red$}
\label{sec:irred2}

This is a technical condition and we first prefer to address the
relating result to this notion. The definition will then become
clear.

\begin{Proposition}
\label{prop:length:2}
Let $\Su=\Su(\pi,\lambda)$ be a half-translation surface. 
Let us assume that $\pi$ satisfies the condition $\Red$. 
Then there exists a full Lebesgue measure set of admissible vectors 
$\lambda$ such that any vertical separatrix $\eta$ on $\Su$ 
satisfies $|\eta| \not =  2\cdot |\gamma(\pi)|$.
\end{Proposition}

\begin{Definition}
We say that $\pi$ does not satisfy the condition $\Red$ if there
exists a decomposition of $\pi$ (up to exchange lines) into the
following way (in terms of table).
$$
\pi= \left(
\begin{array}{cc|c|cc}
Y_1' & & Y_1'' & & Y_1''' \\
Y_2' & 0 & Y_2'' & 0 & Y_2'''
\end{array}
\right)
$$
with
$$
\begin{array}{ccc}
\begin{cases}
\forall i \in Y_1' \Rightarrow \pi(i) \in Y_1''' \sqcup Y_2' \\
\forall i \in Y_1'' \Rightarrow \pi(i) \in Y_1'' \sqcup Y_2''
\end{cases} & \textrm{and} &
\begin{cases}
\forall j \in Y_2' \Rightarrow \pi(j) \in Y_2''' \sqcup Y_1' \\
\forall j \in Y_2'' \Rightarrow \pi(j) \in Y_2'' \sqcup Y_1''
\end{cases}
\end{array}
$$
\end{Definition}

\begin{Example}
\label{example:perm:condition:red}
The permutation $\left( \begin{smallmatrix}
1 & 2 & 2 & 3 & 3 & 1 \\ 0 & 0 \end{smallmatrix} \right)$ 
does not satisfy the condition $\Red$. Indeed, 
$Y_1'=\{1\}$, $Y_1''=\{2,2,3,3\}$, $Y_1'''=\{1\}$ and 
$Y_2'=Y_2''=Y_2'''=\emptyset$ is a decomposition. \\
The permutation of Example~\ref{example:perm} 
satisfies the condition $\Red$.
\end{Example}

\begin{proof}[Proof of Proposition~\ref{prop:length:2}]
As usual we will assume that $\gamma(\pi)$ has length~$1$.
Recall that an admissible vector $\lambda$ satisfy a linear equation
given by Equation~(\ref{eq:rationnal}). Let us consider the set $E$ of
admissible vectors $\lambda$ such that the entries $\lambda_i$ 
defines a codimension one subspace. This is obviously a full
Lebesgue measure set. Let us choose any $\lambda \in E$. We will
prove that if there exists a vertical separatrix of length~$2$,
then the above decomposition of $\Red$ occurs.

Let us assume that one has found a vertical separatrix
$\eta \in \Su(\pi,\lambda)$ of length $2$. A straightforward
computation shows that 
one of the two cases presented by Figure~\ref{fig:irred2:proof:idoc}
and Figure~\ref{fig:irred2:proof} have to occurs.

\begin{figure}[htbp]
\begin{center}
\psfrag{i0}{$i_0$} \psfrag{i1}{$i_1$} \psfrag{j0}{$j_0$}
\psfrag{j1}{$j_1$} \psfrag{e}{$\varepsilon$}
\psfrag{s}{$\Su(\pi,\lambda)$} \psfrag{g}{$\gamma(\pi)$}
\psfrag{et}{$\eta$}

\includegraphics[width=13cm]{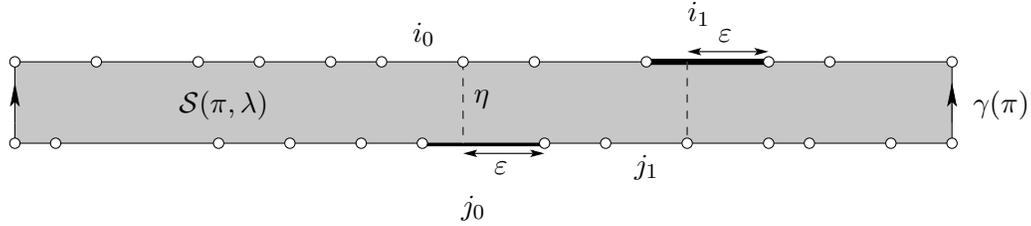}
\end{center}
\caption{ \label{fig:irred2:proof:idoc} A separatrix of length $2$
(here the canonical separatrix $\gamma(\pi)$ has length~$1$). The
two corresponding horizontal intervals of length $\lambda_{j_0+r}$
and $\lambda_{i_1}$, numbered by $j_0$ and $i_1$, are glued by a
translation.}
\end{figure}


\begin{figure}[htbp]
\begin{center}
\psfrag{i0}{$i_0$} \psfrag{i1}{$i_1$} \psfrag{j0}{$j_0$}
\psfrag{j1}{$j_1$} \psfrag{e}{$\varepsilon$}
\psfrag{s}{$\Su(\pi,\lambda)$} \psfrag{g}{$\gamma(\pi)$}
\psfrag{et}{$\eta$}

\includegraphics[width=13cm]{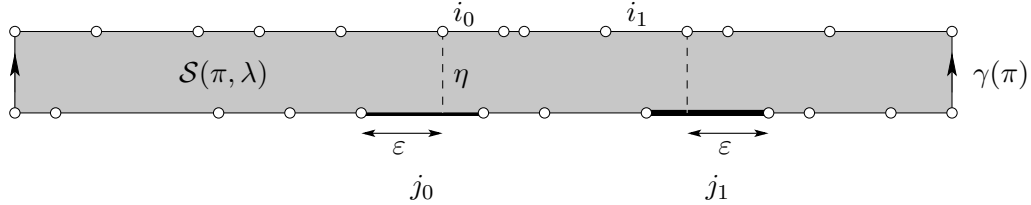}
\end{center}
\caption{ \label{fig:irred2:proof} A separatrix of length $2$
(here the canonical separatrix $\gamma(\pi)$ has length $1$). The
two corresponding horizontal intervals of length $\lambda_{j_0+r}$
and $\lambda_{j_1+r}$, numbered by $j_0$ and $j_1$, are glued by a
central symmetry.}
\end{figure}

Euclidean elementary geometry on the half-translation surface
$\Su$ on Figure~\ref{fig:irred2:proof:idoc} gives
\begin{equation*}
\sum_{i=1}^{i_0} \lambda_i = \sum_{j=1}^{j_0} \lambda_{j+r} -
\varepsilon
\qquad \textrm{and} \qquad
\sum_{i=1}^{i_1} \lambda_i - \varepsilon = \sum_{j=1}^{j_1}
\lambda_{j+r}
\end{equation*}
Adding these two formulas and recalling that
$\lambda_{i_1}=\lambda_{j_0+r}$, we get
\begin{equation}
\label{eq:irred2:idoc}
2 \sum_{i=1}^{i_0} \lambda_i + \sum_{i=i_0+1}^{i_1-1} \lambda_i = 2
\sum_{j=1}^{j_0-1} \lambda_{j+r} + \sum_{j=j_0}^{j_1}
\lambda_{j+r}
\end{equation}

The same argument on Figure~\ref{fig:irred2:proof} produces similar
equalities. More precisely
\begin{equation*}
\sum_{i=1}^{i_0-1} \lambda_i = \sum_{j=1}^{j_0-1} \lambda_{j+r} +
\varepsilon
\qquad \textrm{and} \qquad
\sum_{i=1}^{i_1} \lambda_i = \sum_{j=1}^{j_1} \lambda_{j+r} -
\varepsilon
\end{equation*}
Adding these two formulas and recalling that $\lambda_{j_0+r}=\lambda_{j_1+r}$, we get
the new one:
\begin{equation}
\label{eq:irred2}
2 \sum_{i=1}^{i_0-1} \lambda_i + \sum_{i=i_0}^{i_1} \lambda_i = 2
\sum_{j=1}^{j_0-1} \lambda_{j+r} + \sum_{j=j_0+1}^{j_1-1}
\lambda_{j+r} + 2 \lambda_{j_0+r}
\end{equation}
Recall that by assumption on $\lambda$, there exists exactly
one rational relation between $\lambda_i$. This is the following one.
\begin{equation}
\label{eq:admissible:vect}
\sum_{i=1}^{r} \lambda_i = \sum_{j=1}^l \lambda_{r+j}
\end{equation}
It is easy to check that Equation~(\ref{eq:irred2:idoc})
can not occur: on the right part there is a term $\lambda_{j_0+r}$ and
on the left part, the corresponding term does not appear.

Therefore only the second case arise. Comparing the coefficients
of $\lambda_i$, and forcing the terms to cancel (using
Equation~(\ref{eq:admissible:vect})) this leads to the fact that
$\pi$ does not satisfy $\Red$; with corresponding sets
$$
\begin{array}{lll}
Y'_1=(1,\dots,i_0)   & Y''_1=(i_0+1,\dots,i_1-1) & Y'''_1=(i_1,\dots,r) \\
Y'_2=(1,\dots,j_0-1) & Y''_2=(j_0+1,\dots,j_1-1) &
Y'''_2=(j_1+1,\dots,l)
\end{array}
$$
Therefore if $\lambda$ belongs to $E$ and $\pi$ satisfies $\Red$,
any vertical separatrix $\eta$ has length different from $2$.
Proposition~\ref{prop:length:2} is proven.
\end{proof}

\subsubsection{Irreducibility}

\begin{Definition}
We say that $\pi$ is {\it irreducible} if $\pi$ is weakly
irreducible and satisfies the condition $\Red$.
\end{Definition}

\begin{Proposition}
\label{prop:geo:permu}
Let us consider a surface $\Su(\pi,\lambda)$ with an irreducible
permutation $\pi$. Then there exists a full Lebesgue measure set of
admissible vectors
$\lambda$ such that any vertical separatrix $\eta \subset \Su(\pi,\lambda)$,
different from $\gamma(\pi)$, satisfies $|\eta| \geq 3\cdot |\gamma(\pi)|$.
\end{Proposition}

\begin{proof}[Proof of Proposition~\ref{prop:geo:permu}]
The proof is obvious using Lemma~\ref{rk:length:1} and
Proposition~\ref{prop:length:2}: the length of any separatrix is a
positive integer different from $1$ and $2$ for a full Lebesgue
measure set.
\end{proof}

\begin{Remark}
Proposition~\ref{prop:geo:permu} is related to Keane's i.d.o.c
property. More precisely, For ``true'' permutations, weak
irreducibility coincides with the 
classical definition. This means that $\pi\{1,\dots,k\} \not =
\{1,\dots,k\}$ for any $k=1,\dots,r-1$. The Keane's property
asserts that, for irreducible permutations, for almost all
$\lambda$, any separatrix $\eta \not = \gamma(\pi)$ has infinite
length. It will be interesting to have an analogous result for
generalized permutations. 
\end{Remark}
\subsection{Irreducibility and weak irreducibility}
\label{sec:tech:cond}

\begin{Definition}
Let $\pi$ be a type $(r,l)$ generalized permutation. 
We say that $\pi$ satisfies the Condition~$(*)$ if there exists
only one element $i_0 \leq r$ (respectively $j_0 \geq r+1$) such
that $\pi(i_0) \leq r$ (respectively $\pi(j_0) \geq r+1$).
\end{Definition}

We end this section by the following obvious lemma.

\begin{Lemma}
\label{lm:weirr:irr}
Under condition $(*)$, weak irreducibility implies irreducibility.
\end{Lemma}

\subsection{Irreducibility and ``Breaking up a singularity''}
\label{sec:fundamental:observations:1}

Let $\pi$ be a generalized permutation and $\lambda$ any
admissible vector. Recall (Section~\ref{sec:generalized:permutations})
that for a horizontal separatrix $\beta$, we denote
by $\beta^1$ and $\beta^2$ the two corresponding intervals on the
cylinder $\CYL$.

\begin{Proposition}
\label{prop:fundamental:1}
Let $\pi$ be a generalized permutation and $\Su=\Su(\pi,\lambda)$ the 
suspended flat surface associated to $\lambda$. If $\pi$ is
irreducible then there exists a full Lebesgue measure set of $\lambda$
such that $\gamma(\pi)\in \Su(\pi,\lambda)$ has multiplicity one.
\end{Proposition}

\begin{proof}
It follows from Proposition~\ref{prop:geo:permu} and properties of
\^homologous separatrix detailed in
Proposition~\ref{prop:saddle:multiplicity}.
\end{proof}

\begin{Proposition}
\label{prop:fundamental:2}
Let $\pi$ be a generalized permutation and $\Su=\Su(\pi,\lambda)$ the 
suspended flat surface associated to $\lambda$. We denote by $\beta$
and $\eta$ any two horizontal separatrices. Assume that one of the two 
following holds.
\begin{itemize}
\item $\beta^1$ and $\beta^2$ are located in two {\it different}
horizontal boundaries of the cylinder $\CYL$.

\item all intervals $\beta^1,\ \beta^2$ and $\eta^1,\ \eta^2$ are
located in the {\it same} horizontal boundary of $\CYL$.
\end{itemize}
Then there exists a full Lebesgue measure set of $\lambda$ such
that $\beta$ has multiplicity one.
\end{Proposition}

\begin{proof}[Proof of Proposition~\ref{prop:fundamental:2}]
In the first case there is no condition on the horizontal parameter
$|\beta^1| = |\beta^2|=|\beta|$ (see Equation~(\ref{eq:admissible:vect}).
So one can and do choose the length
of $\beta$ in the flat metric, arbitrary small with respect to the
other length of horizontal separatrices. Therefore
Proposition~\ref{prop:saddle:multiplicity} applies.

In the second case there is only one
linear relation on the length of $\gamma$; This is 
Equation~(\ref{eq:rationnal}). Therefore one gets:
$$
|\beta^1| + |\beta^2| + |\eta^1| + |\eta^2| + \dots = \dots
$$
In the left part of this equality, the terms $\beta^i,\ \eta^i$
survives. In particular, we can choose $|\beta^1| + |\beta^2| +
|\eta^1| + |\eta^2|$ arbitrary small and hence $|\beta^1| =
|\beta^2|=|\beta|$ arbitrary small with respect to the length of
other horizontal separatrices. We are done.
\end{proof}

As direct corollary, one gets the following useful result.

\begin{Theorem}
\label{theo:fundamental:1}
Let us assume that either $\gamma(\pi)$ or $\beta$ (in above
propositions) is a saddle connection connecting two different
singularities; not two poles. Let us also assume that one of the
assumptions of above Proposition~\ref{prop:fundamental:2} 
holds. Then the surface
$\Su(\pi,\lambda)$ is obtained from the surgery ``breaking up a
singularity'' on a surface in a lower dimensional stratum for a
full Lebesgue measure set of $\lambda$.
\end{Theorem}

\subsection{Irreducibility and ``Bubbling a handle''}
\label{sec:fundamental:observations:2}

\begin{figure}
\begin{center}
\psfrag{p}{$P$} \psfrag{gp}{$\gamma(\pi)$} \psfrag{b}{$\beta$}
\psfrag{e1}{$\eta^1$} \psfrag{e2}{$\eta^2$} \psfrag{4p}{$4\pi$}

\includegraphics[width=14cm]{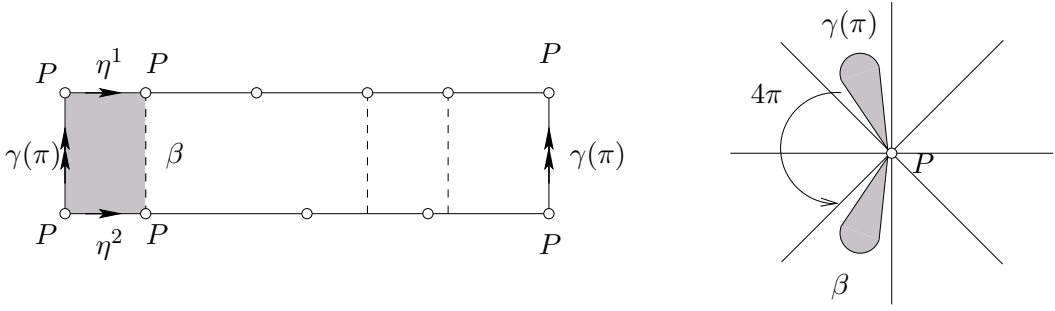}
\end{center}
\caption{ \label{fig:simple:cyl} On the left, the figure
represents a flat surface $\Su(\pi,\lambda)$ with $\pi$ satisfying
$\pi(1)=r+1$. In the vertical direction, one can easily see that
there is a simple cylinder. The boundary component of this
cylinder is $\gamma(\pi) \sqcup \beta$. On the figure presented on the
right, we have represented the diagram of the vertical foliation
of $\Dq=\D z^2$. In this example, the angle of this cylinder is
$4\pi$ (or $6\pi$ is we consider the complementary angle). If the
generalized permutation $\hat \pi$ is irreducible, one can choose
horizontal parameters so that $\gamma(\pi)$ has
multiplicity one. In this case this surface is obtained from a
surface in genus $g-1$, where $g=$genus($\Su$), by ``bubbling a
handle''.
}
\end{figure}

\begin{Notation}
Let $\pi$ be a generalized permutation of the set
$\{1,\dots,r+l\}$ satisfying the condition $\pi(1)=r+1$. We will 
denote by $\hat \pi$ the {\it restricted} generalized
permutation of the set
$\{\widehat{1},2\dots,r,\widehat{r+1},r+2,r+l\}$ (where $\hat i$
means that we forgot the element $i$). In terms of table, this gives
$$
\pi = \left(
\begin{array}{cc}
0 & A \\
0 & B
\end{array} \right) \qquad \textrm{and} \qquad
\hat \pi = \left(
\begin{array}{c}
A \\
B
\end{array} \right)
$$
\end{Notation}
Clearly, the surface $\Su(\pi,\lambda)$, with $\pi$ as above,
possesses a simple cylinder in the vertical direction (see
Figure~\ref{fig:simple:cyl} and Section~\ref{sec:sub:tools}). This
cylinder is filled by regular vertical closed geodesics; each
boundary component is a single vertical separatrix.

Using Theorem~\ref{theo:fundamental:1}, one deduces the following
useful result.

\begin{Theorem}
\label{theo:fundamental:2}
Let $\Su(\pi,\lambda) \in \Q_g(4g-4)$ be a point. Let us
assume that $\pi(1)=r+1$. If $\hat \pi$ is irreducible then
the surface $\Su$ is obtained from the surgery ``bubbling a handle''
on a surface in a the stratum $\Q_{g-1}(4g-8)$ for a full Lebesgue
measure set of $\lambda$.
\end{Theorem}

\begin{proof}[Proof of Theorem~\ref{theo:fundamental:2}]
In the vertical direction, the surface $\Su$ has a simple
cylinder, for any $\lambda$. Let us remove it. We obtain a
half-translation surface with boundaries. Each boundary component
is a single geodesic circle: $\gamma(\pi)$ and one other say
$\beta$. By construction, they have the same length. Let us remove
the singularity and glue these two geodesics {\it segments}
together. We obtain a closed half-translation surface $\Su'$ of
genus $g-1$. The induced quadratic differential has two
singularities say $P_1$ and $P_2$ of multiplicities $k_1,\ k_2$
depending of the conical angle between $\gamma(\pi)$ and $\beta$.
We denote this angle by $k\pi$. By definition
$$
k_1=k-2 \qquad \textrm{and} \qquad k_2=4g-6-k.
$$
By assumptions $\hat \pi$ is irreducible. Thus
Theorem~\ref{theo:fundamental:1} implies that one can choose $\lambda$
in such a way that $\gamma(\pi)$ has multiplicity one. Applying
Theorem~\ref{theo:collapse:singu} we can collapse $\gamma(\pi)$ to
a point. Therefore we obtain a closed flat surface $\Su''$ of
genus $g-1$ with a unique singularity. By construction, ``bubbling
a handle'' at the unique zero of $\Su''$ with the appropriate angle $k\pi$,
we get the surface $\Su$. The theorem is proven.
\end{proof}

We end this section with the following lemma.

\begin{Lemma}
\label{lm:cc:irred:oplus}
The three following statements hold.

\begin{enumerate}

\item
\label{lm:pt:1}
 The component $\Q^{irr}(-1,3,6)$ is adjacent to the component 
$\Q^{irr}(-1,9)$.

\item
\label{lm:pt:2}
The component $\Q^{irr}(-1,3,3,3)$ is adjacent to the component 
$\Q^{irr}(-1,3,6)$.

\item
\label{lm:pt:3}
$\Q^{irr,I}(12) = \Q(8) \oplus 2$ and $\Q^{irr,II}(12) = \Q(8) \oplus 6$.

\item
\label{lm:pt:4}
$\Q^{irr}(-1,9) = \Q(-1,5) \oplus 3$.

\end{enumerate}
\end{Lemma}

\begin{proof}[Proof of the Lemma~\ref{lm:cc:irred:oplus}]
We first prove the two first assertions. Let $\Su(\pi,\lambda) \in 
\Q^{irr}(-1,3,6)$ be a point with $\pi$ as in 
Definition~\ref{def:irr}. Observing the horizontal foliation, it is easy 
to see that there exists a multiplicity one saddle connection connecting the 
two zeroes. Then the first assertion deduces from
Proposition~\ref{theo:collapse:singu}. One proves the second assertion
in the same manner.

\noindent Now let us concentrate with the third assertion. Let $\Su_1=\Su(\pi_1,\lambda),
\Su_2=\Su(\pi_2,\lambda)$ be two representative elements of 
$\Q^{irr,I}(12),\Q^{irr,II}(12)$ with $\pi_1,\pi_2$ as in 
Definition~\ref{def:irr} and $\lambda=(1,\dots,1)$. \\
We detail the first case, that is $\Su_1$; the second case is similar. 
The vertical foliation on $\Su_1$ is completely periodic and decomposes the 
surface into two cylinders (see Figure~\ref{fig:q12irr}). 
\begin{figure}[htbp]
\begin{center}
\psfrag{1}{$1$} \psfrag{5}{$5$} 
\psfrag{2}{$2$} \psfrag{6}{$6$} 
\psfrag{3}{$3$} \psfrag{7}{$7$} 
\psfrag{4}{$4$} 
\includegraphics[width=8cm]{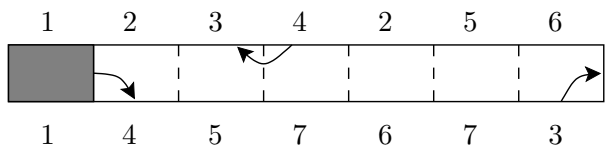}
\end{center}
\caption{ \label{fig:q12irr} 
A vertical decomposition on a half-translation surface representing the 
component $\Q^{irr,I}(12)$.
}
\end{figure}
One of the cylinder is a simple cylinder and the angle between the separatrix 
loops which form the boundary is $2\pi$. Hence $\Su_1 = \Q(8) \oplus 2$.
The proof of $\Su_2 = \Q(8) \oplus 6$ is similar and left to the 
reader. \\
The last statement is a direct verification and left to the reader.
The lemma is proven.
\end{proof}

\noindent We have now all necessary tools to prove Main
Theorem~\ref{theo:main:1} and Main Theorem~\ref{theo:main:2}.

\section{The minimal stratum}
\label{sec:minimal}

This section is devoted to the so-called {\it minimal} stratum
$\Q(4g-4)$ in case $g\geq 3$. We will prove this stratum is connected
except for genus  $4$; for that we will show there exist at most two
connected components. 

\noindent The proof is done by induction on the genus of the surfaces. The 
step of induction is given by Theorem~\ref{theo:simple:cyl}. The
initialization of the induction is reduced to the proof of the
connectedness of the stratum $\Q_5(16)$ which we establish by a
direct argument.

\begin{Theorem}
\label{theo:classi:minimal}
Any connected component of the stratum $\Q(4g-4)$ is described by
the following list.

\begin{itemize}

\item The stratum $\Q(8)$ is connected.

\item The stratum $\Q(12)$ possesses at most two connected components
--- corresponding to $\Q^{irr,I}(12)$ and $\Q^{irr,II}(12)$.

\item Any other stratum $\Q(4g-4)$, in genus $g \geq 5$, is non-empty
  and connected.

\end{itemize}
\end{Theorem}

\begin{Remark}
The stratum $\Q(4g-4)$ is empty for $g\leq 2$
(see~\cite{Masur:Smillie:2}).

\noindent Zorich (see~\cite{Zorich:computation}) has proved that the
stratum $\Q(12)$ is non-connected. The proof uses the Extended Rauzy
classes.
\end{Remark}

\subsection{Step of induction}
\label{sec:step:minimal}

\begin{Theorem}
\label{theo:simple:cyl}
Let $\mathcal C$ be a connected component of the stratum
$\Q(4g-4)$ in genus $g \geq 4$. Then there exist a
half-translation surface $(\Su,\Dq) \in \Q(4g-8)$ and $s\in \mathbb
N^\ast$ such that the surgery ``bubbling a handle'' at the unique
singularity of $\Dq$ in $\Su$ (with discrete parameter $s$ and
arbitrary continuous parameters) produces surfaces belonging to the 
component $\mathcal C$. In other words the map
\begin{gather*}
\oplus \ :\ \pi_0(\Q_{g-1}(4g-8)) \times \mathbb N^\ast \to
\pi_0(\Q_g(4g-4)) \\
(\mathcal C', \ s ) \mapsto \mathcal C' \oplus s
\end{gather*}
is onto for any $g\geq 4$.
\end{Theorem}

\begin{Remark}
Geometrically the previous statement is equivalent to find a
half-translation surface in $\mathcal C$ with a multiplicity one
simple cylinder. In order to do that we will use
Theorem~\ref{theo:fundamental:2}.

\noindent Kontsevich and Zorich have obtained a similar result in the
particular case of Abelian differentials. Their proof uses Rauzy
classes. Here we first give an independent (geometric) proof of their
result. We then give the proof in full generality.
\end{Remark}

Hence Theorem~\ref{theo:simple:cyl} is equivalent to the following.

\begin{Theorem}
Let $\mathcal C$ be a connected component of the stratum
$\Q(4g-4)$ in genus $g \geq 4$. Then there exists 
$\Su(\pi,\lambda) \in \mathcal C$ with $\hat \pi$ irreducible. 
\end{Theorem}

\begin{proof}[Proof of Theorem~\ref{theo:simple:cyl} for
Abelian differentials]
Let $\mathcal C \subset \HA(2g-2)$ be a connected component in genus
$g\geq 2$. In order to find a surface with a multiplicity one simple
cylinder in 
$\mathcal C$ we will use Section~\ref{sec:fundamental:observations:2}.  
Therefore it is sufficient to find a surface
$(\Su,\omega)=\Su(\pi,\lambda) \in \mathcal C$ 
so that $\hat \pi$ is irreducible.

Let $\pi$ be a ``true'' permutation in the symmetric group $S_r$. We
will assume that the surface $\Su(\pi,\lambda)\in \mathcal C$ has no marked point
that is $\pi(i+1) \not = \pi(i) +1$ for all $i=1,\dots,r$ with the
``dummy'' condition $\pi(r+1):=\pi(1)$. We will show there exists a
permutation $\pi_1$ in the class of $\pi$ for the cyclic order with
$\pi_1(1)=1$ and such that $\hat \pi_1$ is irreducible. This will
prove the theorem; Indeed the surface $\Su(\pi_1,\lambda) \in \mathcal
C$.

The differential $\omega$ has a unique singularity thus one can 
assume, using the cyclic order, that $\pi(1)=1$. If $\hat \pi$ is
irreducible then the theorem holds with $\pi_1=\pi$. Otherwise let us
assume that the restricted permutation $\hat \pi$ is reducible.
Then by definition, there exists $2 \leq i_0 < r$ such that
$$
\pi(\{ 2,\dots,i_0\})=\{ 2,\dots,i_0\}.
$$
Let us consider the following new set: $\pi(i_0+1,\dots,r) = (A_1\
r\ A_2)$. With these notations, we have
\begin{multline*}
\left(
\begin{array}{cccccccc}
1 & 2 & \dots & i_0 & i_0+1 & \dots & r  \\
1 & \dots & \dots & \pi(i_0) & A_1 & r & A_2
\end{array}
\right) =
\pi \sim \pi_1=\\
= \left(
\begin{array}{cccccccc}
r & 1 & 2 & \dots & i_0 & i_0+1 & \dots & r-1  \\
r & A_2 & 1 & \pi(2) & \dots & \dots & \pi(i_0) & A_1
\end{array}
\right).
\end{multline*}
It is easy to see that $A_2 \not = \emptyset$: otherwise the
corresponding surface $\Su(\pi_1,\lambda_1)$ will possess a
marked point (see above). In particular $\Su(\pi,\lambda)$ possesses
also a marked point which is a contradiction. 
Thus if $\hat \pi_1$ is reducible, it is easy to see that the
corresponding invariant set (after renumbering elements) $\pi_1(\{
2,\dots,i_0'\})=\{ 2,\dots,i_0'\}$ will satisfy the condition $i_0'
\geq i_0 + 1 > i_0$. The set $\{0,\dots,r \}$ is finite, thus the
theorem follows by repeating finitely many times this
process. Theorem~\ref{theo:simple:cyl} for Abelian differentials is
proven.
\end{proof}

\noindent To clarify the situation, we decompose the proof into several
steps. We first prove the theorem in a weaker version: we add the
additional assumption that there exists a surface
$\Su(\pi,\lambda)\in \mathcal C$ such that $\pi$ satisfies the
condition~$(*)$ (see Section~\ref{sec:tech:cond}). This
corresponds to Proposition~\ref{prop:combi:quadratic}. We then
prove Lemma~\ref{lm:the:1} and Lemma~\ref{lm:the:2} to get 
this additional condition.

\begin{Proposition}
\label{prop:combi:quadratic}
Let $\pi$ be a generalized permutation such that $\Su(\pi,\lambda)
\in \Q(4g-4)$ with $g\geq 4$. Let us assume that $\pi$ satisfies the
condition~$(*)$. Then there exists $\pi_1 \sim \pi$ such that
$\hat \pi_1$ is irreducible.
\end{Proposition}

We first prove the proposition for four particular permutations. 
This corresponds to Lemma~\ref{lm:ex:minimal}.

\begin{Notation}
A generalized permutation $\pi$ is an ordered partition of $X=\{
0,\dots,r+l \}$ into two ordered lists, $X=Y_1 
\sqcup Y_2$. In the present paper we shall always consider only
those generalized permutations, for which each of $Y_1,\ Y_2$
contains at least one entry of multiplicity two. The permutation
satisfies the condition~$(*)$ so that is each set $Y_1,\ Y_2$
contains {\it exactly} one entry of multiplicity two. Up to
re-labeling, one can assume for the next that the two particular
elements are $1\in Y_1$ and $2\in Y_2$. Finally, up to
cyclic order, one can always put $\pi$ into the form $\pi=
\bigl( \begin{smallmatrix} 0 & A & C  \\ 0 & B & D \end{smallmatrix}\bigr)$.

\noindent Note that if $\hat \pi$ is reducible, it involves one of the
two following decomposition cases:

\begin{center}

\begin{enumerate}
\item \label{reduc:cas1}
\begin{tabular}{l}
$\pi(A)=B$ or $\pi(C)=D$
\end{tabular}

\item \label{reduc:cas2}
\begin{tabular}{l}
For all $i \in A$, with $i\not =1$ one have $\pi(i) \in B$; $1\in C$, $\pi(1) \in C$. \\
For all $j \in B$, with $j\not =2$ one have $\pi(j) \in A$;
$2\in D$, $\pi(2) \in D$.
\end{tabular}

\end{enumerate}
\end{center}
In addition, we assume that this decomposition is minimal: we do
not have a decomposition into smaller sets $A',B',C',D'$. This
condition, in case of ``true'' permutation, is equivalent to ask that
$i_0$ is minimal (see the proof versus Abelian differentials). 
\end{Notation}

\begin{Lemma}     
\label{rk:proof:versa}
If $\hat \pi$ involves reducibility of type~(\ref{reduc:cas1}) then
one can find $\pi_1 \sim \pi$ such that $\hat \pi_1$ is irreducible.
\end{Lemma}

\begin{proof}
Imitating the proof of Theorem~\ref{theo:simple:cyl} versus Abelian
differentials, one can easily see that if $\hat \pi$ involves
reducibility of type~(\ref{reduc:cas1}) then one can find $\pi_1
\sim \pi$ with $\hat \pi_1$ irreducible.
\end{proof}

\begin{Lemma}
\label{lm:ex:minimal}
Let $\pi$ be a generalized permutation satisfying the
condition~$(*)$. Assume also that one can put $\pi$ into one
of the four following forms:
$$
\begin{array}{cccc}
\bigl( \begin{smallmatrix} 0 & A & 1 & C & |&3 & 1  \\ 0 & B & 2 & D & |&3 & 2 \end{smallmatrix}\bigr) &
\bigl( \begin{smallmatrix} 0 & A & 1 & C & |&1  \\ 0 & B & 2 & D & |&2 \end{smallmatrix}\bigr) &
\bigl( \begin{smallmatrix} 0 & A & 1 & C & |&3 & 1  \\ 0 & B & 2 & D & |&2 & 3 \end{smallmatrix}\bigr) &
\bigl( \begin{smallmatrix}  0 & A & 1 & C & |&4 &3 & 1 \\ 0 & B & 2 & D & |&4 &2 & 3 \end{smallmatrix}\bigr)
\end{array}
$$
Let us assume that $\hat \pi$ is reducible (with corresponding marked invariant minimal
sets) and $\Su(\pi,\lambda)\in \Q(4g-4)$ with $g\geq 4$. Then there exists
$\pi_1 \sim \pi$ such that $\hat \pi_1$ is irreducible.
\end{Lemma}

\begin{proof}[Proof of Lemma~\ref{lm:ex:minimal}]
Here we address the proof for the first class of permutations; the
others being completely similar. By assumption on minimality of the
decomposition, at least one of the two sets $C$ or $D$ is non-empty. Up to a
permutation of the lines, we assume $C\not = \emptyset$.
Up to re-labeling, let us denote $C:=(C \ 4)$ with $\pi(4) \in B \sqcup D$.
Depending the value of $\pi(4)$, we also put $B:=(B_1 \ 4 \ B_2)$ or
$D:=(D_1\ 4\ D_2)$. Therefore $\pi$ is equivalent to
$$
\begin{array}{ccc}
\pi_B=\bigl( \begin{smallmatrix} 4 & 3 & 1 & 0 & A & 1 & C &  \\
4 & B_2 & 2 & D & 3 & 2 & 0 & B_1  \end{smallmatrix}\bigr) \textrm{ if } \pi(4)\in B
& \textrm{or }
\pi_D=\bigl( \begin{smallmatrix} 4 & 3 & 1 & 0 & A & 1 & C  \\
4 & D_2 & 3 & 2 & 0 & B & 2 & D_1 \end{smallmatrix}\bigr) \textrm{ if } \pi(4)\in D
\end{array}
$$
This involves the two following cases:

\vskip 3mm

\noindent {\bf Case $1$.} If $B_2$ or $D$ are non-empty, the
restricted permutation $\hat \pi_B$ is eventually reducible but in this
case it involves decomposition of type~(\ref{reduc:cas1});
thus the lemma follows from Lemma~\ref{rk:proof:versa}.
Therefore let us assume that $B_2 = D = \emptyset$. If $B_1$ is empty
then the surface $\Su(\pi,\lambda)$ belongs to the stratum
$\Q(8)$ and hence it has genus $3$ which is a contradiction.
Thus up to re-labeling, one can put $B_1:=(B\ 5)$ with $\pi(5) \in A \sqcup C$.
As above depending the value of $\pi(5)$, we also put $A:=(A_1 \ 5 \
A_2)$ or $C:=(C_1\ 5\ C_2)$. Therefore $\pi$ is equivalent to
$$
\begin{array}{ccc}
\pi_A=\bigl( \begin{smallmatrix}5 & A_2 & 1 & C & 4 & 3 & 1 & 0 & A_1   \\
5 & 4 & 2 & 3 & 2 & 0 & B & && \end{smallmatrix}\bigr) \textrm{ if } \pi(5)\in A
& \textrm{or }
\pi_C=\bigl( \begin{smallmatrix} 5 & C_2 & 4 & 3 & 1 & 0 & A & 1 &C_1  \\
5 & 4 & 2 & 3 & 2 & 0 & B & && \end{smallmatrix}\bigr) \textrm{ if } \pi(5)\in C
\end{array}
$$
Now we easily see that each of these permutations $\hat \pi_A$
and $\hat \pi_C$ is eventually reducible but it then involves 
decomposition of Type~(\ref{reduc:cas1}). Thus
Lemma~\ref{rk:proof:versa} applies. 

\vskip 3mm

\noindent {\bf Case $2$.} The discussion of this case is
completely similar to the above one, depending the dichotomy
$\pi(D_2) \subseteq A$ or not. We do not give the complete details
here. The Lemma~\ref{lm:ex:minimal} is proven.
\end{proof}

\begin{proof}[Proof of Proposition~\ref{prop:combi:quadratic}]
Let $\pi$ be a generalized permutation, satisfying condition~$(*)$,
such that $\Su(\pi,\lambda) \in \Q(4g-4)$ with $g\geq 4$.
If $\hat \pi$ is reducible then, according to previous notation and
Lemma~\ref{rk:proof:versa}, one can assume that decomposition
of Type~(\ref{reduc:cas2}) arises. We have
$\pi=\bigl( \begin{smallmatrix} 0 & A & C  \\ 0 & B & D \end{smallmatrix}\bigr)$.
with $1\in A$ and $\pi(1)\in C$, $2\in B$ and $\pi(2)\in D$.
Let us introduce the following new notations: $C = (C_1\ 1\ C_2)$ and
$D = (D_1\ 2\ D_2)$. Then we put $\pi$ into the following form
$\pi=\bigl( \begin{smallmatrix}0 & A &|& C_1 & 1 & C_2  \\
0 & B&| & D_1 & 2 & D_2 \end{smallmatrix}\bigr)$.

\begin{Claim}
Either there exists $\pi_1 \sim \pi$ with $\hat \pi_1$ irreducible
or $\pi$ can be put to one of the two forms
(with $\pi(1)\in C$ and $\pi(2) \in D$):
$$
\begin{array}{ccc}
\bigl( \begin{smallmatrix} 0 & A & C & 1  \\
 0 & B & D & 2 \end{smallmatrix}\bigr) & \textrm{ or } &
\bigl( \begin{smallmatrix} 0 & A & C & 1 & 3  \\
 0 & B & D & 3 & 2 \end{smallmatrix}\bigr)
\end{array}
$$
\end{Claim}

\begin{proof}[Proof of the Claim]
Obviously, if $C_1$ and $C_2$ are empty, the first form arises.
Thus, by symmetry, let us assume $C_2\not = \emptyset$. Up to re-labeling,
we put $C_2:=(C_2 \ 3)$ with $\pi(3)\in D_1 \sqcup D_2$.
As usual, depending the value of $\pi(3)$ we put
$D_1 :=(D'_1 \ 3 \ D''_1)$ or $D_2 := (D'_2 \ 3 \ D''_2)$. Thus, each
case involves the two new permutations in the class of $\pi$:
$$
\begin{array}{ccc}
\pi_{D_1}=\bigl( \begin{smallmatrix} 3 & 0 & A & C_1 & 1 & C_2 &  \\
 3 & D''_1 & 2 & D_2 & 0 & B & D'_1 \end{smallmatrix}\bigr) & \textrm{ or } &
\pi_{D_2}=\bigl( \begin{smallmatrix} 3 & 0 & A & C_1 & 1 & C_2 & \\
 3 & D''_2 & 0 & B & D_1 & 2 & D'_2 \end{smallmatrix}\bigr)
\end{array}
$$
For the first permutation: If $D''_1=D_2=\emptyset$ then the
permutation is of the second class of the lemma thus we are done. 
Otherwise it is easy to see that the restricted permutation
$\hat \pi_{D_1}$ is eventually reducible but then the invariant set is
larger than $A$ and we are done by repeating finitely many times
this process.\medskip

\noindent For the second permutation: We have $D''_2\not = \emptyset$. As above,
the restricted permutation $\hat \pi_{D_2}$ is eventually 
reducible but then the new invariant set is larger than $A$ and we are also done
by repeating finitely many times this process. The claim is proven.
\end{proof}

Now we will consider the two lists $C,\ D$ of the claim.
Each permutation of the claim produces two new classes of permutations.
Here we do not give the details but the algorithm is
completely similar to the one described above. Namely, the following holds:

\begin{Claim}
Either there exists $\pi_1 \sim \pi$ with $\hat \pi_1$ irreducible
or $\pi$ can be put to one of the four forms:
$$
\begin{array}{cccc}
\bigl( \begin{smallmatrix} 0 & A & 1  \\
0 & B & 2  \end{smallmatrix}\bigr) &
\bigl( \begin{smallmatrix}0 & A & 3 & 1  \\
 0 & B & 3 & 2 \end{smallmatrix}\bigr) &
\bigl( \begin{smallmatrix} 0 & A & 3 & 1  \\
0 & B & 2 & 3  \end{smallmatrix}\bigr) &
\bigl( \begin{smallmatrix} 0 & A & 4 & 3 & 1 \\
0 & B & 4 & 2 & 3  \end{smallmatrix}\bigr)
\end{array}
$$
\end{Claim}

Now, let us remark that each of the permutations of the previous claim
falls in the list of Lemma~\ref{lm:ex:minimal}, and for those, we have already
proved the proposition. Proposition~\ref{prop:combi:quadratic} is
proven.
\end{proof}

\begin{Lemma}
\label{lm:the:1}
Let $\mathcal C_0$ be a connected component of the minimal stratum
$\Q(4g-4)$. Then there exist two sequences of connected components
\begin{itemize}

\item $\mathcal C_i \subset \Q(4g-4)$ for $i=0,\dots,m$ ($m \geq 1$)

\item $\mathcal C^j \subset \Q(k_j,4g-4-k_j)$ for $j=1,\dots,m$,
$k_j \geq 1$

\end{itemize}
and a flat surface $\Su=\Su(\pi,\lambda) \in \mathcal C_m$ such that:
$$
\mathcal C_i \cup \mathcal C_{i+1} \subset \overline{\mathcal
C^{i+1}},\ \forall i=0,\dots,m-1 \qquad \textrm{ and } \pi
\textrm{ satisfies the condition } (*).
$$
\end{Lemma}

\begin{proof}[Proof of Lemma~\ref{lm:the:1}]
Let $\Su(\pi,\lambda) \in \mathcal C_0$ be a point. If $\pi$
satisfies the condition $(*)$ we are done. Otherwise let us denote
$\eta_1,\dots,\eta_m$ the set of horizontal separatrices such that
$\eta_i^1$ and $\eta_i^2$ belong inside a same boundary component of
the cylinder Cyl$(\Su)$. By assumption $m\geq 3$. Let us 
``break up'' the unique zero into two zeroes to obtain a new surface,
say $\Su'_1$, which belong to a component $\mathcal C^1$. By
Proposition~\ref{prop:fundamental:2} and the fact that $m\geq 3$,
the saddle connection $\eta_1$ on $\Su'_1$ has multiplicity one so
one can it to a point. Thanks to this process we get a new
surface $\Su_1$ in the minimal stratum, which belong to a
component $\mathcal C_1$ (eventually different from $\mathcal C_0$).
By construction we have
$$
C_0 \cup C_1 \subset \overline{\mathcal C^{1}}
$$
Repeating inductively this process on saddle connection $\eta_i$
for $i=2,\dots,m-1$ we obtain the following diagram
$$
\xymatrix{
\mathcal C^1 \ar@{-}[d] \ar@{-}[rd] & \mathcal C^2 \ar@{-}[d] \ar@{-}[rd] & \dots &
\mathcal C^{i+1} \ar@{-}[d] \ar@{-}[rd] & \dots & \dots  &
\mathcal C^m \ar@{-}[d] \ar@{-}[rd]  \\
\mathcal C_0  & \mathcal C_1  & \dots &
\mathcal C_i  & \mathcal C_{i+1}  & \dots &
\mathcal C_{m-1}  & \mathcal C_m
}
$$
with $\mathcal C_i \cup \mathcal C_{i+1} \subset
\overline{\mathcal C^{i+1}}$.

Finally, on the surface $\Su_m(\pi_m,\lambda_m) \in \mathcal C_m$, by
construction there exists a single pair of horizontal separatrices
$\eta_m$ and $\eta_{m-1}$ such that $\eta_i^1$ and $\eta_i^2$
belong to the same boundary component of the cylinder Cyl$(\Su)$.
In other word $\pi_m$ satisfies the
condition~$(*)$. Lemma~\ref{lm:the:1} is proven.
\end{proof}

\begin{Lemma}
\label{lm:the:2}
Let $\mathcal C_0 \subset \Q(4g-4)$ be a connected component.
Let us assume there exist two components $\mathcal C_1 \subset \Q(4g-4)$ and
$\mathcal C^1 \subset \Q(k,4g-4-k)$ such that
$$
\mathcal C_0 \cup \mathcal C_1 \subset \overline{\mathcal C^1}
$$
Let us also assume there exists a flat surface $(\Su,\Dq) \in
\mathcal C_1$ with a multiplicity one simple cylinder. Then there
also exists a flat surface $(\Su',\Dq') \in \mathcal C_0$ with a
multiplicity one simple cylinder.
\end{Lemma}

\begin{proof}[Proof of Lemma~\ref{lm:the:2}]
The proof is obvious using description of surfaces in terms of
separatrices diagrams.
\end{proof}

\begin{Corollary}
\label{cor:condition}
For each connected component $\mathcal C$ of $\Q(4g-4)$ with $g\geq 4$,
there exists a surface $\Su(\pi,\lambda)\in \mathcal C$ such that $\pi$ satisfies
the condition~$(*)$.
\end{Corollary}

\begin{proof}
It follows from Lemma~\ref{lm:the:1} and Lemma~\ref{lm:the:2}.
\end{proof}

\begin{proof}[Proof of Theorem~\ref{theo:simple:cyl}]
It follows obviously from Corollary~\ref{cor:condition} and 
Proposition~\ref{prop:combi:quadratic}.
\end{proof}

\subsection{Connectedness of the minimal stratum}

We are ready to prove Theorem~\ref{theo:classi:minimal}.
First we directly show that $\Q(4g-4)$ is connected for
$g=3,5$ and has at most two components for $g=4$. We then prove the
theorem inductively on $g$.

\begin{Lemma}
\label{q8:connexe}
\textrm{ }

The stratum $\Q(8)$ is connected.

The stratum $\Q(12)$ has at most two connected components.
\end{Lemma}

\begin{proof}[Proof of Lemma~\ref{q8:connexe}]
See Lemma~\ref{lm:appendix:q8} and Lemma~\ref{lm:appendix:q12} in
the appendix.
\end{proof}

\begin{Lemma}
\label{q16:connexe}
The stratum $\Q(16)$ is connected.
\end{Lemma}

\begin{proof}[Proof of Lemma~\ref{q16:connexe}]
Thanks to previous lemma let $\mathcal C_0=\Q(8)$ be the
unique connected component of the stratum $\Q(8)$. Let $\mathcal C_1 =
\mathcal C_0 \oplus 2 \oplus 2$ be the component of $\Q_5(16)$ obtaining by 
bubbling two handles on a surface of $\mathcal C_0$ with an angle of $2\pi$.
Recall that
$$
\mathcal C_0 \oplus 2 = \Q^{irr,I}(12)  \qquad \textrm{and} \qquad
\mathcal C_0 \oplus 6 = \Q^{irr,II}(12)
$$
Let $\mathcal C'$ be any component of $\Q(16)$. By
Theorem~\ref{theo:simple:cyl}, there exists $s_0$ such that
$\mathcal C' = \mathcal C \oplus s_0$ where $\mathcal C$ is a
component of $\Q(12)=\Q^{irr,I}(12) \cup \Q^{irr,II}(12)$. Using the
properties of the map $\oplus$
(Proposition~\ref{prop:properties:oplus}), we directly get the relations
$$
\mathcal C_0 \oplus 2 \oplus s = \mathcal C_0 \oplus 6 \oplus s =
\mathcal C_1 \textrm{ for any } s=1,\dots,8.
$$
Thus $\mathcal C'=\mathcal C \oplus s_0 = \mathcal C_1$ which
proves the lemma.
\end{proof}

\begin{proof}[Proof of Theorem~\ref{theo:classi:minimal}]
We process by induction, initialization being given by
Lemma~\ref{q16:connexe}. Let us fix $g > 5$. Let us assume that
$\Q(4g'-4)$ is connected for all genera $5 \leq g' < g$. We have
to show that $\Q(4g-4)$ is connected. We denote by $\mathcal
C_{g-1}= \Q(4(g-1)-4)$ the {\it unique} connected component of
this stratum. We also define $\mathcal C_g \subseteq \Q(4g-4)$ by
\begin{equation*}
\mathcal C_g = \mathcal C_{g-1} \oplus 1
\end{equation*}
By Theorem~\ref{theo:simple:cyl} the map
\begin{gather*}
\oplus \ :\ \pi_0(\Q_{g-1}(4g-8)) \times \mathbb N^\ast \to
\pi_0(\Q_g(4g-4)) \\
(\mathcal C', \ s ) \mapsto \mathcal C' \oplus s
\end{gather*}
is onto. But the stratum $\Q(4(g-1)-4)$ is connected and $s$ can be chosen in
$\{1,\dots,2g-2\}$ (up to consider the complementary angle). Thus we obtain
a (onto) map
\begin{gather*}
\oplus \ :\mathbb \{1,\dots,2g-2\} \to \pi_0(\Q_g(4g-4)) \\
s \mapsto \mathcal C_{g-1} \oplus s
\end{gather*}
In order to end the proof, it remains to show that
\begin{equation}
\label{eq:desire}
\mathcal C_g = \mathcal C_{g-1} \oplus s \qquad \textrm{for any s}.
\end{equation}
Now always by Theorem~\ref{theo:simple:cyl}, there exists $r_0$ such that
\begin{equation*}
\mathcal C_{g-1} = \mathcal C_{g-2} \oplus r_0 \qquad \textrm{with } \mathcal C_{g-2}
\subseteq \Q(4(g-2)-8)
\end{equation*}
(this stratum is non-empty because $g-2 \geq 4$). But recalling
that the stratum $\Q(4(g-1)-4)$ ($g-1 \geq 5$) is connected, we
also have
$$
\mathcal C_{g-2} \oplus s = \mathcal C_{g-1} \qquad \textrm{for
any } s
$$
Using properties of the map $\oplus$, this yields 
$$
\mathcal C_g = \mathcal C_{g-1} \oplus 1 = (\mathcal C_{g-2}
\oplus s) \oplus 1 = (\mathcal C_{g-2} \oplus 1) \oplus s =
\mathcal C_{g-1} \oplus s \qquad \textrm{for any } s
$$
Thus we get the desired relation~(\ref{eq:desire}).
Theorem~\ref{theo:classi:minimal} is proven.
\end{proof}

\section{Adjacency of the strata}
\label{sec:adjacency}

\begin{Theorem}
\label{theo:adja}
Let $\mathcal C \subset \Q_g(k_1,\dots,k_n)$ be a connected
component with any $g \geq 3$, $n \geq 2$. Let us make the additional
assumption that $\mathcal C$ is neither an 
hyperelliptic component nor an irreducible component. Then there
exists a connected component $\mathcal C_0 \subset \Q(4g-4)$ 
such that $\mathcal C_0 \subset \overline{\mathcal C}$.
\end{Theorem}

We will deduce Theorem~\ref{theo:adja} from the following one.

\begin{Theorem}
\label{theo:anexe:1}
Let $\mathcal C \subset \Q_g(k_1,\dots,k_n)$ be a connected component
with $g \geq 1$ and $n \geq 2$. Let us assume that $\mathcal C$ is
none of one of the following:
\begin{itemize}
\item the hyperelliptic connected component $\Q^{hyp}(4(g-k) -
6,4k+2)$ or $\Q^{hyp}(-1,-1,4g-2)$
\item a component of $\Q(-1,5)$
\item the irreducible component $\Q^{irr}(-1,9)$
\end{itemize}
Then there exists a flat surface $(\Su,\Dq) \in \Q_g$ in a lower
dimensional stratum and a surgery ``breaking up a singularity''
at a singularity of $\Dq$ on $\Su$ such that the resulting surface
$(\Su',\Dq')$ belongs to~$\mathcal C$.
\end{Theorem}

\begin{Remark}
Surprisingly, the answer is quite difficult and we find that some
components, in small genera, are not hyperelliptic neither adjacent
to the minimal stratum. Note also that the corresponding statement for
Abelian differentials is trivial.
\end{Remark}

\subsection{Link with our main result}
\label{sec:link:main:result}

Following the description of the topology of the strata of the moduli
space in local coordinates (see~\cite{Kontsevich:Zorich}), one has

\begin{Theorem}[Kontsevich]
\label{theo:def:local}
For any $(\Su,\Dq) \in \Q(4g-4)$ with $g\geq 3$ there exists a small
open set $U(\Su,\Dq)$ of $(\Su,\Dq)$ in the whole space $\Q_g$ such that
$$
\Q(k_1,\dots,k_n) \cap U(\Su,\Dq)
$$
is non-empty and connected for any $k_i$ with $\sum k_i=4g-4$.
\end{Theorem}

Combining this result with Theorem~\ref{theo:adja} we get an
upper bound for the number of connected components of any stratum.
Let us denote by $r_g$ the number of connected components of the
minimal stratum $\Q(4g-4)$. Then for any $g\geq 3$, 
Theorem~\ref{theo:adja} and Theorem~\ref{theo:def:local} together
give 
$$
\begin{array}{lll}
1 \ \leq & \# \left\{ \begin{array}{l} \textrm{components
    of a statum which do not contain}  \\
\textrm{a hyperelliptic component neither an exceptional component}
    \end{array} \right\} & \leq \ r_g \\ 
2 \ \leq & \# \{ \textrm{components of statum which contain a
  hyperelliptic component} \} & \leq \ r_g + 1 \\
1 \ \leq  & \# \{ \textrm{components of statum which
  contain an exceptional component} \}& \leq \ r_g + 1
\end{array}
$$
In Section~\ref{sec:minimal} we have proved that $r_g=1$ for $g \geq 5$
which implies main Theorem~\ref{theo:main:1}.

\subsection{Strategy of the proof of Theorem~\ref{theo:anexe:1}}

\noindent Let $\mathcal C$ be a connected component of the stratum
$\Q_g(k_1,\dots,k_n)$ with $g \geq 1$ and $n \geq 2$. 
According to Theorem~\ref{theo:collapse:singu} one has
to construct a surface $(\Su,\Dq) \in \mathcal C$ with a
multiplicity one separatrix between two {\it different
singularities}, not two poles. In order to do that we will use the
criterion given by Proposition~\ref{prop:fundamental:2} and 
Proposition~\ref{prop:fundamental:1}. \medskip

\noindent Let $(\Su(\pi,\lambda),\Dq)$ be a point in $\mathcal C$. Either
Proposition~\ref{prop:fundamental:2} applies and then we have done or 
assumptions of the proposition does not hold. In this last case we get
some restrictions on the combinatorics of $\pi$.   \medskip

\noindent Now let us consider the vertical foliation on
$(\Su(\pi,\lambda),\Dq)$. As above either
Proposition~\ref{prop:fundamental:1} applies and we have done or
assumptions of the proposition again does not hold. In this last case we get
also new restrictions on the combinatorics of $\pi$. In particular a
simple computation shows that the permutation $\pi$ is completely
determined. It corresponds to a ``hyperelliptic'' or ``irreducible''
permutation. These two cases are avoid; indeed the component $\mathcal
C$ is not hyperelliptic neither irreducible. \medskip

\noindent The proof is decomposed into several cases. Recall that $n$ denote the
number of singularities. First we consider the general case $n \geq
4$. We then prove the case $n=3$. Finally we conclude with the
holomorphic case $n=2$ and the meromorphic case $n=2$: the stratum
$\Q(-1,4g-3)$.

\begin{Remark}
The corresponding statement of Theorem~\ref{theo:anexe:1} versus
Abelian differentials is trivial. Indeed
Proposition~\ref{prop:fundamental:2} applies directly (all transitions
functions are given by translations).
\end{Remark}

\subsection{Proof of Theorem~\ref{theo:anexe:1} in case $n \geq 4$}
\label{sec:proof:n>=4}

Let $\mathcal C$ be a connected component and $\Su(\pi,\lambda) \in
\mathcal C$. We denote the horizontal cylinder of $\Su$
by $\CYL$. The boundaries components of $\CYL$ are denoted by
$I$ and $J$. One associates to each saddle connection $\eta$ two 
intervals $\eta^1,\ \eta^2$ on $I \sqcup J$.
\begin{TheoClaim}
\label{claim:1}
There exists a saddle connection $\alpha$ between a {\it zero}
$P_1$ and another singularity $P_2$.
\end{TheoClaim}

\begin{proof}[Proof of the claim]
The assumption $g \geq 1$ implies there exists at least a zero
$P_1$. If there exists a saddle connection attached to this zero we
are done. Otherwise any separatrix emanating from $P_1$ 
is actually a separatrix loop.

Exchanging the role of $I$ and $J$ if necessary, the previous assertion 
means that any saddle connection located in $I$ is a separatrix 
loop emanating from $P_1$ (see Figure~\ref{fig:theo:case4:1}a).

We claim there exists a separatrix attached to a zero, say 
$P'_1 \not = P_1$, located on $J$. 
Otherwise it will mean that any separatrix in $J$ is
attached to a pole. The only possibility to obtain that is presented
in Figure~\ref{fig:theo:case4:1}b. In particular the stratum is 
$\Q(k,-1,-1)$ contradicting the assumption $n\geq 4$.
\begin{figure}[htbp]
\begin{center}

\psfrag{C}{$\scriptstyle \CYL$}
\psfrag{P1}{$\scriptstyle P_1$}
\psfrag{PP}{$\scriptstyle P_1'$}  \psfrag{PPP}{$\scriptstyle P_1''$}
\psfrag{dots}{$\scriptstyle \dots\dots\dots\dots$}
\psfrag{I}{$\scriptstyle I$}
\psfrag{J}{$\scriptstyle J$}
\psfrag{g}{$\scriptstyle \gamma(\pi)$}

  \subfigure[]{\epsfig{figure=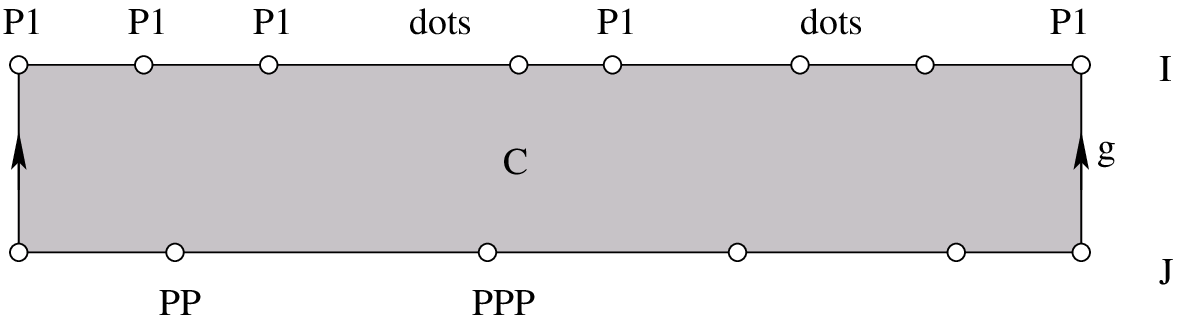,width=7.0cm}} \qquad 
  \subfigure[]{\epsfig{figure=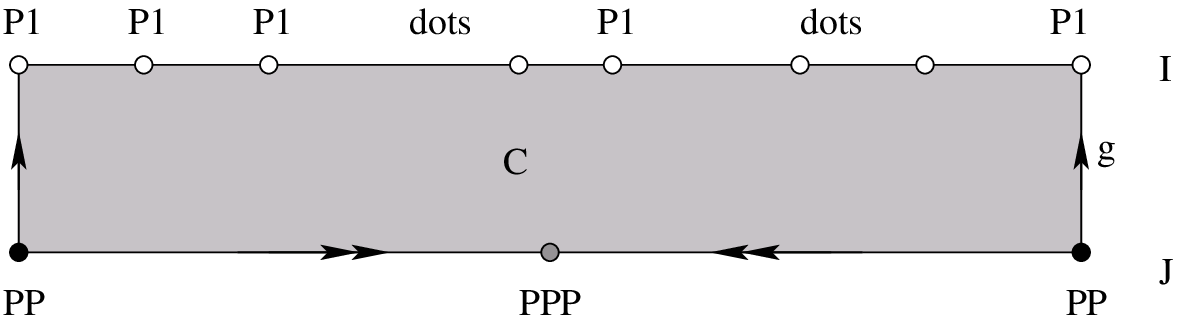,width=7.0cm}}

\end{center}
\caption{
\label{fig:theo:case4:1}
}
\end{figure}
If there exists no saddle connection attached to $P'_1$ then any 
separatrix located on $J$ is a separatrix loop emanating from 
$P'_1$. It implies that the surface has only two zeroes. 
This is again a contradiction and the claim is proven.
\end{proof}

\begin{figure}[htbp]
\begin{center}

\psfrag{C}{$\scriptstyle \CYL$}
\psfrag{P1}{$\scriptstyle P_1$}
\psfrag{P2}{$\scriptstyle P_2$}  \psfrag{P3}{$\scriptstyle P_3$}
\psfrag{dots}{$\scriptstyle \dots\dots\dots\dots$}
\psfrag{I}{$\scriptstyle I$}
\psfrag{J}{$\scriptstyle J$}
\psfrag{g}{$\scriptstyle \gamma(\pi)$}
\psfrag{I1}{$\scriptstyle I_1$} \psfrag{I2}{$\scriptstyle I_2$} 
\psfrag{J1}{$\scriptstyle J_1$} \psfrag{J3}{$\scriptstyle J_3$} 
\psfrag{d}{$\scriptstyle \dots$}
\psfrag{a}{$\scriptstyle \alpha$} \psfrag{b}{$\scriptstyle \beta$}
\psfrag{e}{$\scriptstyle \eta$}

  \subfigure[]{\epsfig{figure=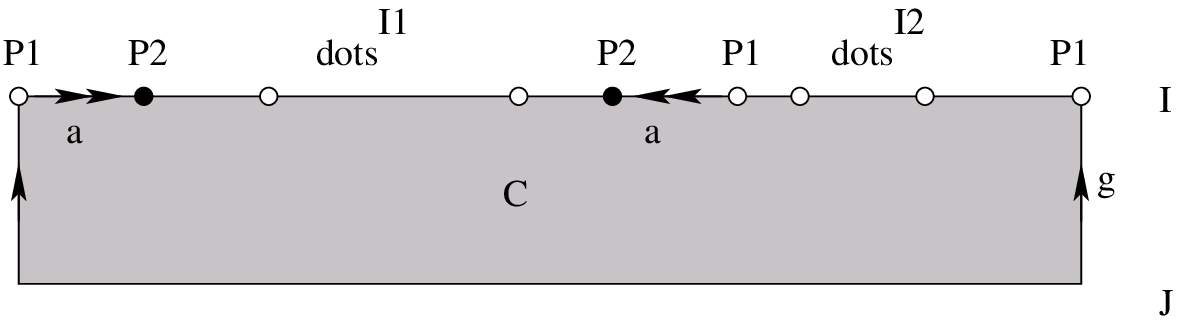,width=7.0cm}} \qquad
  \subfigure[]{\epsfig{figure=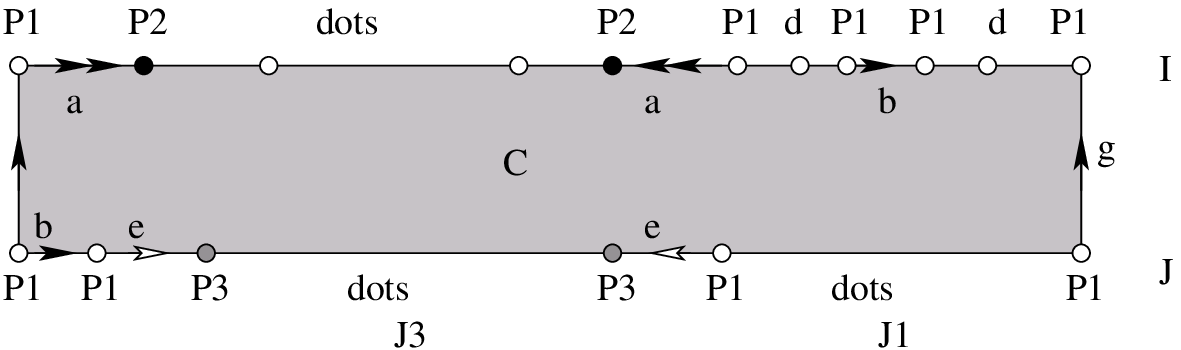,width=7.0cm}}\qquad
  \subfigure[]{\epsfig{figure=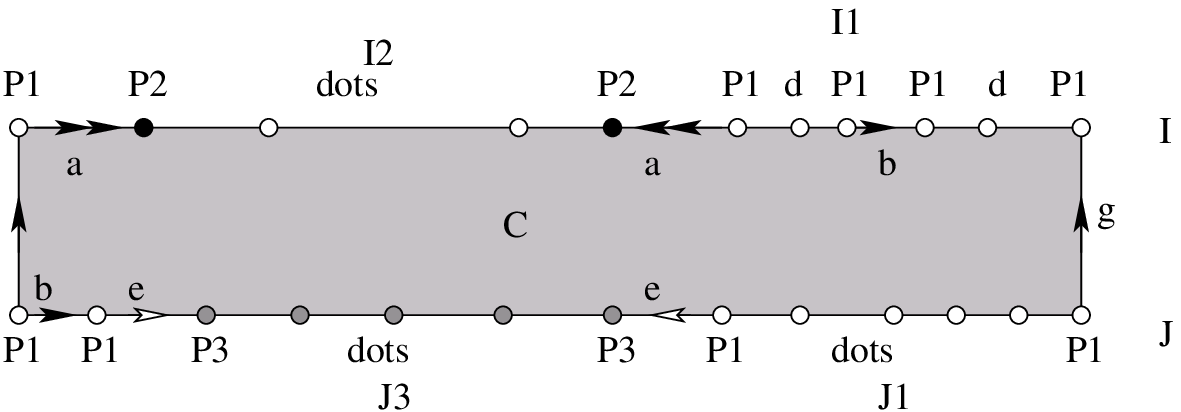,width=7.0cm}}

\end{center}
\caption{
\label{fig:theo:case4:2}
}
\end{figure}

Let $\alpha$ be a saddle connection between the zero $P_1$ and the
singularity $P_2$ (we assume that $\alpha^1 \subset I$). If
$\alpha^2 \subset J$ then Proposition~\ref{prop:fundamental:2} 
implies that $\alpha$ has multiplicity one and the theorem is proven.
Hence let us assume that
$$
\alpha^1,\ \alpha^2 \subset I
$$
(we refer to Figure~\ref{fig:theo:case4:2}a). \\
The complement of $\alpha^1,\alpha^2$ in $I$ has two connected components. 
Following the figure, we will refer to these two components as 
$I_1,I_2$. The interval $I_i$ is a union of separatrices.

\begin{TheoClaim}
\label{claim:2}
Either the theorem is proved or 
any separatrix $\tau$, with $\tau^1 \subset I_i$, is a separatrix loop
attached to $P_i$. In addition $\tau^2 \subset J$ (eventually
$I_2=\emptyset$). 
\end{TheoClaim}

\begin{proof}[Proof of the claim]
Obviously one has: 
Either there exists a saddle connection $\tau$ between two singularities 
(not two poles) with $\tau^1 \subset I_1$ 
or any separatrix of $I_1$ is a separatrix loop attached to $P_1$.

The first case is the conclusion of the theorem. Thus let us 
consider the second case. Let $\tau$ be a separatrix 
with $\tau^1 \subset I_1$. Then $\tau^2 \subset I \sqcup J$. 
Recalling that $\alpha^1,\ \alpha^2 \subset I$ and applying 
Proposition~\ref{prop:fundamental:2} we get that if $\tau^2 \subset I$ 
then $\alpha$ has multiplicity one and therefore the theorem is again proved. 
Otherwise $\tau^2 \subset J$ giving the second conclusion of the claim.

Repeating the same argument for separatrices $\tau$ with 
$\tau^1 \subset I_2$, we get the claim.
\end{proof}

Let us assume that the second conclusion of Claim~\ref{claim:2} holds. 
The singularity $P_1$ is a zero hence there exists a separatrix
loop $\beta$ attached to $P_1$ with $\alpha_1 \subset I_1$. Thanks to 
Claim~\ref{claim:2} $\alpha^2 \subset J$ (see Figure~\ref{fig:theo:case4:2}b). 
In particular there exists a saddle connection $\eta$ between $P_1$ and a singularity 
$P_3 \not = P_1$ (possibly $P_3=P_2$) in such a way that $\eta^1 \subset J$.
If $\eta^2 \subset I$ then thanks to Proposition~\ref{prop:fundamental:2}, 
$\eta$ has multiplicity one and the theorem is proved. 
Thus one can and do assume $\eta^2 \subset J$.
The complement of $\eta^1,\eta^2$ in $J$ has two connected components. 
Following the figure, we will refer to these two components as 
$J_1,J_3$. The interval $J_i$ is a union of separatrices.

\begin{TheoClaim}
\label{claim:3}
Either the theorem is proved or 
any separatrix $\tau$, with $\tau^1 \subset J_i$, is a separatrix loop
attached to $P_i$. In addition $\tau^2 \subset I$.
\end{TheoClaim}

\begin{proof}[Proof of the claim]
The proof parallels the one of previous Claim~\ref{claim:2}.
\end{proof}

Finally thanks the previous claims, either the theorem is proved or 
$\Su$ can be putted into the form prescribed by
Figure~\ref{fig:theo:case4:2}c. The last case implies that $\Su$ has at
most three singularities: $P_1,P_2,P$ contradicting $n\geq 4$. Therefore 
Theorem~\ref{theo:anexe:1} in case $n \geq 4$ is proven.

\subsection{Proof of Theorem~\ref{theo:anexe:1} in case $n=3$}

The proof is similar to the case $n\geq 4$ with a refinement. We first
prove the following lemma.

\begin{Lemma}
\label{lm:permu:2:ele}
Let $\mathcal C \subset \Q_g(k_1,k_2,k_3)$ be a connected component 
with $g \geq 1$. Then there exists $(\Su(\pi,\lambda),\Dq) \in
\mathcal C$ a half-translation surface such that one of the two
following holds.

\begin{itemize}

\item $(\Su(\pi,\lambda),\Dq)$ has a multiplicity one saddle
  connection between two different singularities (not two poles)

\item The singularity pattern is exactly $(k_1,k_2,k_3)=(-1,-1,4g-2)$ and 
the combinatorics of $\pi$ has the form
$$
\pi = \left( \begin{array}{ccccc}
1 & \dots & r & 0_1 & 0_1  \\
0_2 & 0_2 & \pi(1) & \dots & \pi(r)
\end{array} \right)
\qquad \textrm{or} \qquad
\pi = \left( \begin{array}{cc}
A & \\
0 & 0
\end{array} \right) \textrm{ with } \pi(A)=A
$$
\end{itemize}
\end{Lemma}

\begin{proof}[Proof of Lemma~\ref{lm:permu:2:ele}]
Let $\Su(\pi,\lambda)$ be a point in $\mathcal C$. Imitating the proof of
Claim~\ref{claim:1} one gets that either there exists a saddle 
connection $\alpha$ between a zero and another singularity or 
one can put the surface $\Su$ into the form prescribed by 
Figure~\ref{fig:theo:case4:1}b. The last case produces a surface
$\Su(\pi,\lambda)$ belonging to the stratum $\Q(k,-1,-1)$ with $\pi = 
\left( \begin{smallmatrix} A & \\ 0 & 0 \end{smallmatrix} \right)$ 
and $\pi(A)=A$; and so the lemma is proven. \medskip

\noindent Therefore let us assume $\Su$ possesses a saddle 
connection $\alpha$ between a zero and another singularity. Let us
again apply conclusions of Claim~\ref{claim:2} and  
Claim~\ref{claim:3}. This leads to the following dichotomy. 
Either $\Su$ possesses a multiplicity one saddle
connection, and so the lemma is proven, or $\Su$ can put put into the
form prescribed by Figure~\ref{fig:theo:case4:2}c. \medskip

\noindent Finally one can assume the following situation: Any separatrix $\tau$
belonging to $I_1$ or $J_1$ is a separatrix loop attached to $P_1$. 
In addition $\tau^1 \subset I_1$ and $\tau^2 \subset J_1$. Moreover 
$n=3$ thus $P_3 \not = P_2$. \medskip

\noindent If $P_2$ is not a pole then there exists a separatrix loop $\tau$ 
attached to $P_2$ with $\tau^1 \subset I_2$. If $\tau^2 \subset J$ 
then $\tau^2 \subset J_3$ implying $P_2=P_3$ which is
contradiction. Therefore $\tau^2 \subset I$. Now
Proposition~\ref{prop:fundamental:2} implies (indeed
$\alpha^1,\alpha^2 \subset I$) that $\alpha$ has multiplicity one and
so the lemma is proven. \\
The same proof holds if $P_3$ is not a pole. \medskip

\noindent Therefore let us assume that $P_2$ and $P_3$ are two poles. 
Then the stratum containing $\Su$ is $\Q(k,-1,-1)$ and the permutation
$\pi$ has the combinatorics $\left( \begin{smallmatrix} 1 & \dots & r
  & 0_1 & 0_1  \\ 0_2 & 0_2 & \pi(1) & \dots & \pi(r)
\end{smallmatrix} \right)$. The lemma is proven.
\end{proof}

Therefore in order to prove Theorem~\ref{theo:anexe:1} in case $n=3$
we have to analyse the combinatorics of two permutations given in 
Lemma~\ref{lm:permu:2:ele}.

\begin{NoNumberLemma}
Let $\pi=\bigl( \begin{smallmatrix} 2 & \dots & r & 0_1 & 0_1 & 1 \\
0_2 & \pi(1) & \pi(2) & \dots & \pi(r) & 0_2
\end{smallmatrix}\bigr)$. Assume $\pi$ is different (up to the cyclic
order) from $\Pi_1(r,0)$. Then 
there exists $\pi' \sim \pi$ such that $\pi'$ is irreducible and
$\gamma(\pi')$ is a saddle connection between the zero and a pole.
\end{NoNumberLemma}

\begin{proof}[Proof of the lemma]
First of all, let us remark that the separatrix $\gamma(\pi)$
is actually a saddle connection.
Let us also remark that $\pi$ satisfies the condition~$(*)$,
therefore by Lemma~\ref{lm:weirr:irr}, weakly irreducibility implies
irreducibility for any permutation in the class of $\pi$. \medskip

\noindent The above generalized permutation $\pi$ is reducible if and only
if $\pi(r) = 1$. One observes that 
$$
\pi = \bigl( \begin{smallmatrix} 2 & \dots & r & 0_1 & 0_1 & 1 \\
0_2 & \pi(1) & \pi(2) & \dots & \pi(r) & 0_2 \end{smallmatrix}\bigr) \sim \pi'=
\bigl( \begin{smallmatrix} 3 & \dots & r & 0_1 & 0_1 & 1 & 2 \\
0_2 & \pi(1) & \pi(2) & \dots & \pi(r-1) & 1 & 0_2
\end{smallmatrix}\bigr).
$$
It is easy to check that $\pi'$ is irreducible if and only if
$\pi(r-1)=2$.

\noindent Therefore repeating finitely many times this process, one
shows that either there exists an irreducible permutation $\pi'$ in
the class of $\pi$ or
$$
\pi(i) = r - i + 1 \qquad \textrm{ for all } i=1,\dots r.
$$
The last equations mean that $\pi \sim \Pi_1(r,0)$ which is a
contradiction. The lemma is proven.
\end{proof}

\begin{NoNumberLemma}
Let $\pi = \left( \begin{smallmatrix} A & \\ 0 & 0 \end{smallmatrix}
\right)$ with $\pi(A)=A$. Assume $\pi$ is different (up to the cyclic
order) from $\Pi_2(r,1)$. Then there exists $\pi' \sim \pi$ such
that $\pi'$ is irreducible and $\gamma(\pi')$ is a saddle connection
between the zero and a pole.
\end{NoNumberLemma}

\begin{proof}
It is easy to see that for each permutation $\pi'$ in the class of $\pi$,
the separatrix $\gamma(\pi')$ is actually a saddle connection between the
zero and a pole. If $\pi$ is weakly reducible then by definition, we have
$\pi=\bigl( \begin{smallmatrix} B&C \\ 0&0 \end{smallmatrix}\bigr)$
with $\pi(B)=C$. Let us decompose the two lists $B,\ C$ into the following way:
$B = (1\ B_2)$ and $C = (C_1\ 1\ C_2)$. With these notations,
$\pi\sim\bigl( \begin{smallmatrix} 1 & C_2 & 1 & B_2 & C_1 & \\ 0 & & & & 0
\end{smallmatrix}\bigr)$. This permutation is weakly reducible if and only if $C_1 =
\emptyset$. \medskip

\noindent One can repeat finitely many times this process, with
$B:=B_2$ and $C:=C_2$ for the new sets, to get the following
statement. Either there exists a weakly irreducible permutation $\pi'
\sim \pi$ or $B=C=(1\ 2\ \dots\ r)$. The last case means that 
$\pi=\Pi_2(r,1)$ which is a contradiction. Therefore we have proved
that there exists a weakly irreducible $\pi' \sim \pi$. Form this
permutation it is easy to check that there exists $\pi'' \sim \pi'$
irreducible. Hence $\pi'' \sim \pi$ and the lemma is proven.
\end{proof}

Now we are ready to prove Theorem~\ref{theo:anexe:1} in case $n=3$. 
Let $\mathcal C$ be any component of any stratum
$\Q_g(k_1,k_2,k_3)$. Applying 
Lemma~\ref{lm:permu:2:ele} combining with the two above lemma 
one gets the following dichotomy. There exits $\Su \in \mathcal C$
such that either $\Su$ possesses a multiplicity one saddle connection
between two singularities (not two poles), and so the theorem is
proved, or $\Su=\Su(\pi,\lambda)$ with $\pi = \Pi_1(r,0)$ or $\pi =
\Pi_2(r,1)$. These last two permutations correspond to hyperelliptic
permutations and in particular it implies that $\mathcal C$ is itself a
hyperelliptic component of $\Q(-1,-1,4g-2)$ which is a
contradiction. The theorem is proven.

\subsection{Proof of Theorem~\ref{theo:anexe:1} in holomorphic case $n=2$}

We first prove the following lemma.

\begin{Lemma}
\label{lm:permu:2:ele:2}
Let $\mathcal C \subset \Q_g(k_1,k_2)$ be a connected component with
any $k_1,\ k_2 > 0$ and $g\geq 2$. Then there exists 
$(\Su(\pi,\lambda),\Dq) \in \mathcal C$ a half-translation surface 
such that one of the two following holds.

\begin{itemize}

\item $(\Su(\pi,\lambda),\Dq)$ has a multiplicity one saddle
connection between the two zeroes.

\item The combinatorics of $\pi$ has the form
$$
\pi = \left( \begin{smallmatrix} 0_1 & 1 & \dots & r & 0_1 & r+1 & \dots & r+l \\
0_2 & \sigma_1(1) & \dots & \sigma_1(r) & 0_2 & \sigma_2(r+1) & \dots & \sigma_2(r+l)
\end{smallmatrix} \right) \textrm{ or }
\pi = \left( \begin{smallmatrix} A \\ B \end{smallmatrix} \right)
\textrm{ with } \pi(A)=A,\ \pi(B)=B
$$
where $\sigma_1,\sigma_2$ are ``true'' permutations of the sets
$\{1,\dots,r\}$ and $\{r+1,\dots,r+l\}$. 

\end{itemize}
\end{Lemma}

\begin{proof}[Proof of Lemma~\ref{lm:permu:2:ele:2}]
The proof is similar to the one of Lemma~\ref{lm:permu:2:ele:2}. 
It parallels the proof of Claim~\ref{claim:1}, Claim~\ref{claim:2}
and Claim~\ref{claim:3}. 
\end{proof}

Therefore in order to prove Theorem~\ref{theo:anexe:1} in holomorphic
case $n=2$ we have to analyse the two permutations of
Lemma~\ref{lm:permu:2:ele:2}.

\begin{NoNumberLemma}
Let $\pi = \left( \begin{smallmatrix}
0_1 & 1 & \dots & r & 0_1 & r+1 & \dots & r+l \\
\sigma_1(r) & 0_2 & \sigma_2(r+1) & \dots & \sigma_2(r+l) &
0_2 & \sigma_1(1) & \dots & \sigma_1(r-1)
\end{smallmatrix} \right)$. Let us assume that $\pi \not =
\Pi_1(r,l)$. Then there exists $\pi' \sim \pi$ such that $\pi'$ is
irreducible and $\gamma(\pi')$ is a saddle connection between the
two zeroes.
\end{NoNumberLemma}

\begin{proof}[Proof of the lemma]
First let us remark that the separatrix $\gamma(\pi)$ on $\Su(\pi,\lambda)$
is actually a saddle connection. Note also that $\pi$ satisfies
the Condition~$(*)$; Thus the permutation $\pi$ is irreducible if and 
only if $\pi$ is weakly irreducible. \medskip

\noindent One checks directly that $\pi$ is reducible if and only if
$\sigma_1(r)=1$. In this case one can put $\pi$ into the form
\begin{multline*}
\pi=\left( \begin{smallmatrix}
0_1 & 1 & \dots & r & 0_1 & r+1 & \dots & r+l \\
\sigma_1(r) & 0_2 & \sigma_2(r+1) & \dots & \sigma_2(r+l) &
0_2 & \sigma_1(1) & \dots & \sigma_1(r-1)
\end{smallmatrix} \right) \sim \pi'=\\
=\left( \begin{smallmatrix}
0_1 & 1 & 2 & \dots & r & 0_1 & r+1 & \dots & r+l \\
\sigma_1(r-1) & 1 & 0_2 & \sigma_2(r+1) & \dots & \sigma_2(r+l) &
0_2 & \sigma_1(1) & \dots & \sigma_1(r-2)
\end{smallmatrix} \right).
\end{multline*}
Repeating this process to $\pi'$ this yields to the following
dichotomy. Either there exists an irreducible permutation $\pi'$ with
$\pi'\sim \pi$, and so the lemma is proved, or 
$$
\left\{
\begin{array}{ll}
\sigma_1(i) = r-i+1    &  \textrm{for } i=1,\dots,r \\
\sigma_2(j) = 2r+l-j+1 &  \textrm{for } j=r+1,\dots,r+l
\end{array} \right.
$$
These last equalities mean that $\pi = \Pi_1(r,l)$ which is a
contradiction. The lemma is proven.
\end{proof}

\begin{NoNumberLemma}
Let $\pi=\left( \begin{smallmatrix} A \\ B \end{smallmatrix} \right)$
with $\pi(A)=A,\ \pi(B)=B$. Let us also assume that $\pi \not =
\Pi_2(r,l)$. Then there exists $\pi' \sim \pi$ such that $\pi'$ is
irreducible and $\gamma(\pi')$ is a saddle connection between the
two zeroes.
\end{NoNumberLemma}

\begin{proof}
The proof is completely similar to the previous proof and left to the
reader.
\end{proof}

Now we are ready to prove Theorem~\ref{theo:anexe:1} in holomorphic
case $n=2$. Let $\mathcal C$ be any component of any stratum
$\Q_g(k_1,k_2)$. Applying Lemma~\ref{lm:permu:2:ele:2} combining with
the two above lemma one gets the following dichotomy. There exits $\Su
\in \mathcal C$ such that either $\Su$ possesses a multiplicity one
saddle connection between two zeroes, and so the theorem is
proved, or $\Su=\Su(\pi,\lambda)$ with $\pi = \Pi_1(r,l)$ or $\pi =
\Pi_2(r,l)$. These last two permutations correspond to hyperelliptic
permutations and in particular it implies that $\mathcal C$ is itself
a hyperelliptic component of $\Q(4(g-k)-6,4k+2)$ (where $r=2k+1$ and
$l=2(g-k)-3$) which is a contradiction. The theorem is proven.

\subsection{Proof of Theorem~\ref{theo:anexe:1} in case $\Q(-1,4g-3)$}

The proof is decomposed in two steps. We recall that a simple 
cylinder is a metric cylinder with boundary components consisting 
of two {\it single} homologous separatrices. A simple cylinder has
multiplicity one if its boundary has multiplicity one.

\begin{Proposition}
\label{prop:contrac:ik:1}
Let $\mathcal C \subset \Q(-1,4g-3)$ with $g \geq 3$. Then there
exists $\Su(\pi,\lambda) \in \mathcal C$ a half-translation surface
such that one of the two followings holds.

\begin{itemize}

\item $\Su$ has a multiplicity one saddle connection between the zero
  and the pole.

\item $\Su$ has a multiplicity one simple cylinder (i.e. 
$\mathcal C = \mathcal C' \oplus s$ with $s \in \{1,\dots,2g \}$
and $\mathcal C' \subset \Q(-1,4(g-1)-3)$ is component).

\end{itemize}
\end{Proposition}

\begin{Proposition}
\label{prop:contrac:ik:2}
The following assertions hold.

\begin{enumerate}

\item 
\label{prop:point:1}
For any $s \in \{1,2,4\}$, there exists a flat surface
$(\Su,\Dq) \in \mathcal  \Q(-1,5)\oplus s$ with a  multiplicity one saddle
connection.

\item 
\label{prop:point:2}
Any connected component of $\Q(-1,13)$ possesses a half-translation 
surface equipped with a multiplicity one saddle connection. 

\item 
\label{prop:point:3}
If $(\Su,\Dq)$ has a multiplicity one saddle connection then
$(\Su,\Dq) \oplus s$ has also a multiplicity one saddle connection by
bubbling a handle.

\end{enumerate}
\end{Proposition}

We address the proof of the second proposition in the appendix
(see Proposition~\ref{prop:contrac:ik:2:appendix}).
Let us first show how these two propositions give the theorem and then
let us prove the first proposition.

\begin{proof}[Proof of Theorem~\ref{theo:anexe:1}]
It is done by induction on $g$. The theorem is already proven in  
case $g=4$
(Proposition~\ref{prop:contrac:ik:2},~(\ref{prop:point:2})). Let
$g > 4$ be any integer and let $\mathcal C \subset \Q(-1,4g-3)$ be a
connected component. One has to prove that $\mathcal C$ possesses a
surface with a multiplicity one saddle connection. \medskip

\noindent Thanks to Proposition~\ref{prop:contrac:ik:1}, either 
$\mathcal C$ possesses a surface with a multiplicity one saddle
connection, and then the theorem is proved, or $\mathcal C$ contains a
surface with a multiplicity one simple cylinder. In the last
eventuality we have $\mathcal C = \mathcal C'
\oplus s$ where $\mathcal C' \subset \Q(-1,4(g-1)-3)=\Q(-1,4g'-3)$ is
a connected component. We have $4 \leq g' < g$ thus by assumption
$\mathcal C'$ contains a surface $\Su'$ with a multiplicity one saddle
connection. Then we conclude by
Proposition~\ref{prop:contrac:ik:2},~(\ref{prop:point:3}): $\Su' 
\oplus s \in \mathcal C$ is the required surface and the theorem is
proven. \medskip  

\noindent It remains to prove the theorem for genus $g=3$ case. Again
form Proposition~\ref{prop:contrac:ik:1}, if
$\mathcal C \subset \Q(-1,9)$ has no surface with multiplicity one
saddle connection then $\mathcal C$ has the form $\Q(-1,5) \oplus s$
with $s=1,2,3,4$. Moreover
Lemma~\ref{lm:cc:irred:oplus},~(\ref{lm:pt:4}) implies 
$\Q(-1,9)^{irr}=\Q(-1,5)\oplus 3$. Therefore the theorem follows from 
Proposition~\ref{prop:contrac:ik:2},~(\ref{prop:point:1}) which
examines cases $s=1,2,4$. The theorem is proven.
\end{proof}

\begin{proof}[Proof of Proposition~\ref{prop:contrac:ik:1}]
Let $\mathcal C \subset \Q(-1,4g-3)$ be a connected component with $g
\geq 3$ and $\Su(\pi,\lambda) \in \mathcal C$. The surface $\Su$ has a
pole and a unique zero $P$. Hence one can put $\Su$ into the form
prescribed by Figure~\ref{fig:pole:zero}.
\begin{figure}[htbp]

\psfrag{C}{$\scriptstyle \CYL$}
\psfrag{P}{$\scriptstyle P$}
\psfrag{PO}{$\scriptstyle \textrm{pole}$}  
\psfrag{dots}{$\scriptstyle \dots\dots\dots\dots$}
\psfrag{I}{$\scriptstyle I$}
\psfrag{J}{$\scriptstyle J$}
\psfrag{g}{$\scriptstyle \gamma(\pi)$} \psfrag{a}{$\scriptstyle \alpha$}

\begin{center}
\includegraphics[width=7cm]{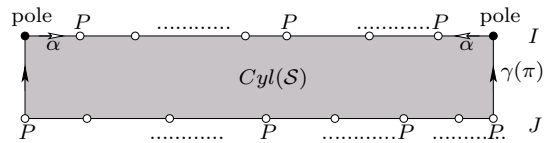}
\end{center}
\caption{
\label{fig:pole:zero}
A half-translation surface with a unique zero and a unique pole. 
}
\end{figure}
\noindent Let $\alpha$ be the horizontal saddle connection between the pole and
the zero.

\begin{Claim}
Either $\alpha$ has multiplicity one or $\Su(\pi,\lambda) = 
\Su(\pi',\lambda')$ with $\pi'=\bigl( \begin{smallmatrix} 0 & 1 &
  \dots & r & 0 \\ A \end{smallmatrix}\bigr)$ and 
$\pi'(\{1,\dots,r\}) \subset A$.
\end{Claim}

\begin{proof}[Proof of the claim]
Apply Proposition~\ref{prop:fundamental:2}.
\end{proof}

\noindent Therefore thanks to above claim let us assume that $\Su(\pi,\lambda)
\in \mathcal C$ with $\pi$ a permutation given by previous claim. We 
define $\sigma\sim \pi$ by 
$\sigma=\bigl( \begin{smallmatrix} 1 & 2 & \dots & r & 0 & 0 \\ 
1 & A \end{smallmatrix}\bigr)$ and $\sigma(\{1,\dots,r\}) \subset A$.
Note that the vertical direction on $\Su(\sigma,\lambda)$ 
decomposes the surface into (at least) one simple cylinder 
(see also Figure~\ref{fig:simple:cyl}). The reducibility of 
$\hat \sigma$ involves one of the two followings possible 
decompositions.
$$
\begin{array}{lll}
\scriptstyle{Type }(1) &
\sigma=\bigl( \begin{smallmatrix} 1 & 2 & \dots & k & |&k+1 & \dots & r & 0 & 0 \\
1 & A_1 & k & A_2&| &  A' \end{smallmatrix}\bigr) &
\qquad \textrm{ with } \sigma(\{1,\dots,r\}) = A_1 \sqcup A_2 \\
\scriptstyle{Type }(2) &
\sigma=\bigl( \begin{smallmatrix} 1 & 2 & \dots & r & 0 &|& 0 \\
1 &  &  A &  & |& A' \end{smallmatrix}\bigr) &
\qquad \textrm{ with } \sigma(A') \subset A
\end{array}
$$

\begin{Claim}
If $\hat \sigma$ involves reducibility of Type~$(1)$ 
then there exists $\sigma_1\sim \sigma$ such that $\hat \sigma_1$ 
is irreducible.
\end{Claim}

\begin{proof}[Proof of the claim]
The proof parallels the one of Theorem~\ref{theo:simple:cyl} versus
Abelian differentials (Section~\ref{sec:step:minimal}).
\end{proof}

\noindent For reducible permutations $\sigma$ of Type~$(2)$, let us denote $A =
(A_1 \ k \ A_2)$ with $\sigma(A_2) \subseteq A'$ and $2 \leq k \leq r$
(eventually $A_2 = \emptyset$).

\begin{Claim}
If $\hat \sigma$ involves reducibility of Type~$(2)$ 
then either there exists $\sigma_1\sim \sigma$ such that 
$\hat \sigma_1$ is irreducible or $k=r$.
\end{Claim}

\begin{proof}[Proof of the claim]
With above notations, if $k \not = r$ then one has 
$$
\sigma = \bigl( \begin{smallmatrix} 1 & 2 & \dots & k & \dots & r & 0 & 0 \\
1 &  A_1 & k & A_2 &  & A' \end{smallmatrix}\bigr) \sim \sigma_1 = 
\bigl( \begin{smallmatrix} k & k+1 & \dots & r & 0 & 0 & 1 & 2 & \dots & k-1  \\
k & A_2 & A' & 1 & A_1 \end{smallmatrix}\bigr). 
$$
If $\hat \sigma_1$ is irreducible then the claim is proven. Otherwise
it is easy to see that $\hat \sigma_1$ involves reducibility of
Type~$(1)$. Thus previous claim applies.
\end{proof}
\noindent Finally the last case to consider is the following one: 
$\hat \sigma$ is reducible with decomposition of Type~$(2)$ and
$k=r$. Reporting these data into the permutation $\pi$ one gets
$\pi= \bigl( \begin{smallmatrix} 0 & 1 & \dots & r & 0  \\
r & A_2 & A' & 1 & A_1 \end{smallmatrix}\bigr)$ with $\pi(A_2)
\subseteq A'$. 

\begin{Claim}
The vertical separatrix $\gamma(\pi)$ on $\Su(\pi,\lambda)$ is 
a saddle connection. Moreover the above permutation $\pi$ is reducible
if and only if
$$
\pi= \bigl( \begin{smallmatrix} 0 && 1 & \dots & \dots & r & 0  \\
r & A' & 1 & \pi(2) & \dots & \pi(r-1) & A'' \end{smallmatrix}\bigr)
$$
with $A''=\pi(A')$.
\end{Claim}

\begin{proof}[Proof of the claim]
Let $\pi= \bigl( \begin{smallmatrix} 0 & 1 & \dots & r & 0  \\
r & A_2 & A' & 1 & A_1 \end{smallmatrix}\bigr)$ with $\pi(A_2)
\subseteq A'$. 
A direct observation shows that if $\pi$ is reducible then $A_2 =
\emptyset$. The same approach shows that if $\pi$ is reducible then 
$A_1 =\pi(2) \ \dots \ \pi(r-1) \ A''$ with $\pi(A')=A''$. The claim 
is proven.
\end{proof}

\noindent The conclusion of previous claims is the following one. We
have proved that $\mathcal C$ contains a surface $\Su(\pi,\lambda)$
such that either $\Su$ satisfied the conclusions of the proposition or 
$$
\pi= \bigl( \begin{array}{ccccccc} 0 &1 & \dots & & \dots & r & 0  \\
r & A' & 1 & \pi(2) & \dots & \pi(r-1) & A'' \end{array}\bigr) 
\qquad \textrm{ and } \pi(A') = A''.
$$
Then we have to study the two last cases $r\geq 2$ and $r=1$. This
corresponds to the two next lemma which end the proof of
Proposition~\ref{prop:contrac:ik:1}. 
\end{proof}

\begin{NoNumberLemma}
Let $r \geq 2$ be any integer. Let 
$\pi= \bigl( \begin{smallmatrix} 0 && 1 & \dots & \dots & r & 0  \\
r & A' & 1 & \pi(2) & \dots & \pi(r-1) & A'' \end{smallmatrix}\bigr)$ 
be a permutation with $\pi(A') = A''$. We also assume that 
$\Su(\pi,\lambda) \in \Q(-1,4g-3)$ with $g \geq 3$. \\
Then there exists a permutation $\pi'\sim \pi$ such that 
$\Su(\pi',\lambda)$ has a multiplicity one saddle connection.
\end{NoNumberLemma}

\begin{proof}[Proof of the lemma]
Let us assume that $r \geq 3$, then $r-1 \not = 1$. Let us introduce
some notations to clarify the situation. Let $C$ stand for
$(\pi(2)\ \dots \ \pi(r-1))$. One decomposes $C$ into the following way:
$C=(C_1 \ r-1 \ C_2)$. Equipped with these notations:
$\pi= \bigl( \begin{smallmatrix} 0 & 1 & \dots & & r-1 & r & 0  \\
r & A' & 1 & C_1 & r-1 & C_2 & A'' \end{smallmatrix}\bigr) \sim \sigma =
\bigl( \begin{smallmatrix} 0 & 1 & \dots & \dots & r-1 & r & 0 \\
r-1 & C_2 & A'' & r & A' & 1 & C_1 \end{smallmatrix}\bigr)$. A direct 
observation shows that $\sigma$ is irreducible and $\gamma(sigma)$ is 
a saddle connection. The lemma then follows from 
Theorem~\ref{theo:fundamental:1}. \medskip

\noindent Now let us assume that $r = 2$. The list $A''$ is non-empty.
Let us denote $A''=(3 \ B)$ with $\pi(3)\in A'$. We also 
use $A'=(B_1 \ 3 \ B_2)$. Thus 
$\pi= \bigl( \begin{smallmatrix} 0 & 1 & 2 & 0  \\
2 & A' & 1 & A'' \end{smallmatrix}\bigr) \sim \sigma =
\bigl( \begin{smallmatrix} 0 & 1 & 2 && 0  \\
B & 2 & B_1 & 3 & B_2 & 1 & 3 \end{smallmatrix}\bigr)$. 
This permutation is reducible if and only if $B=B_1=B_2=\emptyset$. 
This case consists of $\pi=\bigl( \begin{smallmatrix} 0 & 1 & 2& 0  \\
2 & 3 & 1 & 3 \end{smallmatrix}\bigr)$ and then $\Su(\pi,\lambda) 
\in \Q(-1,5)$ contradicting $g\geq 3$. Hence $\sigma$ is irreducible 
and the lemma again follows from Theorem~\ref{theo:fundamental:1}.
\end{proof}

\noindent Next lemma studies $r=1$ case.

\begin{NoNumberLemma}
Let $\pi= \bigl( \begin{smallmatrix} 0 & 1 &
  0   \\ A' & 1 & A'' \end{smallmatrix}\bigr)$ be a permutation 
with $\pi(A')=A''$. Then one of the following assertions holds.
\begin{itemize}

\item There exists $\pi' \sim \pi$ such that $\Su(\pi',\lambda)$ 
has a multiplicity one saddle connection.

\item The combinatorics of $\pi$ is given by $A'=A''=(2 \ 3 \ \dots \ l)$.
Moreover there exists $\pi' \sim \pi$ and an admissible
vector $\lambda_0$ such that $\Su(\pi',\lambda_0)$ has a multiplicity
one simple cylinder.

\end{itemize}
\end{NoNumberLemma}

\begin{proof}[Proof of the lemma]
Let $\pi$ the permutation $\pi= \bigl( \begin{smallmatrix} 0 & 1 & 0   \\
A' & 1 & A'' \end{smallmatrix}\bigr)$. The list $A'$ is non-empty
thus let us introduce $A' = (2 \ A'_2)$ with $\pi(2)\in A''$. Hence
we can also introduce lists $A''_1,\ A''_2$ such that
$A'' = (A''_1 \ 2 \ A''_2)$. According to these notations, one has
$\pi= \bigl( \begin{smallmatrix} 0 & & 1 & & & 0   \\
2 & A'_2 & 1 & A''_1& 2 & A''_2 \end{smallmatrix}\bigr) \sim
\pi'=\bigl( \begin{smallmatrix} 0 & & 1 & & & 0   \\
2 & A''_2 & 2 & A'_2 & 1 & A''_1 \end{smallmatrix}\bigr)$.
A straightforward computation shows that the permutation $\pi'$ is reducible
if and only if $A''_1  = \emptyset$. Thus, repeating inductively
this process with $A' := A'_2$ and $A'' := A''_2$, we get 
either $\Su(\pi',\lambda)$ has a multiplicity one saddle connection or
$$
\pi \sim \pi' = \left(
\begin{array}{ccccccccccc}
0 & & & & & 1 & & & & & 0 \\
3 & 4 & \dots & l & 2 & 3 & 4 & \dots & l & 1 & 2
\end{array}
\right)
$$
In order to finish the proof of the lemma, we have to present 
an admissible vector $\lambda_0$
such that the vertical foliation on $\Su(\pi',\lambda_0)$ is
completely periodic and decomposes the surface with at least a
multiplicity one simple cylinder. Let us consider the admissible
vector (for $\pi'$):
\begin{equation}
\label{eq:adm:vect}
\lambda_0 = \left( \ (l-1)\alpha \ ;\ \alpha \ ;\ (l-1)\alpha \ ;\
\underbrace{\alpha \ ;\ \dots \ ;\ \alpha}_{2l-1 \textrm{ times}} \ \right)
\qquad \textrm{ for any } \alpha \in \mathbb R^{+}
\end{equation}
It is easy to see that the vertical foliation on $\Su(\pi',\lambda_0)$
decomposes the surface into $g-1$ cylinders 
(see Figure~\ref{fig:vertical:simple:i9}). 
One checks that the vertical cylinder corresponding to the horizontal 
interval numbered $l$ is a multiplicity one simple cylinder (see the figure).
\begin{figure}[htbp]
\begin{center}
\psfrag{0}{$0_h$}  \psfrag{3}{$3_h$}
\psfrag{1}{$1_h$}  \psfrag{4}{$4_h$}
\psfrag{2}{$2_h$}  \psfrag{5}{$5_h$}

\psfrag{a}{$0_v$}  \psfrag{d}{$3_v$}
\psfrag{b}{$1_v$}  \psfrag{e}{$4_v$}
\psfrag{c}{$2_v$}  \psfrag{f}{$5_v$}
\psfrag{g}{$6_v$}  \psfrag{h}{$7_v$}
\psfrag{i}{$8_v$}

\psfrag{00}{$0$}  \psfrag{30}{$3$}
\psfrag{10}{$1$}  \psfrag{40}{$4$}
\psfrag{20}{$2$}  \psfrag{50}{$5$}

\psfrag{u1}{deforming continuously}
\psfrag{u2}{vertical parameters}

\psfrag{c1}{$C_1$}  \psfrag{c2}{$C_2$}
\psfrag{e}{$\eta$}  \psfrag{S}{$\Su(\pi',\lambda)$}

\includegraphics[width=10cm]{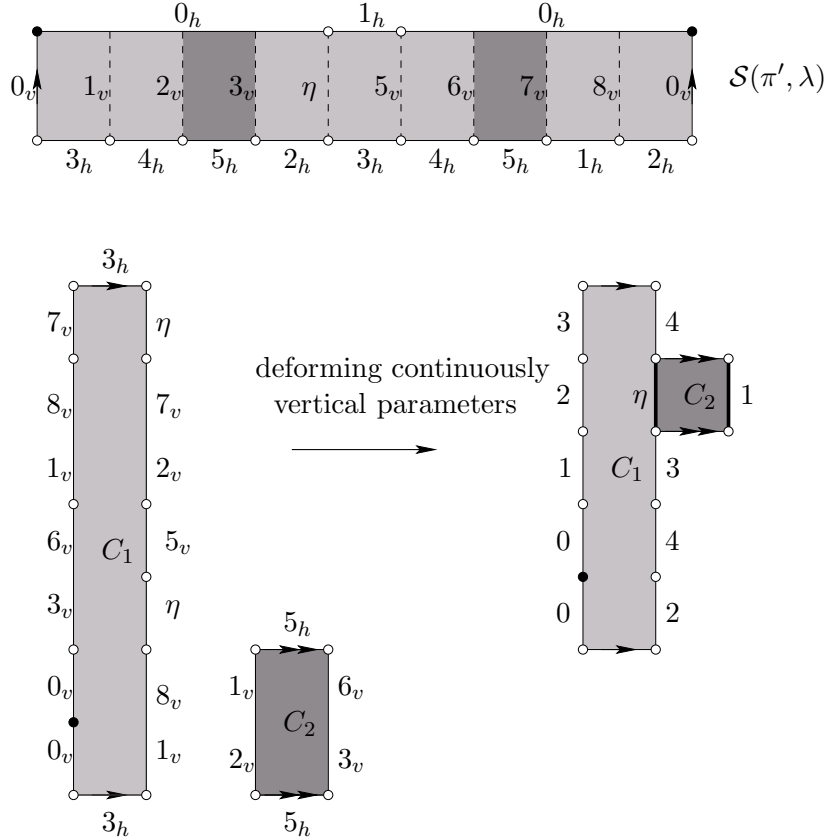}
\end{center}
\caption{
\label{fig:vertical:simple:i9}
A half-translation surface $\Su'=\Su(\pi'\lambda_0)$ suspended 
over the permutation $\pi'$ and the admissible vector $\lambda_0$ 
(see Equation~(\ref{eq:adm:vect}). The vertical
foliation decomposes $\Su'$ in two cylinders $C_1$ and $C_2$. 
The cylinder $C_2$ is a simple cylinder: the conical angle of its 
boundary $\eta$ is $\pi$. Moreover the length of $\eta$ can be
chosen arbitrary small with respect to other vertical parameters: 
$C_2$ has multiplicity one.
}
\end{figure}
In Figure~\ref{fig:vertical:simple:i9}, we present a complete description for
the surface given by the case $r=5$ (i.e. $g=3$). This completes the proof 
of the lemma.
\end{proof}


Now we can deduce from Theorem~\ref{theo:anexe:1} our main result,
namely Theorem~\ref{theo:adja}.

\subsection{Proof of Theorem~\ref{theo:adja}}

We first establish 

\begin{Proposition}[Reformulation of Theorem~\ref{theo:anexe:1} in
holomorphic case $n=2$]
\begin{multline*}
\mathcal C \subseteq \Q(k_1,k_2), \ k_i >0,
\textrm{ is hyperelliptic} \iff \\ \forall  \ (\Su,\Dq) \in \mathcal
C, \textrm{ any saddle connection on $\Su$ has multiplicity at least } 2
\end{multline*}
\end{Proposition}

\noindent Through the proof of Theorem~\ref{theo:anexe:1} in case
$n=3$ and $n=4$ one gets a similar characterization of hyperelliptic
components (compare with Lemma~\ref{lm:characteriz:hyper}).

\begin{Proposition}
\label{prop:cor:theo}
Hyperelliptic components of strata with $3$ and $4$
singularities are characterized by the following:
\begin{enumerate}

\item $\mathcal C \subseteq \Q(k_1,k_2,k_3)$ is hyperelliptic
$\iff \exists i_0 \not = i_1 \in \{1,2,3 \}, \forall \ (\Su,\Dq) \in
  \mathcal C$, \\
any saddle connection on $\Su$ between $P_{k_i}$ and $P_{k_j}$, with
$i\in \{i_0,i_1\}$ and $j\not \in \{i_0,i_1\}$
has multiplicity at least $2$.

\item $\mathcal C \subseteq \Q(k_1,k_2,k_3,k_4)$ is
  hyperelliptic
$\iff \exists i_0 \not = i_1 \in \{1,2,3,4\}$ \\
$ \forall  \ (\Su,\Dq) \in  \mathcal C$,
any saddle connection on $\Su$ between $P_{k_i}$ and $P_{k_j}$,
with $i\in \{i_0,i_1\}$ and $j\not \in \{i_0,i_1\}$, has multiplicity
at least $2$.

\end{enumerate}
\end{Proposition}

\begin{Proposition}
\label{prop:collapse:non-hyper}
Let $\mathcal C$ be a connected component of a $\Q(k_1,\dots,k_n)$
with $n\geq 5$. Then for each $i\not = j$, there exists a half-translation
surface $\Su \in \mathcal C$ with a multiplicity one saddle connection
between two singularities $P_{k_i},P_{k_j} \in \Su$ of multiplicities
$k_i$ and $k_j$.
\end{Proposition}

\begin{proof}[Proof of Theorem~\ref{theo:adja}]
We discuss the theorem following the different values of $n$ (the number
of singularities). First let us note that if we prove the theorem in
case $n=2,3,4,5$ then the theorem follows for any $n\geq 6$ because  
the locus of hyperelliptic and irreducible components is located on
strata with $2,3$ and $4$ singularities. \medskip

\noindent The theorem in case $n=2$ corresponds to
Theorem~\ref{theo:anexe:1} which is already proved. \medskip

\noindent Thus let us assume that $n=3$. Let $\mathcal C \subseteq
\Q(k_1,k_2,k_3)$ be a
non-hyperelliptic component. Recall that all $k_i$ are non-zero.
Up to re-organize $k_i$, one can assume $k_1 \leq k_2 \leq
k_3$. Proposition~\ref{prop:cor:theo} 
with $i_0=1,\ i_1=2$ gives a half-translation surface
$(\Su,\Dq)\in \mathcal C$ with a multiplicity one saddle connection
between $P_{k_3} \in \Su$ and $P_{k_i} \in \Su$ for $i\in \{1,2\}$. In other
terms there exists a component $\mathcal C'$ of the strata
$\Q(k_1+k_3,k_2) \sqcup \Q(k_2+k_3,k_1)$ such that $\mathcal C'
\subset \overline{\mathcal C}$. There is three possibilities:
$\mathcal C'$ is hyperelliptic, irreducible or ``regular''
(neither hyperelliptic or irreducible).

\noindent If $\mathcal C'$ is regular then Theorem~\ref{theo:anexe:1}
leads to the result. Indeed on can find $\mathcal C_0 \subseteq \Q(4g-4)$ a
component with $\mathcal C_0\subset \overline{\mathcal C'}$.
But $\overline{\mathcal C'} \subset \overline{\mathcal C}$
therefore $\mathcal C_0 \subset \overline{\mathcal C}$ and we are
done.

\noindent Assume that $\mathcal C'$ is hyperelliptic.
Applying Proposition~\ref{prop:adja:two:singu}, one can connect the
hyperelliptic component $\mathcal C'$ to another non-hyperelliptic
component $\mathcal C''$ of $\Q(k_i+k_3,k_j)$ passing through the
stratum $\Q(k_1,k_2,k_3)$. Therefore there exists a component
$\mathcal C_2 \subseteq \Q(k_1,k_2,k_3)$ such that $\mathcal C' \sqcup
\mathcal C'' \subset \overline{\mathcal C_2}$. Thus
$\mathcal C' \subset \overline{\mathcal C} \sqcup \overline{\mathcal
  C_2}$. The assumption $k_1 \leq k_2 \leq k_3$ implies 
$k_i+k_3 \not = k_j$ thus Corollary~\ref{cor:continuous} implies
$\mathcal C_2 = \mathcal C$. Therefore $\mathcal C$ is adjacent to a
non-hyperelliptic component of $\Q(k_i+k_3,k_j)$ and we are again done.

\noindent Finally let us assume that $\mathcal C'$ is irreducible;
that is $\mathcal C'=\Q^{irr}(-1,9)$. Recall that the component
$\mathcal C \not =\Q^{irr}(-1,3,6)$ therefore the theorem follows 
from Proposition~\ref{prop:adja:irre} and
Corollary~\ref{cor:continuous}. \medskip

\noindent The proofs for cases $n=4,5$ are similar to the above
discussion. Theorem~\ref{theo:adja} is proven.
\end{proof}

\section{Proof of Main Theorem~\ref{theo:main:1} and Main
  Theorem~\ref{theo:main:2}} 
\label{sec:proof}

\begin{proof}[Proof of Main Theorem~\ref{theo:main:1}]
As mentioned in Section~\ref{sec:link:main:result},
Main Theorem~\ref{theo:main:1}
follows from Theorem~\ref{theo:adja} on adjacency of strata and
Theorem~\ref{theo:classi:minimal} on classification of minimal
strata.
\end{proof}

\begin{proof}[Proof of Main Theorem~\ref{theo:main:2}]
Theorem~\ref{theo:main:2} is already proved for genus $0$ case (see 
Proposition~\ref{prop:particular:cases}). The genus $3$ case is also proved
(see Theorem~\ref{theo:adja} and Section~\ref{sec:link:main:result}). \medskip

\noindent Thus let us consider genera $1,2$ cases. One can first
prove, by a direct argument on finitely many strata, 
Theorem~\ref{theo:main:2} for $g=1,2$ and $n\leq 5$ cases. Indeed the
number of such a strata of $\Q_2$ and of $\Q_1$ with $n=2,3,4,5$ is
finite. We can check, using Rauzy classes, that the theorem holds. In 
particular we get that any stratum of $\Q_2$ with $n=5$ is
connected and any stratum of $\Q_1$ with $n=5$ is also connected. 

Now let us prove that any stratum of $\Q_1$ with $n\geq 6$
is also connected. Let $\mathcal C_1,\mathcal C_2$ be two components
of $\Q(k_1,\dots,k_6)$. Let us assume that $k_1\leq \dots \leq
k_6$. Then $k_5+k_6\not = k_i$ for any $i=1,\dots,4$. By
Proposition~\ref{prop:collapse:non-hyper} there exists two components
$\mathcal C_0$ and $\mathcal C'_0$ of $\Q(k_1,\dots,k_4,k_5+k_6)$ with
$\mathcal C_0 \subset \overline{\mathcal C_1}$ and $\mathcal C'_0 \subset
\overline{\mathcal C_2}$. Any stratum of $\Q_1$ with $n = 5$ being connected,
one has $\mathcal C_0 = \mathcal C'_0$ hence $\mathcal C_0 \subset
\overline{\mathcal C_1} \cup \overline{\mathcal C_2}$. Therefore
Corollary~\ref{cor:continuous} applies and $\mathcal C_1
= \mathcal C_2$; so that any stratum of $\Q_1$ with $n=6$ is
connected. Above argument extends to any stratum of $\Q_1$ with $n\geq
6$. Theorem~\ref{theo:main:2} in case $g=1$ is therefore proven.  The same
argument holds for the genus $2$ case. \medskip

\noindent We now prove Theorem~\ref{theo:main:2} in case $g=4$. We
have already
proved that any component of a stratum of $\Q_4$ which is not
irreducible neither hyperelliptic is adjacent to $\Q^{irr,I}(12)$ or
$\Q^{irr,II}(12)$ (Theorem~\ref{theo:adja}). In order to prove the
theorem it suffices to show that those components are in fact
adjacent to $\Q^{irr,I}(12)$. By Theorem~\ref{theo:anexe:1}, any
non-irreducible or 
non-hyperelliptic component is adjacent to a non-irreducible or
non-hyperelliptic component of a stratum of $\Q_4$ with $n=2$. Let
$\mathcal C$ be any non-irreducible or non-hyperelliptic component of
$\Q(k_1,k_2)$ with $k_1+k_2 = 12$. One has to prove that $\mathcal C$
is adjacent to $\Q^{irr,I}(12)$. By Theorem~\ref{theo:adja}, $\mathcal
C$ is adjacent to $\Q^{irr,I}(12)$ or $\Q^{irr,II}(12)$. Let us assume
that second case holds. By Proposition~\ref{prop:adja:12} there
exists $\mathcal C' \subset \Q(k_1,k_2)$ such that $\Q^{irr,II}(12)
\subset \overline{\mathcal C'}$ and $\Q^{irr,I}(12) \subset
\overline{\mathcal C'}$. Therefore $\Q^{irr,II}(12) \subset
\overline{\mathcal C'} \cup \overline{\mathcal C}$ and
Corollary~\ref{cor:continuous} implies $\mathcal C =
\mathcal C'$. Theorem~\ref{theo:main:2} is proven.
\end{proof}

\appendix

\section{Connectedness of particular strata}

This section is devoted to the computation of the set of connected 
components of some strata in low dimensional case.

\begin{Lemma}
\label{lm:appendix:q8}
The stratum $\Q(8)$ is connected.
\end{Lemma}

\begin{proof}
Let $(\Su(\pi,\lambda),\Dq)$ be a genus three half-translation surface,
with a single singularity. The conical angle around this singularity
is $10\pi$. We will show that the possibilities for the
combinatorics of the gluing maps of the set of horizontal separatrix 
loops is very restrictive.

Recall that to each permutation $\pi$ one associates two lines
of the table representing $\pi$
(see~\ref{sec:generalized:permutations}). There are obviously
two possibilities for these two lines: either the numbers of elements
are equal or they are different. One defines $\mathcal A_1$ 
the set of generalized permutations corresponding to the first case. 
We also define $\mathcal A_2$ the set of generalized
permutations corresponding to the second case. \medskip

\noindent By assumption $(\Su,\Dq) \in \Q(8)$, hence there are exactly five
separatrices loops. Therefore the number of elements of lines for
permutations in $\mathcal A_1$ is $(5,5)$ and for permutations in
$\mathcal A_2$ is $(6,4)$. A straightforward computation shows 
the two sets $\mathcal A_1,\ \mathcal A_2$ are very simple. More
precisely, up to cyclic order, one has $\# \mathcal A_1=4$ and $\#
\mathcal A_2=3$.
$$
\begin{array}{ll}
\textrm{up to cyclic order:} & \mathcal A_1 = \left\{
\bigl( \begin{smallmatrix}5 & 3 & 5 & 2 & 4  \\1 & 2 & 1 & 3 & 4
\end{smallmatrix}\bigr), \
\bigl( \begin{smallmatrix}5 & 4 & 5 & 2 & 3  \\ 1 & 2 & 1 & 3 & 4
\end{smallmatrix}\bigr), \
\bigl( \begin{smallmatrix}5 & 4 & 5 & 3 & 2  \\ 1 & 2 & 1 & 3 & 4
\end{smallmatrix}\bigr), \
\bigl( \begin{smallmatrix}5 & 3 & 5 & 3 & 4  \\1 & 2 & 1 & 2 & 4
\end{smallmatrix}\bigr) \right\} \\
\textrm{up to cyclic order:} & \mathcal A_2 = \left\{
\bigl( \begin{smallmatrix}5 & 2 & 5 & 3 & 4 & 2  \\ 1 & 3 & 1 & 4
\end{smallmatrix}\bigr),\
\bigl( \begin{smallmatrix}3 & 5 & 4 & 2 & 5 & 2  \\1 & 3 & 1 & 4
\end{smallmatrix}\bigr),\
\bigl( \begin{smallmatrix}5 & 3 & 2 & 5 & 4 & 2  \\1 & 3 & 1 & 4
\end{smallmatrix}\bigr) \right\}
\end{array}
$$
This proves that the stratum $\Q(8)$ has at most $7$ connected
components. Now, let us consider surfaces $\Su(\pi,\lambda_2)$ with
$\pi \in \mathcal A_2$ and $\lambda_2 = (1,1,1,1,1,1,2,1,2,1)$. A
direct verification shows that the vertical foliation on surfaces
$\Su=\Su(\pi,\lambda_2)$ decomposes $\Su$ into a single
cylinder. Therefore, we get a permutation encoding this cylinder. One
can check that this permutation belongs to the set $\mathcal
A_1$. Thus procedure connects permutations of $\mathcal A_2$ to 
permutations of $\mathcal A_1$. \medskip

\noindent To conclude, let us consider surfaces $\Su(\pi,\lambda_1)$ with $\pi
\in \mathcal A_1$ and $\lambda_1 = (1,1,1,1,1,1,1,1,1,1)$. These
surfaces are arithmetic surfaces. It is easy to check that all of these
surfaces belong to the same $\textrm{PSL}_2(\mathbb Z)$-orbit. The
lemma is proven.
\end{proof}

\begin{Lemma}
\label{lm:appendix:q12}
The stratum $\Q(12)$ has at most two connected components.
\end{Lemma}

\begin{proof}
One has to show that $\Q(12) = \Q^{irr,I}(12) \cup
\Q^{irr,II}(12)$. Thanks to previous lemma, let $\mathcal C_0$ be the
unique connected component of the stratum $\Q(8)$. 
Theorem~\ref{theo:simple:cyl} implies that any component $\mathcal C$
of the stratum $\Q(12)$ has the following form.
\begin{equation}
\label{eq:q12}
\mathcal C = \mathcal C_0 \oplus s \textrm{ with } s=1,\dots,6=2g
\end{equation}
We now recall the construction of the two irreducible components 
$\Q^{irr,I}(12),\ \Q^{irr,II}(12)$. Let 
$\pi_1 = \bigl(
\begin{smallmatrix} 1 & 2 & 3 & 4 & 2 & 5 & 6  \\
1 & 4 & 5 & 7 & 6 & 7 & 3 \end{smallmatrix}\bigr)$ and $\pi_2 = \bigl(
\begin{smallmatrix} 1 & 2 & 3 & 4 & 3 & 5 & 6   \\
1 & 5 & 7 & 4 & 2 & 6 & 7 \end{smallmatrix}\bigr)$ be two
permutations. Let $\Su_i=\Su(\pi_i,\lambda_0)$ be the suspended flat
surfaces with admissible vector $\lambda_0 = (1,1,1,1,1,1,1)$. By 
definition $\Q^{irr,I}(12)$ is the component containing
$\Su_1$ and $\Q^{irr,II}(12)$ is the component containing $\Su_2$. 
Let us recall (see Lemma~\ref{lm:cc:irred:oplus}) that the vertical
foliation on $\Su_i$ produces a multiplicity one simple cylinder which
gives:
$$
\begin{array}{l}
\Q^{irr,I}(12) = \mathcal C_0 \oplus 2 \\
\Q^{irr,II}(12) = \mathcal C_0 \oplus 6 
\end{array}
$$
\noindent We will show:
$$
\begin{array}{l}
\mathcal C_0 \oplus 4 = \mathcal C_0 \oplus 1 = \mathcal C_0 \oplus 5
= \mathcal C_0 \oplus 2 = \Q^{irr,I}(12) \\
\mathcal C_0 \oplus 3 = \mathcal C_0 \oplus 6 = \Q^{irr,II}(12)
\end{array}
$$
Combining with Equation~(\ref{eq:q12}) this will give the
lemma. \medskip

\noindent Let $\mathcal C_0 \oplus s$, $s=1,\dots,6$ be any component
of $\Q(12)$.

\medskip
\noindent {\bf case s=3.}

Let us consider the new permutation $\pi_2 \sim \pi'_2 = \bigl(
\begin{smallmatrix} 5 & 6 & 1 & 2 & 3 & 4 & 3   \\
5 & 7 & 4 & 2 & 6 & 7 & 1 \end{smallmatrix}\bigr)$. Obviously, $\hat
\pi'_2$ is irreducible thus one gets a (vertical) multiplicity one
simple cylinder on $\Su(\pi'_2,\lambda_0)$. This cylinder has angle
$3\pi$. In other words $\mathcal C_0 \oplus 3 =
\Q^{irr,II}(12)$.

\medskip
\noindent {\bf case s=4.}

Let us consider the new permutation $\pi_1 \sim \pi'_1 = \bigl(
\begin{smallmatrix}5 & 6 & 1 & 2 & 3 & 4 & 2   \\ 5 & 7 & 6 & 7 & 3 & 1 & 4
 \end{smallmatrix}\bigr)$. Obviously, $\hat
\pi'_1$ is irreducible thus one gets a (vertical) multiplicity one
simple cylinder on $\Su(\pi'_1,\lambda_0)$. This cylinder has angle
$4\pi$ or in other words $\mathcal C_0 \oplus 4 = \Q^{irr,I}(12)$.

\medskip
\noindent {\bf cases s=1 and s=5.}

Now let us consider the new permutation $\sigma= \bigl(
\begin{smallmatrix} 1 & 2 & 3 & 4 & 5 & 6 & 5   \\
1 & 4 & 7 & 3 & 7 & 2 & 6  \end{smallmatrix}\bigr)$. The permutation
$\hat \sigma$ is irreducible therefore the surface
$\Su(\sigma,\lambda_0)$ has a (vertical) multiplicity one simple
cylinder. This cylinder has angle $4\pi$, therefore this
surface $\Su(\sigma,\lambda_0)$ belongs to component $\Q^{irr,I}(12)$.

Let us consider the new permutation $\sigma \sim \sigma'= \bigl(
\begin{smallmatrix} 3 & 4 & 5 & 6 & 5 & 1 & 2  \\
3 & 7 & 2 & 6 & 1 & 4 & 7 \end{smallmatrix}\bigr)$. Obviously, $\hat
\sigma'$ is irreducible thus one gets a (vertical) multiplicity one
simple cylinder on $\Su(\sigma',\lambda_0)$.
This cylinder has angle $\pi$. In other words 
$\mathcal C_0 \oplus 1 = \Q^{irr,I}(12)$.

Let us consider the new permutation $\sigma \sim \sigma''= \bigl(
\begin{smallmatrix} 2 & 3 & 4 & 5 & 6 & 5 & 1   \\
2 & 6 & 1 & 4 & 7 & 3 & 7  \end{smallmatrix}\bigr)$. Obviously, $\hat
\sigma''$ is irreducible thus one gets a (vertical) multiplicity one
simple cylinder on $\Su(\sigma'',\lambda_0)$.
This cylinder has angle $5\pi$ or in other words
$\mathcal C_0 \oplus 5 = \Q^{irr,I}(12)$.

The lemma is proven.
\end{proof}

\begin{Lemma}
\label{lm:-1:5:appendix}
The stratum $\Q(-1,5)$ is connected.
\end{Lemma}

\begin{proof}[Proof of Lemma~\ref{lm:-1:5:appendix}]
Let $\Su(\pi,\lambda) \in \Q(-1,5)$ be any surface.
One can directly check that (up to cyclic order), the combinatorics of
the permutation $\pi$ has only two possibilities.
$$
\pi_1=\bigl(\begin{smallmatrix} 0 & 0 & 1 & 2   \\ 1 & 3 & 2 & 3 \end{smallmatrix}\bigr)
\qquad \textrm{or} \qquad
\pi_2=\bigl(\begin{smallmatrix} 0 & & 1 & & 0  \\ 2 & 3 & 2 & 1 & 3 \end{smallmatrix}\bigr).
$$
Equipped with the two admissible vectors
$$
\lambda_1=(1,1,2,1,2,1,1,1) \qquad \textrm{and} \qquad 
\lambda_2=(2,1,2,1,1,1,1,1),
$$
we obtain two flat surfaces $\Su_i=\Su(\pi_i,\lambda_i)$, $i=1,2$. A
direct computation shows the vertical foliation on $\Su_2$ produces
the surface $\Su_1$. The lemma is proven.
\end{proof}

\begin{Proposition}
\label{prop:contrac:ik:2:appendix}
The following assertions hold:

\begin{enumerate}

\item For any $s \in \{1,2,4\}$, there exists a flat surface
$(\Su,\Dq) \in \mathcal  \Q(-1,5)\oplus s$ with a  multiplicity one saddle
connection.

\item Any connected component of $\Q(-1,13)$ has a half-translation 
surface with a multiplicity one saddle connection. 

\item If $(\Su,\Dq)$ has a multiplicity one saddle connection then
$(\Su,\Dq) \oplus s$ has also a multiplicity one saddle connection by
bubbling a handle.

\end{enumerate}
\end{Proposition}

\begin{proof}[Proof of Proposition~\ref{prop:contrac:ik:2:appendix}]
We consider separately the three cases.

\noindent {\bf Proof of the first point} 
We will present three surfaces in $\Q(-1,5)\oplus s$ (for $s=1,2,4$)
with a multiplicity one saddle connection. For that we will construct
surfaces in $\Q(-1,9)$ with a multiplicity one simple cylinder of angle
$\pi,2\pi$ and $4\pi$. \medskip 

\noindent Thus let $\pi_1 = \bigl( \begin{smallmatrix} 3 & 4 & 0 & 0 & 1 & 2   \\
3 & 5 & 2 & 1 & 4 & 5  \end{smallmatrix}\bigr)$ and $\pi_2 = \bigl(
\begin{smallmatrix} 2 & 3 & 4 & 0 & 0 & 1   \\  2 & 4 & 5 & 1 & 3 & 5
\end{smallmatrix}\bigr)$ be two permutations. The corresponding
surfaces belong to $\Q(-1,9)$. The surfaces $\Su_1=\Su(\pi_1,\lambda)$
and $\Su_2=\Su(\pi_2,\lambda)$ possess a (vertical) multiplicity one
simple cylinder. One can checks that $\hat \pi_i$ is irreducible for 
$i=1,2$ thus, by Theorem~\ref{theo:fundamental:2}, each of the
surfaces $\Su_1$ and $\Su_2$ possesses a multiplicity one simple
cylinder. Thanks to a direct computation this cylinder has angle $\pi$
for $i=1$ and $2\pi$ for $i=2$. In other terms this proves $\Su_1 \in
\Q(-1,5) \oplus 1$ and $\Su_2 \in \Q(-1,5) \oplus 2$. \medskip 

\noindent Let us show that $\Q(-1,5) \oplus 4 = \Q(-1,5) \oplus 1$. One has 
$\pi_1 \sim \pi_3 = \bigl(\begin{smallmatrix} 1 & 2 & 3 & 4 & 0 & 0   \\
1 & 4 & 5 & 3 & 5 & 2  \end{smallmatrix}\bigr)$. The permutation $\hat \pi_3$ 
is irreducible and the angle of the cylinder on $\Su(\pi_3,\lambda)$
is $4\pi$. Then surface $\Su(\pi_3,\lambda)$ belongs to $\Q(-1,5)
\oplus 4$ and $\Q(-1,5) \oplus 1$. \medskip

\noindent It remains to end the proof to find a multiplicity one
saddle connection on each surface $\Su_1,\ \Su_2$. This is done with
the following remark. Let us consider the permutations
$$
\pi_1 = \bigl(
\begin{smallmatrix} 3 & 4 & 0 & 0 & 1 & 2   \\
3 & 5 & 2 & 1 & 4 & 5  \end{smallmatrix}\bigr) \sim \pi'_1 =
\bigl(
\begin{smallmatrix} 0 & 1 & 2 & 3 & 4 & 0   \\
1 & 4 & 5 & 3 & 5 & 2  \end{smallmatrix}\bigr) \qquad \textrm{and} \qquad
\pi_2 = \bigl(
\begin{smallmatrix} 2 & 3 & 4 & 0 & 0 & 1   \\
 2 & 4 & 5 & 1 & 3 & 5 \end{smallmatrix}\bigr) \sim \pi'_2 =
\bigl(
\begin{smallmatrix} 0 & 1 & 2 & 3 & 4 & 0   \\
4 & 5 & 1 & 3 & 5& 2 \end{smallmatrix}\bigr).
$$
For $i=1,2$ the vertical separatrix loop $\gamma(\pi'_i)$ is a saddle
connection and $\pi'_i$ is irreducible. Therefore
Proposition~\ref{prop:fundamental:1} applies. Hence first point of the 
proposition is proven.  \bigskip

\noindent {\bf Proof of the second point}

One has to prove that each component of $\Q(-1,13)$ has a
half-translation surface with a multiplicity one saddle
connection. Thanks to Proposition~\ref{prop:contrac:ik:1} we are 
reduced to consider components of the form $\mathcal C \oplus s$ with
$\mathcal C \subset \Q(-1,9)$. The previous proof shows that each
component of $\Q(-1,9)$ 
(different form $\Q^{irr}(-1,9)$) possesses a surface with a 
multiplicity one saddle connection. Finally one has to consider
components of the form
$$
\Q^{irr}(-1,9) \oplus s \qquad \textrm{with } s=1,\dots,6
$$

\noindent Let us recall that $\Q^{irr}(-1,9)= \Q(-1,5) \oplus 3$.
Moreover any component $\Q(-1,5) \oplus s$, with $s=1,2,4,5$, possesses
a flat surface with a multiplicity one saddle connection. Using
properties of the map $\oplus$ (see
Proposition~\ref{prop:properties:oplus}), this yields to 
$\Q(-1,9)^{irr} \oplus s = \Q(-1,5) \oplus s \oplus 3$ which 
proves the proposition for $s=1,2,4,5$. \medskip

\noindent The case $s=6$ is reduced to the case $s=3$.
$$
\Q(-1,9)^{irr} \oplus 6 = \Q(-1,5) \oplus 6 \oplus 3 = \Q(-1,9)^{irr} \oplus 3.
$$
We finish the proof of the second statement using the second property
of the map $\oplus$.
$$
\Q(-1,5) \oplus 3 \oplus 3 = \Q(-1,5) \oplus 1 \oplus 5
$$

\bigskip

\noindent {\bf Proof of the third point}

\noindent The proof is obvious.

\end{proof}

\noindent We finish this section with an independent proof of a
theorem of Masur and Smillie~\cite{Masur:Smillie:2}.
The original proof uses algebraic geometry. Here we only use
combinatorics of generalized permutations.

\begin{NoNumberTheorem}[Masur,~Smillie]
The following strata
$$
 \Q(\emptyset),\ \Q(1,-1)\ (\text{in genus } g=1)
 \quad \text{and}\quad
 \Q(4),\ \Q(1,3)\ (\text{in genus } g=2)
$$
are empty.
\end{NoNumberTheorem}

\begin{proof}[Proof of the theorem]
Let us assume that the stratum $\Q(4)$ is non-empty. Thus there
exists a genus two half-translation surface with a single zero. Applying  
``breaking up a singularity'' to this zero into two zeroes of
order tow we get a point inside the stratum $\Q(2,2)$
(see~\cite{Lanneau} or~\cite{Masur:Zorich},~\cite{Eskin:Masur:Zorich}). 
By construction, this new surface has a multiplicity one saddle connection 
so it belongs to a non-hyperelliptic component of $\Q(2,2)$. 
Now we will prove that this stratum is connected and equals to its
hyperelliptic component which leads to a contradiction. \medskip

\noindent Let us consider a point $\Su(\pi,\lambda) \in \Q(2,2)$. As
usual, a direction computation shows that one can put (up to cyclic
order) $\pi$ into one of the two following forms $\bigl(
\begin{smallmatrix}1 & 2 & 1 & 3    \\  4 & 3 & 4 & 2
\end{smallmatrix}\bigr)$ or $\bigl(\begin{smallmatrix} 1 & 2 & 1 & 2   \\
3 & 4 & 3 & 4  \end{smallmatrix}\bigr)$. According to 
Lemma~\ref{lm:hyper:pi:n-m}, each of these two permutations gives rise to 
surfaces into the component $\Q^{hyp}(2,2)$. 
Therefore $\Q(2,2)=\Q^{hyp}(2,2)$ is connected and hence the 
stratum $\Q(4) = \emptyset$. \medskip

\noindent Using the same approach, we prove that $\Q(1,-1)$ and
$\Q(1,3)$ are empty (one considers respectively the connected hyperelliptic
strata $\Q(-1,-1,2)$ and $\Q(1,1,2)$).\medskip

\noindent The stratum $\Q(\emptyset)$ is empty because the quotient of a quadratic 
differential $\Dq$ by $\D z^2$ is a holomorphic function on the torus thus 
constant. Therefore $\Dq = \omega^2$ which is a contradiction.
\end{proof}

\section{Deformations of hyperelliptic and irreducible components}

\begin{Proposition}
\label{prop:adja:two:singu}
Let $(\Su,\Dq) \in \Q^{hyp}(4(g-k)-6,4k+2)$ be a point with 
$0 \geq k \geq g-2$ and $g \geq 2$. 
Let $(k_1,k_2)$ be a positive partition of $4k+2$. 
Then there exists a continuous path $\rho : [0,1]
\longrightarrow \Q_g$ of the interval $[0,1]$ into the {\it whole} moduli
space $\Q_g$ such that:
\begin{itemize}

\item $\rho(0) = (\Su,\Dq)$

\item $\rho(t) \in \Q(4(g-k)-6,k_1,k_2)$ $\qquad \forall \ 0 < t < 1$.

\item $\rho(1) \in \Q(4(g-k)-6,4k+2) \ \backslash \ \Q^{hyp}(4(g-k)-6,4k+2)$.

\end{itemize}
\end{Proposition}

\begin{proof}
We remark that it is sufficient to prove the proposition for a particular
point of the component $\Q^{hyp}(4(g-k)-6,4k+2)$. We first claim:

\begin{Claim}
Let us fix $r = 2k+1$ and $l = 2(g-k)-3$.
Let $a$ be any integer with $2 \leq a \leq r+1$. Let us consider
the generalized permutation
$$
\Pi_1(r,l,a) =
\left( \begin{array}{cccccccccccc}
0_1 & 0_3 & 1 & \dots & r & 0_1 & r+1 & \dots & r+l \\
r+l & \dots & r+1 & 0_2 & r & \dots & a & 0_3 & a-1 & \dots & 1 & 0_2
\end{array} \right).
$$
Then for any admissible vector $\lambda$ one has 
$\Su(\Pi_1(r,l,a),\lambda) \in \Q(\  a-2,\ 4k+4-a,\ 4(g-k)-3 \
)$. Moreover, the horizontal saddle connection labeled $0_3$ has
multiplicity one. The resulting surface $\tilde S$ obtained by
shrinking this saddle connection to a point is
$\Su(\Pi_1(r,l),\lambda)$ which belongs to the hyperelliptic
component.
\end{Claim}

\begin{proof}[Proof of the claim]
A straightforward calculation of the angle of the conical
singularities located at end-points of the intervals.
\end{proof}

The proposition follows from above claim taking $a=k_1+2$.
\end{proof}

One can easily have similar results on deformations of hyperelliptic
components of other strata. Here we present a similar result
concerning irreducible component. The proof is just based on the
deformations of an adequate generalized permutation.

\begin{Proposition}
\label{prop:adja:irre}
Let $(\Su,\Dq) \in \Q^{irr}(-1,9)$ be a point. Let $(k_1,k_2)$ be any pairs in the list
$\{ (-1,10),(1,8),(2,7),(3,6),(4,5) \}$. Then there exists a continuous path
$\rho : [0,1] \longrightarrow \Q_3$ of the interval
$[0,1]$ into the {\it whole} moduli space $\Q_3$ such that:
\begin{itemize}

\item $\rho(0) = (\Su,\Dq)$

\item $\rho(t) \in \Q(-1,k_1,k_2)$ $\qquad \forall \ 0 < t < 1$.

\item $\rho(1) \in \Q(-1,9) \ \backslash \ \Q^{irr}(-1,9)$.

\end{itemize}
\end{Proposition}

\begin{proof}
The proof is similar to the previous one, deforming the permutation
$$
\left( \begin{array}{cccccc}
0 & 1 & 2 & 3 & 4 & 0\\
4& 3& 2& 5&1&5
\end{array} \right).
$$
representing the irreducible component of $\Q(-1,9)$.
\end{proof}

\begin{Proposition}
\label{prop:adja:12}
Let $(\Su,\Dq) \in \Q^{irr,II}(12)$ be a point. Let $(k_1,k_2)$ be any
pair in the list $\{ (-1,13),(1,11),(2,10),(3,9),(4,8),(5,7),(6,6)
\}$. \\
Then there exists a continuous path $\rho : [0,1] \longrightarrow
\Q_4$ of the interval $[0,1]$ into the {\it whole} moduli space $\Q_4$
such that:
\begin{itemize}

\item $\rho(0) = (\Su,\Dq)$

\item $\rho(t) \in \Q(k_1,k_2)$ $\qquad \forall \ 0 < t < 1$.

\item $\rho(1) \in \Q^{irr,I}(12)$

\end{itemize}
\end{Proposition}

\begin{proof}
The proof parallels the one of Proposition~\ref{prop:adja:two:singu}
deforming the permutation $\left( \begin{smallmatrix} 1 & 2 & 3 & 4 &
  3 & 5 & 6\\ 1 & 5& 7& 4& 2&6 & 7 \end{smallmatrix} \right)$.
\end{proof}


\end{document}